\documentclass[a4paper,reqno]{amsart}
\usepackage[english]{babel}
\usepackage{amsmath,amssymb,amsfonts,amsthm}

\usepackage{scrtime}

\usepackage{enumitem}

\usepackage{hyperref}
\hypersetup{colorlinks = true, urlcolor = magenta, linkcolor = blue, citecolor = red}

\setlist[enumerate]{leftmargin=*}
\setlist[itemize]{leftmargin=*}

\newtheorem{thm}{Theorem}[section]
\newtheorem{prp}[thm]{Proposition}
\newtheorem{lem}[thm]{Lemma}
\newtheorem{cor}[thm]{Corollary}

\newtheorem{lthCZ}{Theorem}

\newtheorem{lthHardy}{Theorem}

\newtheorem{lthTrees}{Theorem}

\theoremstyle{definition}
\newtheorem{dfn}[thm]{Definition}
\newtheorem{rem}[thm]{Remark}

\numberwithin{equation}{section}

\def\Xint#1{\mathchoice
   {\XXint\displaystyle\textstyle{#1}}%
   {\XXint\textstyle\scriptstyle{#1}}%
   {\XXint\scriptstyle\scriptscriptstyle{#1}}%
   {\XXint\scriptscriptstyle\scriptscriptstyle{#1}}%
   \!\int}
\def\XXint#1#2#3{{\setbox0=\hbox{$#1{#2#3}{\int}$}
     \vcenter{\hbox{$#2#3$}}\kern-.5\wd0}}

\def\dashint{\Xint-}

\DeclareFontFamily{U}{matha}{\hyphenchar\font45}
\DeclareFontShape{U}{matha}{m}{n}{
      <5> <6> <7> <8> <9> <10> gen * matha
      <10.95> matha10 <12> <14.4> <17.28> <20.74> <24.88> matha12
      }{}
\DeclareSymbolFont{matha}{U}{matha}{m}{n}
\DeclareFontSubstitution{U}{matha}{m}{n}

\DeclareFontFamily{U}{mathx}{\hyphenchar\font45}
\DeclareFontShape{U}{mathx}{m}{n}{
      <5> <6> <7> <8> <9> <10>
      <10.95> <12> <14.4> <17.28> <20.74> <24.88>
      mathx10
      }{}
\DeclareSymbolFont{mathx}{U}{mathx}{m}{n}
\DeclareFontSubstitution{U}{mathx}{m}{n}

\DeclareMathDelimiter{\vvvert}{0}{matha}{"7E}{mathx}{"17}

\newcommand{\rmI}{\mathrm{I}}
\newcommand{\rmII}{\mathrm{II}}

\newcommand{\loc}{\mathrm{loc}}
\newcommand{\sloc}{\mathrm{sloc}}
\newcommand{\fin}{\mathrm{fin}}
\newcommand{\comp}{\mathrm{comp}}
\newcommand{\ord}{\mathrm{ord}}

\newcommand{\Adm}{\mathcal{A}}
\newcommand{\cE}{\mathcal{E}}
\newcommand{\cM}{\mathcal{M}}
\newcommand{\cQ}{\mathcal{Q}}
\newcommand{\cS}{\mathcal{S}}

\newcommand{\Fth}{\widehat}

\newcommand{\Ch}{\mathfrak{C}}

\newcommand{\defeq}{\mathrel{:=}}
\newcommand{\eqdef}{\mathrel{=:}}
\newcommand{\tc}{\,:\,}

\newcommand{\NN}{\mathbb{N}}
\newcommand{\ZZ}{\mathbb{Z}}
\newcommand{\RR}{\mathbb{R}}
\newcommand{\CC}{\mathbb{C}}

\newcommand{\Npos}{\mathbb{N}_+}

\newcommand{\BMO}{\mathrm{BMO}}

\newcommand{\lH}{\mathfrak{h}}

\newcommand{\chr}{\mathbf{1}}

\DeclareMathOperator{\Span}{span}

\DeclareMathOperator*{\esssup}{ess\,sup}

\newcommand{\Tree}{\mathfrak{T}}

\newcommand{\scc}{\mathfrak{s}}
\newcommand{\prd}{\mathfrak{p}}
\newcommand{\lev}{\ell}
\newcommand{\dist}{d}
\DeclareMathOperator{\diam}{diam}
\newcommand{\Rz}{\mathcal{R}}
\newcommand{\opL}{\mathcal{L}}
\newcommand{\czz}{c_{00}}
\newcommand{\opSh}{\Sigma}
\DeclareMathOperator{\supp}{supp}
\newcommand{\modE}{\mathcal{E}}

\newcommand{\sobolev}[2]{L^{#1}_{#2}}

\title[Nondoubling Calder\'on--Zygmund theory]{Nondoubling Calder\'on--Zygmund theory\\across geometries}

\author{Jos\'e M. Conde Alonso}
\address[J.M. Conde Alonso]{Departamento de Matem\'aticas \\ Universidad Aut\'onoma de Madrid \\ Ciudad Universitaria de Cantoblanco \\ 
28049 Madrid \\ Spain}
\email{jose.conde@uam.es}

\author{Alessio Martini}
\address[A. Martini]{Dipartimento di Scienze Matematiche ``G.L. Lagrange'' \\ Politecnico di Torino \\ Corso Duca degli Abruzzi 24 \\ 10129 Torino \\ Italy}
\email{alessio.martini@polito.it}

\author{Maria Vallarino}
\address[M. Vallarino]{Dipartimento di Scienze Matematiche ``G.L. Lagrange'' \\ Politecnico di Torino \\ Corso Duca degli Abruzzi 24 \\ 10129 Torino \\ Italy}
\email{maria.vallarino@polito.it}

\thanks{The first-named author was partially supported by the 2025 Leonardo Grant for Scientific Research and Cultural Creation \texttt{LEO25-1-19094-CBB-MAT-169} (BBVA Foundation). The BBVA Foundation accepts no responsibility for the opinions, statements and contents included in the project and/or the results thereof, which are entirely the responsibility of the authors. The second- and third-named authors are members of Gruppo Nazionale per l'Analisi Mate\-matica, la Probabilit\`a e le loro Applicazioni (GNAMPA) of Istituto Nazionale di Alta Matematica (INdAM)}

\begin{document}

\begin{abstract}
We develop a Calderón--Zygmund theory without the doubling property, a metric or a martingale structure. Our measure spaces carry an exhausting family of admissible sets with recursive child partitions. A stopping-time construction gives a decomposition of an integrable function, with a bounded good part and bad terms that cancel on the parents of the stopping sets. We introduce the corresponding atomic Hardy and BMO spaces, establish their duality and interpolation properties, and prove weak-type $(1,1)$ and Hardy-space bounds for singular integral operators under a H\"ormander condition and a size estimate between each stopping set and its parent. We apply this framework to arbitrary flow trees, whose associated flow measures need not be either locally or globally doubling. As a result, we obtain endpoint bounds for Riesz transforms and spectral multipliers satisfying a scale-invariant Sobolev condition of order greater than $3/2$ in this setting.
\end{abstract}

\maketitle

\section{Introduction}

The study of the boundedness properties of integral operators near the $L^1$-endpoint is intimately tied to the decomposition of functions into pieces with additional properties: cancellation, localization, higher integrability, and so on. For example, the Calder\'on--Zygmund decomposition splits an integrable function at a prescribed height $\lambda$ into a bounded good part and a sum of localized mean-zero bad parts. The construction is usually carried out over the dyadic filtration of $\RR^n$, which is doubling. Indeed, given a dyadic cube $Q$, both $2Q$ (the centered dilate) and $\Fth Q$ (the dyadic parent) have measure comparable to that of $Q$. The doubling property enters at a decisive point: a maximal stopping set has measure comparable to that of its parent, and the average there of an integrable function cannot be much larger than the stopping height. This is what places the good part in $L^\infty$. 

In a nondoubling setting the measure of a dyadic child may be arbitrarily small compared with that of its parent. The average on a stopping child may then be arbitrarily large even when the average on the parent is at most $\lambda$, and the classical replacement of the function by its average on the child no longer produces a bounded good part. This has been overcome in several scenarios: first, on $\RR^n$ a Besicovitch-type covering can be used to write a non-dyadic Calderón--Zygmund decomposition suitable for studying singular integrals when the underlying measure has polynomial growth, see \cite{To01}. Second, on a general martingale space one can modify the classical Gundy decomposition to write down a three-term decomposition \cite{LSMP}, which can be further refined to drop one term \cite{CP19,CPW}. All these approaches exploit the same idea: since one cannot force the averages of the function $f$ to be both large and small in the stopping sets, one has to impose cancellation on the parent of the stopping set instead. This creates a second difficulty for the study of singular integrals: cancellation controls the kernel far from the support of a bad function, but it does not control the intermediate region between the stopping set and the set on which the cancellation takes place. On the other hand, even if the stopping sets are disjoint --- usually by construction --- their parents need not be. This destroys part of the orthogonality in the decomposition, yielding worse integrability estimates for the good part than in the classical case.  

We aim to build a useful Calderón--Zygmund decomposition with minimal assumptions, distilling what is necessary and getting rid of everything else. At the same time, we shall show below the applicability of our construction by finding a scenario where no prior approaches work. We draw inspiration in the locally doubling case from \cite{HS,V}. Our starting point is that the decomposition that we are looking for only needs a recursive stopping structure. Let $(X,\cE,m)$ be a measure space and let $\Adm\subseteq\cE$ be a family of sets of finite positive measure. Each $Q\in\Adm$ has an at most countable family of children $\Ch(Q)\subseteq\Adm$ which forms a partition of $Q$; the descendants of $Q$ generate the measurable subsets of $Q$ up to null sets; and $X$ admits an increasing exhaustion by elements of $\Adm$. We call $(X,\cE,m,\Adm)$ a \emph{Calderón--Zygmund space} (and we abbreviate it as CZ space in the body of the paper). We impose no metric structure, no bound on the number of children, and no uniform estimate for $m(\Fth Q)/m(Q)$ when $Q\in\Ch(\Fth Q)$.

\subsection*{A Calder\'on--Zygmund decomposition with a bounded good part and weak-type estimates} Fix a Calderón--Zygmund space. If a stopping child $Q$ has parent $\Fth Q$, our bad functions are not required to be supported on $Q$. Instead, for a function $a$ supported on $Q$, we take
\[
b=a-\left(\int_X a\,dm\right)\frac{\chr_{\Fth Q}}{m(\Fth Q)}
\]
to be the mean-zero atom associated with $a$; here $\chr_S$ denotes the characteristic function of a set $S$. The function $b$ is supported on $\Fth Q$ and has integral zero. Since the correcting constant is spread over the parent, its size is governed by the average on $\Fth Q$, even when $m(Q)/m(\Fth Q)$ is arbitrarily small. A first stopping-time construction gives a good term in every $L^p$ space with finite $p$. Iterating that construction at decreasing heights yields the following stronger conclusion.

\begin{lthCZ}\label{thm:mainCZdec}
Let $(X,\cE,m,\Adm)$ be a Calderón--Zygmund space, let $f\in L^1(X)$, and let
\[
\lambda>\frac{\|f\|_1}{m(X)},
\]
with the usual convention when $m(X)=\infty$. There are an at most countable set $J$ and, for each $j \in J$, pairs $Q_j,\Fth Q_j\in\Adm$ with $Q_j\in\Ch(\Fth Q_j)$, functions $a_j\in L^1(X)$ vanishing off $Q_j$, and a decomposition
\[
f=g+\sum_{j\in J}b_j,
\qquad
b_j=a_j-\left(\int_X a_j\,dm\right)
       \frac{\chr_{\Fth Q_j}}{m(\Fth Q_j)},
\]
such that
\[
\|g\|_\infty\lesssim\lambda,
\qquad
\|g\|_1\leq 4\|f\|_1,
\qquad
\sum_{j\in J}m(Q_j)\leq \frac{2\|f\|_1}{\lambda},
\qquad
\sum_{j\in J}\|a_j\|_1\leq 2\|f\|_1.
\]
In particular, $b_j$ vanishes off $\Fth Q_j$ and $\int_X b_j \,dm = 0$ for all $j \in J$, and
\[
\sum_{j\in J}\|b_j\|_1\leq4\|f\|_1.
\]
\end{lthCZ}

Here and in what follows, for two nonnegative quantities $A$ and $B$, we write $A \lesssim B$ to indicate that the estimate $A \leq C B$ holds for some positive finite constant $C$. Moreover, $A \simeq B$ denotes the conjunction of $A \lesssim B$ and $B \lesssim A$. Subscripted variants, such as $\lesssim_p$ or $\simeq_p$, indicate that the implicit constants may depend on the parameter $p$.

The $L^\infty$ estimate is the main novelty of the theorem. For general atomic filtrations, the decomposition from \cite{CPW} gives control of the good part in martingale BMO. Earlier ones \cite{LSMP,CP19} only provided $L^p$-estimates for finite $p$. Theorem~\ref{thm:mainCZdec} upgrades this to the classical conclusion $\|g\|_\infty\lesssim\lambda$ without imposing a condition equivalent to regularity of the filtration. Besides, our construction does not rely on an underlying martingale structure. In addition to its role in the endpoint theory developed here, we believe that this stronger form of Calderón--Zygmund decomposition may have applications in multilinear harmonic analysis; we do not pursue that direction in the present paper.

\subsection*{An abstract weak-type criterion}

To pass from the decomposition to operators, suppose that every $Q\in\Adm$ is assigned a \emph{doubling enlargement} $Q^*$ such that
\begin{equation}\label{eq:extQ}
Q\subseteq Q^*,\qquad m(Q^*)\leq\kappa m(Q)
\end{equation}
for some constant $\kappa \in [1,\infty)$.
The enlargement is an additional geometric datum; it need not be an admissible set or belong to the recursive family. For an operator $T$ with kernel $K_T$, we use an integral H\"ormander condition away from $Q^*$ and a size condition in the region between the enlargements of a child and its parent.

\begin{lthCZ}\label{thm:mainCZwt}
Let $(X,\cE,m,\Adm)$ be a Calderón--Zygmund space with enlargements as in \eqref{eq:extQ}, and let $T$ be bounded on $L^q(X)$ for some $q\in(1,\infty)$. Assume that $T$ has an off-diagonal integral kernel $K_T$, in the sense that
\begin{equation}\label{eq:CZop_offdiag}
Tf(x)=\int_X K_T(x,y)f(y)\,dm(y)
\qquad\text{for a.e. }x\in X\setminus Q^*
\end{equation}
whenever $f \in L^q(X)$ vanishes off $Q \in \Adm$. Assume also that the kernel $K_T$ satisfies
\begin{align}
\label{eq:CZop_horm}
C_H(T)&\defeq
\sup_{Q\in\Adm}\ \esssup_{y,z\in Q}
\int_{X\setminus Q^*}|K_T(x,y)-K_T(x,z)|\,dm(x)<\infty,\\
\label{eq:CZop_size}
C_S(T)&\defeq
\sup_{\Fth Q\in\Adm}\ \sup_{Q\in\Ch(\Fth Q)}\ \esssup_{y\in Q}
\int_{\Fth Q^*\setminus Q^*}|K_T(x,y)|\,dm(x)<\infty.
\end{align}
Then $T$ is of weak type $(1,1)$ and is bounded on $L^p(X)$ for every $p\in(1,q]$. The bounds depend only on $q$, $p$, $\kappa$, $\|T\|_{q\to q}$, $C_H(T)$, and $C_S(T)$.
\end{lthCZ}

The distinction between $Q^*$ and $\Fth Q$ is essential. The parent $\Fth Q$ belongs to the recursive structure and is the set on which the bad term produced by Theorem~\ref{thm:mainCZdec} is supported and has cancellation. Its measure may be arbitrarily large compared with $m(Q)$. By contrast, $Q^*$ is a geometric enlargement whose measure is controlled by $m(Q)$, and it is therefore $Q^*$, rather than $\Fth Q$, that can be placed in the exceptional set of a weak-type argument. Cancellation controls the contribution of the bad term only outside $\Fth Q^*$. This leaves the intermediate region $\Fth Q^*\setminus Q^*$, where cancellation alone gives no estimate; the quantity $C_S(T)$ is designed precisely to control that region. In regular or doubling settings these two roles are often hidden by measure comparability (the two sets $2Q$ and $\Fth Q$ mentioned above; see also Remark \ref{rem:size_smoothness} below). Here their separation is the main point of the theorem.

\subsection*{The abstract \texorpdfstring{$H^1$}{H1}--\texorpdfstring{$\BMO$}{BMO} theory}

The parent-supported cancellation in Theorem~\ref{thm:mainCZdec} also determines the appropriate Hardy space. Inspired by \cite{CP}, given $r\in(1,\infty]$, an $r$-atom is a function
\begin{equation}\label{eq:atoms}
b=a-\left(\int_Xa\,dm\right)\frac{\chr_{\Fth Q}}{m(\Fth Q)},
\end{equation}
where $Q\in\Ch(\Fth Q)$, $a$ vanishes off $Q$, and
\[
\|a\|_r\leq m(Q)^{-1/r'}.
\]
The space $H^{1,r}(X)$ consists of the $\ell^1$ sums of such atoms. When $m(X)<\infty$, the local space $\lH^{1,r}(X)$ is obtained by adjoining the normalized constant atom $\chr_X/m(X)$.

The corresponding BMO seminorms are
\begin{equation}\label{eq:BMOseminorm}
\|h\|_{\BMO_q}
\defeq \sup_{\Fth Q\in\Adm}\sup_{Q\in\Ch(\Fth Q)}
\left(\dashint_Q\left|h-\dashint_{\Fth Q}h\,dm\right|^q\,dm\right)^{1/q},
\end{equation}
where $q \in [1,\infty)$; here
\[
\dashint_S f \,dm \defeq \frac{1}{m(S)} \int_S f \,dm
\]
denotes the mean value of $f$ on a set $S$. Thus oscillation is measured on the child relative to the average on its parent.
We shall denote by $\BMO_q(X)$ the space of functions $h$ for which $\|h\|_{\BMO_q}$ is finite; 
by quotienting out constant functions we obtain the normed space $\BMO_q(X)/\CC$. When $m(X)<\infty$, we equip $\BMO_q(X)$ itself with the norm $\vvvert h \vvvert_{\BMO_q} \defeq \max\left\{ \|h\|_{\BMO_q}, \left| \dashint_X h \,dm \right| \right\}$.

The next results collect the structural, interpolation, and operator endpoint consequences of these definitions. As usual, $p' = p/(p-1)$ denotes the conjugate exponent to $p \in [1,\infty]$.

\begin{lthHardy}\label{thm:mainHardyequiv}
Let $(X,\cE,m,\Adm)$ be a Calderón--Zygmund space.
Then the spaces $H^{1,r}(X)$, $r\in(1,\infty]$, coincide with equivalent norms. Analogous statements hold for the local spaces $\lH^{1,r}(X)$ when $m(X)<\infty$, and also for the spaces $\BMO_{r'}(X)$. Moreover, with respect to the natural duality pairings,
\[
H^{1,r}(X)' = \BMO_{r'}(X)/\CC
\]
isometrically; if $m(X) < \infty$, then also
\[
\lH^{1,r}(X)' = \BMO_{r'}(X)
\]
isometrically.
\end{lthHardy}

When we do not need to specify which of the equivalent norms we are using, we drop the additional indices in the notation for Hardy and BMO spaces, and simply write $H^1(X)$, $\lH^1(X)$ and $\BMO(X)$.

We write $[\cdot,\cdot]_{\theta,p}$ for the usual real interpolation functors.

\begin{lthHardy}\label{thm:mainHardyinterp}
Let $(X,\cE,m,\Adm)$ be a Calderón--Zygmund space.
If $m(X)=\infty$ 
and $1<r<p\leq\infty$, then
\[
L^r(X) = [H^1(X),L^p(X)]_{p'/r',r} \ ;
\]
moreover, 
for every $r\in(1,\infty)$,
\[
L^r(X) = [H^1(X),\BMO(X)/\CC]_{1/r',r} \ .
\]
When $m(X)<\infty$, the corresponding identities hold with $\lH^1(X)$ in place of $H^1(X)$ and with $\BMO(X)$ in place of $\BMO(X)/\CC$.
\end{lthHardy}

\begin{lthHardy}\label{thm:mainHardybdd}
Let $(X,\cE,m,\Adm)$ be a Calderón--Zygmund space.
If $T$ satisfies the assumptions of Theorem~\ref{thm:mainCZwt}, then
\[
T : H^1(X) \to L^1(X),
\qquad
T^* : L^\infty(X) \to \BMO(X)/\CC
\]
are bounded, with bounds depending only on $q$, $\kappa$, $\|T\|_{q\to q}$, $C_H(T)$, and $C_S(T)$.
In the finite-measure case the boundedness holds for the maps $T:\lH^1(X)\to L^1(X)$ and $T^*:L^\infty(X)\to\BMO(X)$.
\end{lthHardy}

The proof of Theorem \ref{thm:mainHardyequiv} goes through a John--Nirenberg inequality, which implies the equality of the spaces $\BMO_r$ for $r \in [1,\infty)$. A functional-analytic argument, based on $L^q$--$L^{q'}$ duality, then gives the duality isomorphisms $(H^{1,q})' = \BMO_{q'}/\CC$ for all $q \in (1,\infty]$, from which the equality of the spaces $H^{1,q}$ follows. In particular, we do not rely on the rather technical constructions of atomic decompositions used elsewhere (see, e.g., \cite[Theorem A]{CW}) to prove the $q$-independence of $H^{1,q}$. Moreover, unlike other duality-based approaches (see, e.g., \cite[Section 6 and Remark 6.2]{CMM}), here we manage to include the endpoint $q=\infty$ (on this, see also the comments before Proposition \ref{prp:sH1duality} below).

Another novelty in our result is that our duality isomorphisms are \emph{isometric}, something that other approaches in the literature do not seem to achieve; for this, it is crucial that atoms and BMO seminorms are defined as in \eqref{eq:atoms}-\eqref{eq:BMOseminorm}, and indeed the result arguably demonstrates the naturality of our definitions. On the other hand, when the underlying Calder\'on--Zygmund space is \emph{regular}, i.e., the measure of an admissible set is comparable to that of any of its children, one can replace those definitions with more ``traditional'' ones and obtain the same Hardy and BMO spaces, up to equivalent norms (see Remarks \ref{rem:regularBMO} and \ref{rem:regularHardy}).

The interpolation statements in Theorem \ref{thm:mainHardyinterp} are obtained directly from the decomposition in Theorem~\ref{thm:mainCZdec}; in particular, the argument does not rely on importing martingale interpolation as a black box. 
Further, the operator estimates in Theorems \ref{thm:mainCZwt} and \ref{thm:mainHardybdd} show that the same separation between the parent $\Fth Q$ and the controlled enlargement $Q^*$ governs both weak-type and Hardy-space endpoints.

\subsection*{Application scenarios and endpoint estimates on flow trees}

Our abstract theory applies to a variety of geometric settings, thus providing a unified framework for a number of results about singular integrals and Hardy--BMO spaces available in the literature. Besides doubling metric measure spaces \cite{CW71,CW}, this includes globally nondoubling settings such as $ax+b$ groups \cite{HS,V} and hypergroups \cite{MPl}, as well as rank-one extensions of Carnot groups \cite{MOV,DMLMV25} and Damek--Ricci spaces \cite{V07,LM26,BeLi}. Actually, all of the mentioned settings can be given the structure of a regular Calder\'on--Zygmund space, where the measure of an admissible set $Q$ is comparable to that of its parent $\Fth Q$. Crucially, however, such measure comparability is not required by our abstract theory, whose applicability thus extends to fully nondoubling situations. Indeed, our theory can be applied to recover the weak-type $(1,1)$ estimates for Calderón--Zygmund operators on $\mathbb{R}^n$ with a measure of polynomial growth, which is the the main result in \cite{To01}. This can be done via \cite{CP19}, which constructs a filtration that gives rise to a (possibly non-regular) Calder\'on--Zygmund space structure. More generally, Theorem \ref{thm:mainCZdec} is an improvement over the martingale-based Calderón--Zygmund decompositions from  \cite{LSMP,CP19,CCP2022,CPW}.

Beyond these known examples, our methods allow us to address problems which could not be tackled via existing approaches. As an example of these new applications, we shall show how our theory can be used to prove endpoint estimates for natural singular integrals on flow trees. A flow tree $(\Tree,\mu)$ is an oriented, locally finite tree equipped with a positive measure satisfying the flow identity
\[
\mu(x)=\sum_{y\in\scc(x)}\mu(y),
\]
where $\scc(x)$ is the set of successors of $x$. Such measures retain a strong conservation law along the orientation, but may be highly nondoubling and need not satisfy any local doubling estimate. Earlier work developed harmonic analysis for nondoubling flow measures \cite{LSTV}. Results on homogeneous trees and subsequent transference arguments yielded $L^p$ bounds, for $1<p<\infty$, for the Riesz transform and spectral multipliers associated with the flow Laplacian \cite{MSV,MSTV}. The weak-type and Hardy-space endpoints in \cite{MSTV}, however, required the local doubling property. On a flow tree, the admissible sets are singletons and truncated descendant trapezoids. Their recursive partitions alternate between vertical subdivision and subdivision according to the successor relation. Although a child may carry an arbitrarily small fraction of the mass of its parent, every admissible set $Q$ has an enlargement $Q^*$ satisfying $\mu(Q^*)\leq3\mu(Q)$. Moreover, the abstract kernel assumptions follow from two estimates that refer only to the metric and the level function (see Proposition \ref{prp:treeCZop} below).

The required kernel bounds can be verified for the first-order Riesz transform and for spectral multipliers of the flow Laplacian. This gives the following endpoint results without any local doubling assumption.

\begin{lthTrees}\label{thm:maintrees}
Let $(\Tree,\mu)$ be an arbitrary flow tree and let $\nabla$ and $\opL= \frac{1}{2} \nabla^* \nabla$ be the flow gradient and flow Laplacian thereon.
\begin{enumerate}[label=\textnormal{(\alph*)}]
\item The Riesz transform
\[
\Rz=\nabla\opL^{-1/2}
\]
is of weak type $(1,1)$ and is bounded on $L^p(\Tree)$ for every $p\in(1,\infty)$. At the Hardy--BMO endpoints,
\[
\Rz : H^1(\Tree) \to L^1(\Tree),
\qquad
\Rz^* : L^\infty(\Tree) \to \BMO(\Tree)/\CC
\]
are bounded.
\item Let $s>3/2$, let $\chi\in C_c^\infty(3/4,5/4)$ satisfy $0\leq\chi\leq1$ and $\chi(1)=1$, and let $F:\RR\to\CC$ be Borel measurable. If
\[
\sup_{0<t\leq1} \|F(t\,\cdot)\chi\|_{\sobolev{2}{s}}
+\sup_{0<t\leq1} \|F(2-t\,\cdot)\chi\|_{\sobolev{2}{s}}<\infty,
\]
then $F(\opL)$ is of weak type $(1,1)$ and is bounded on $L^p(\Tree)$ for every $p\in(1,\infty)$.
%
Under the stronger condition
\[
\sup_{0<t\leq2}\|F(t\,\cdot)\chi\|_{\sobolev{2}{s}}<\infty,
\]
the operator $F(\opL)$ also maps $H^1(\Tree)$ boundedly into $L^1(\Tree)$, and $L^\infty(\Tree)$ boundedly into $\BMO(\Tree)/\CC$.
\end{enumerate}
\end{lthTrees}

The $L^p$ estimates in Theorem~\ref{thm:maintrees} were already known on arbitrary flow trees, while the endpoint estimates had previously been established under local doubling assumptions \cite{MSTV}. Thus the new conclusions are the weak-type $(1,1)$ and Hardy--BMO endpoints in the fully nondoubling setting. The two multiplier hypotheses reflect the two levels of kernel control: the condition adapted to both ends of the spectrum $[0,2]$ of $\opL$ is enough for the weak-type bound, whereas the uniform condition also yields the Hardy and BMO bounds. We refer to \cite{MSTV} for a discussion of the optimality of the regularity threshold $3/2$ for spectral multipliers on flow trees. We also point out that the adjoint Riesz transform $\Rz^*$ in general is not $H^1 \to L^1$ bounded \cite{MSTV}, nor is it of weak type $(1,1)$ \cite{Liu}.

\smallskip
The rest of the paper is organized as follows. Section \ref{sec:sec1} develops the Calderón--Zygmund machinery necessary to prove Theorem \ref{thm:mainCZdec}. It also includes the weak-$L^1$ abstract endpoint theory. Section \ref{sec:sec2} introduces BMO and the atomic Hardy spaces, proves duality and interpolation, and establishes the abstract kernel criterion. Finally, Section \ref{sec:sec3} explores applications of the theory. For flow trees, it constructs admissible trapezoids and their enlargements, proves the tree-specific form of the kernel criterion, and applies it to the Riesz transform and spectral multipliers.

\section{Calderón--Zygmund spaces}\label{sec:sec1}

\subsection{Basic properties of CZ spaces and maximal functions} As we stated in the Introduction, we work with minimal assumptions that allow us to construct a CZ decomposition. The following definition repeats what was said in the Introduction plus a couple of measurability technicalities.

In what follows, if $\cQ$ is a collection of sets, then $\bigcup \cQ$ denotes the union of the collection, i.e., $\bigcup \cQ = \bigcup_{Q \in \cQ} Q$.

\begin{dfn}\label{dfn:CZdef}
A \emph{CZ space} is a measure space $(X,\cE,m)$ equipped with a collection $\Adm \subseteq \cE$ of \emph{admissible sets}, with the following properties.
\begin{enumerate}[label=(\alph*)]
\item\label{en:CZdef_meas} Each $Q \in \Adm$ satisfies $0 < m(Q) < \infty$.
\item\label{en:CZdef_exh} There exists an \emph{exhausting sequence} $(Q_n)_{n \in \NN}$ of admissible sets, satisfying $Q_{n} \subseteq Q_{n+1}$ for all $n \in \NN$ and $X = \bigcup_{n \in \NN} Q_n$.
\item\label{en:CZdef_chil} Each $Q \in \Adm$ is associated with an at most countable disjoint collection $\Ch(Q) \subseteq \Adm$ of \emph{children} of $Q$, such that $Q = \bigcup \Ch(Q)$; correspondingly, for any $k \in \NN$, we define the collection $\Ch^k(Q)$ of \emph{$k$-descendants} of $Q$ by setting
\[
\Ch^{0}(Q) \defeq \{Q\}, \quad \Ch^{k+1}(Q) \defeq \bigcup_{R \in \Ch^k(Q)} \Ch(R),
\]
and let $\Ch^*(Q) \defeq \bigcup_{k \in \NN} \Ch^k(Q)$.
\item\label{en:CZdef_sigma} For all $Q \in \Adm$, the $m$-completion of the $\sigma$-algebra on $Q$ generated by $\Ch^*(Q)$ contains $\cE|_Q \defeq \{ E \cap Q \tc E \in \cE \}$.
\end{enumerate}
\end{dfn}

Let $(X,\cE,m,\Adm)$ be a CZ space.
We denote by $L^0(X)$ the space of measurable functions $f : X \to \CC$, identified up to almost everywhere equality. A function $f \in L^0(X)$ is said to be \emph{well-supported} if, for any $\lambda \in (0,m(X))$, there exists $Q \in \Adm$ with $m(Q) > \lambda$ and $\{ f \neq 0 \} \subseteq Q$.

\begin{prp}\label{prp:CZbasic}
Assume that $(X,\cE,m,\Adm)$ is a CZ space.
\begin{enumerate}[label=(\roman*)]
\item\label{en:CZbasic_sigma} $(X,\cE,m)$ is $\sigma$-finite.
\item\label{en:CZbasic_approx} If $(Q_n)_{n \in \NN}$ is an exhausting sequence, then, for any $f \in L^0(X)$, the sequence $(f_n)_{n \in \NN}$ defined by $f_n \defeq f \chr_{Q_n}$ satisfies:
\begin{enumerate}[label=(\alph*)]
\item $|f_n| \leq |f|$ for all $n \in \NN$;
\item $f_n \to f$ pointwise;
\item each $f_n$ is well-supported.
\end{enumerate}
\item\label{en:CZbasic_span} For any $Q_0 \in \Adm$, the space $\Span \{ \chr_Q \tc Q \in \Ch^*(Q_0) \}$ is dense in $L^1(Q_0)$.
\end{enumerate}
\end{prp}
\begin{proof}
Parts \ref{en:CZbasic_sigma} and \ref{en:CZbasic_approx} are obvious.

As for \ref{en:CZbasic_span}, this follows from the definition of CZ space, namely, we can approximate $f \in L^1(Q_0)$ with the sequence $(f_n)_n$ given by
\begin{equation}\label{eq:martingale_approx}
f_n \defeq \sum_{Q \in \Ch^n(Q_0)} \left(\dashint_Q f \,dm\right) \chr_Q.
\end{equation}
As each $\Ch^n(Q_0)$ is at most countable, clearly each $f_n$ is in the $L^1(Q_0)$-closure of $\Span\{ \chr_Q \tc Q \in \Ch^*(Q_0) \}$, so it is enough to show that $f_n \to f$ in $L^1(Q_0)$.

Now, by construction, if $\Sigma_n$ denotes the $\sigma$-algebra on $Q_0$ generated by $\Ch^n(Q_0)$ for any $n \in \NN$, then the sequence $(f_n)_{n \in \NN}$ is a martingale on the probability space $(Q_0,\cE|_{Q_0},m(Q_0)^{-1} m)$ with respect to the filtration $(\Sigma_n)_n$, and this martingale is closed by $f \in L^1(Q_0)$. In particular, the martingale $(f_n)_n$ is uniformly integrable and converges in $L^1$ to a function $\tilde f \in L^1(Q_0)$, which is measurable with respect to the $\sigma$-algebra $\Sigma$ generated by $\Ch^*(Q_0) = \bigcup_n \Ch^n(Q_0)$, and moreover $(f_n)_n$ is closed by $\tilde f$ too \cite[Theorem 9.22]{K}. Consequently, $\tilde f$ is a version of the conditional expectation of $f$ with respect to $\Sigma$, and $\int_{E} f \,dm = \int_E \tilde f \,dm$ for any $E \in \Sigma$; by altering $E$ by negligible sets, we see that the same equality holds more generally for any $E$ in the $m$-completion of $\Sigma$, so (by the definition of CZ space) for any $E \in \cE|_{Q_0}$. As both $f,\tilde f$ are $\cE|_{Q_0}$-measurable, this proves that $f = \tilde f$ a.e., and thus $f_n \to f$ in $L^1(Q_0)$.
\end{proof}

\begin{dfn}\label{dfn:regularCZ}
A CZ space $(X,\cE,m,\Adm)$ is said to be \emph{regular} if 
\[
\rho(X) \defeq \sup_{\Fth Q \in \Adm} \sup_{Q \in \Ch(\Fth Q)} \frac{m(\Fth Q)}{m(Q)} < \infty.
\]
\end{dfn}

\begin{rem}\label{rem:regularCZ}
Assume that the CZ space $(X,\cE,m,\Adm)$ is regular.
Then, for any $\Fth Q \in \Adm$,
\[
m(\Fth Q) = \sum_{Q \in \Ch(\Fth Q)} m(Q) \geq \frac{\# \Ch(\Fth Q)}{\rho(X)} m(\Fth Q),
\]
and therefore
\[
\sup_{\Fth Q \in \Adm} \# \Ch(\Fth Q) \leq \rho(X).
\]
\end{rem}

In what follows, unless otherwise specified, we shall not assume that a CZ space $(X,\cE,m,\Adm)$ is regular.

\begin{dfn}
Let $(X,\cE,m,\Adm)$ be a CZ space.
For any $Q_0 \in \Adm$ and any measurable function $f$ on $X$, we define the associated maximal function by
\[
\cM_{Q_0} f(x) \defeq \begin{cases}
\sup\limits_{Q \in \Ch^*(Q_0) \tc  x \in Q} \dashint_Q |f| \,dm &\text{if } x \in Q_0,\\
0 &\text{otherwise}.
\end{cases}
\]
\end{dfn}

\begin{prp}\label{prp:Mf}
Let $(X,\cE,m,\Adm)$ be a CZ space and $Q_0 \in \Adm$.
\begin{enumerate}[label=(\roman*)]
\item\label{en:Mf_meas} For any $f \in L^0(X)$, the function $\cM_{Q_0} f$ is measurable and nonnegative.
\item\label{en:Mf_bdd} $\cM_{Q_0}$ is of weak type $(1,1)$ and $L^p(X)$-bounded for all $p \in (1,\infty]$, with universal bounds independent of $Q_0$, and in particular
\[
\|\cM_{Q_0} f\|_\infty \leq \|f\|_\infty, \qquad \sup_{\lambda > 0} \lambda m\{ \cM_{Q_0} f > \lambda \} \leq \|f\|_1.
\]
\item\label{en:Mf_ptw} For any $f \in L^0(X)$,
\[
|f| \chr_{Q_0} \leq \cM_{Q_0} f \quad\text{$m$-almost everywhere.}
\]
\end{enumerate}
\end{prp}
\begin{proof}
\ref{en:Mf_meas}.
Nonnegativity is clear, and measurability follows from the fact that $\cM_{Q_0} f = \sup_n \cM_{Q_0,n} f$, where
\[
\cM_{Q_0,n} f \defeq \sum_{Q \in \Ch^n(Q_0)} \left(\dashint_Q |f| \,dm\right) \,\chr_Q
\]
are measurable functions.

\ref{en:Mf_bdd}. As $\|\cM_{Q_0} f\|_\infty \leq \|f\|_\infty$ by definition, it only remains to check the weak-type $(1,1)$ bound, as the remaining boundedness properties follow by interpolation.

On the other hand, for any $\lambda > 0$, if
\[
\cQ_\lambda \defeq \left\{ Q \in \Ch^*(Q_0) \tc \dashint_Q |f| \,dm > \lambda \right\}
\]
then
\[
\{ \cM_{Q_0} f > \lambda \} = \bigcup \cQ_\lambda = \bigcup \tilde \cQ_\lambda,
\]
where $\tilde\cQ_\lambda$ is the set of maximal elements of $\cQ_\lambda$. As $\tilde\cQ_\lambda$ is an at most countable disjoint collection,
\[
m\{ \cM_{Q_0} f > \lambda \} = \sum_{Q \in \tilde\cQ_\lambda} m(Q) \leq \sum_{Q \in \tilde\cQ_\lambda} \frac{1}{\lambda} \int_Q |f| \,dm \leq \frac{1}{\lambda} \|f\|_1,
\]
as claimed.

\ref{en:Mf_ptw}. As $\cM_{Q_0} f = \cM_{Q_0} (f \chr_{Q_0})$, it is enough to prove the inequality for functions $f$ vanishing outside $Q_0$. In addition, if $\int_{Q_0} |f| \,dm = \infty$, then $\cM_{Q_0} f|_{Q_0} \equiv \infty$, so the inequality holds trivially. So it only remains to prove the inequality for $f \in L^1(Q_0)$.

Now, if $f \in \Span \{ \chr_Q \tc Q \in \Ch^*(Q_0) \}$, then we can write
\[
f = \sum_{Q \in \cQ} c_Q \chr_{Q}
\]
for some finite disjoint collection $\cQ \subseteq \Ch^*(Q_0)$; as
\[
|f(x)| = |c_Q| = \dashint_Q |f| \,dm \leq \cM_{Q_0} f(x)
\]
for any $x \in Q \in \cQ$, while $f(x) = 0$ if $x \notin \bigcup \cQ$, the claimed inequality holds true for any such $f$.

For an arbitrary $f \in L^1(Q_0)$, by Proposition \ref{prp:CZbasic}\ref{en:CZbasic_span} we can find functions $f_n \in \Span \{ \chr_Q \tc Q \in \Ch^*(Q_0) \}$ such that $f_n \to f$ in $L^1$; so, by part \ref{en:Mf_bdd} we also have $\cM_{Q_0}f_n \to \cM_{Q_0} f$ in $L^{1,\infty}$ and, up to replacing $f_n$ with a subsequence, we may also assume that $f_n \to f$ and $\cM_{Q_0}f_n \to \cM_{Q_0} f$ pointwise almost everywhere. The desired inequality for $f$ thus follows from the analogous ones for the $f_n$.
\end{proof}

\subsection{Calder\'on--Zygmund decomposition}
We now proceed with the construction of Calder\'on--Zygmund decompositions of integrable functions on CZ spaces.

\begin{lem}\label{lem:higherint}
Let $(X,\cE,m,\Adm)$ be a CZ space. Let $Q_0 \in \Adm$. Let $\cQ \subseteq \Ch^*(Q_0) \setminus \{Q_0\}$ be a disjoint collection, and for any $Q \in \cQ$ let $\Fth Q \in \Ch^*(Q_0)$ be given so that $Q \in \Ch(\Fth Q)$. Let $\lambda > 0$ and $f \in L^0(X)$ be such that
\begin{equation}\label{eq:stopping_cond}
\dashint_{\Fth Q} |f| \,dm \leq \lambda \qquad \forall Q \in \cQ,
\end{equation}
and define
\[
h \defeq \sum_{Q \in \cQ} \left( \int_{Q} |f| \,dm \right) \frac{1}{m(\Fth Q)} \chr_{\Fth Q}.
\]
Then, for all $p \in [1,\infty)$,
\begin{equation}\label{eq:higherint}
\|h\|_p^p \lesssim_p \lambda^{p-1} \|h\|_1, \qquad \|h\|_1 = \| f\chr_{\bigcup \cQ} \|_1,
\end{equation}
with universal implicit constants.
\end{lem}
\begin{proof}
The right-hand equality in \eqref{eq:higherint} is clear from the definition of $h$ and the disjointness of $\cQ$. So it remains to prove the left-hand inequality; by interpolation, it is enough to prove it for $p \in \Npos \defeq \NN \setminus \{0\}$. For this, we follow the proof of \cite[Lemma 3.3]{LSMP}.

Notice first that, for any $p \in \Npos$,
\begin{equation}\label{eq:hLp}
\begin{split}
\|h\|_p^p &= \int |h|^p \,dm \\
&= \int \sum_{Q_1,\dots,Q_p \in \cQ} \left(\int_{Q_1} |f| \,dm\right) \cdots \left(\int_{Q_p} |f| \,dm\right) \frac{\chr_{\Fth Q_1 \cap \dots \cap \Fth Q_p}}{m(\Fth Q_1) \cdots m(\Fth Q_p)} \,dm \\
&= \sum_{Q_1,\dots,Q_p \in \cQ} \left(\int_{Q_1} |f| \,dm\right) \cdots \left(\int_{Q_p} |f| \,dm\right) \frac{m(\Fth Q_1 \cap \dots \cap \Fth Q_p)}{m(\Fth Q_1) \cdots m(\Fth Q_p)}.
\end{split}
\end{equation}

Assume now that $p > 1$. Recall that all the $\Fth Q_j$ above are elements of $\Ch^*(Q_0)$, so they are disjoint unless they are contained in one another, and in that case their intersection is the smallest one of them. So, up to removing zero terms in the above sum (where $\Fth Q_1 \cap \dots \cap \Fth Q_p = \emptyset$) and reordering, we may assume that $\Fth Q_p \subseteq \Fth Q_1, \dots, \Fth Q_{p-1}$ and therefore $\Fth Q_1 \cap \dots \cap \Fth Q_p = \Fth Q_p$; thus,
\[\begin{split}
\|h\|_p^p 
&\leq p \sum_{\substack{Q_1,\dots,Q_p \in \cQ \\ \Fth Q_p \subseteq \Fth Q_1,\dots, \Fth Q_{p-1}}} \left(\int_{Q_1} |f| \,dm\right) \cdots \left(\int_{Q_p} |f| \,dm\right) \frac{1}{m(\Fth Q_1) \cdots m(\Fth Q_{p-1})} \\
&= p \sum_{Q_1,\dots,Q_{p-1} \in \cQ} \left( \sum_{\substack{Q \in \cQ \\ \Fth Q \subseteq \Fth Q_1,\dots,\Fth Q_{p-1}}}  \int_{Q} |f| \,dm\right) \\
&\qquad\times\left(\int_{Q_1} |f| \,dm\right) \cdots \left(\int_{Q_{p-1}} |f| \,dm\right) \frac{1}{m(\Fth Q_1) \cdots m(\Fth Q_{p-1})}.
\end{split}\]
For similar reasons, if $\Fth Q_1 \cap \dots \cap \Fth Q_{p-1} \neq \emptyset$, then $\Fth Q_1 \cap \dots \cap \Fth Q_{p-1} = \Fth R$ for some $R \in \{Q_1,\dots,Q_{p-1}\} \subseteq \cQ$, so
\[
\sum_{\substack{Q \in \cQ \\ \Fth Q \subseteq \Fth Q_1,\dots,\Fth Q_{p-1}}}  \int_{Q} |f| \,dm 
\leq \int_{\Fth R} |f| \,dm 
\leq \lambda m(\Fth R) = \lambda m(\Fth Q_1 \cap \dots \cap \Fth Q_{p-1}),
\]
where we used the disjointness of $\cQ$ and inequality \eqref{eq:stopping_cond} applied to $R \in \cQ$; thus
\[\begin{split}
\|h\|_p^p 
&\leq p \lambda \sum_{Q_1,\dots,Q_{p-1} \in \cQ} \left(\int_{Q_1} |f| \,dm\right) \cdots \left(\int_{Q_{p-1}} |f| \,dm\right) \frac{m(\Fth Q_1 \cap \dots \cap \Fth Q_{p-1})}{m(\Fth Q_1) \cdots m(\Fth Q_{p-1})} \\
&= p \lambda \|h\|_{p-1}^{p-1},
\end{split}\]
where we used \eqref{eq:hLp} again.

Therefore, inductively we deduce that, for any $p \in \Npos$,
\[
\|h\|_p^p \leq p! \, \lambda^{p-1} \|h\|_1 ,
\]
as desired.
\end{proof}

\begin{lem}\label{lem:stoppingsets}
Let $(X,\cE,m,\Adm)$ be a CZ space. Let $f \in L^1(X)$, $Q_0 \in \Adm$ and $\lambda \in (0,\infty)$ be such that
\[
\dashint_{Q_0} |f| \,dm \leq \lambda.
\]
Then, there exists a disjoint collection $\cQ \subseteq \Ch^*(Q_0) \setminus \{Q_0\}$ and, for any $Q \in \cQ$, a given $\Fth Q \in \Ch^*(Q_0)$ with $Q \in \Ch(\Fth Q)$, satisfying the following properties:
\begin{equation}\label{eq:stopping_basic}
\bigcup \cQ = \{ \cM_{Q_0} f > \lambda \}, \qquad \sup_{Q \in \cQ} \dashint_{\Fth Q} |f| \,dm \leq \lambda.
\end{equation}
In particular,
\begin{equation}\label{eq:stopping_adv}
\sum_{Q \in \cQ} m(Q) = m\left(\bigcup \cQ\right) \leq \|f\|_1 / \lambda, \qquad \|f \chr_{Q_0 \setminus \bigcup \cQ} \|_{\infty} \leq \lambda.
\end{equation}
\end{lem}
\begin{proof}
We shall call $Q \in \Ch^*(Q_0)$ \emph{stopping} if $\dashint_Q |f| \,dm > \lambda$, while $\dashint_{R} |f| \,dm \leq \lambda$ for any ancestor $R \in \Ch^*(Q_0)$ of $Q$ (i.e., any $R \in \Ch^*(Q_0)$ such that $Q \in \Ch^*(R) \setminus \{R\}$). By our assumptions, $Q_0$ itself is not stopping. 
Let $\cQ \subseteq \Ch^*(Q_0) \setminus \{Q_0\}$ be the collection of the stopping sets.

Notice that, if two distinct elements of $\Ch^*(Q_0)$ intersect, then one of them is an ancestor of the other. Then, clearly, $\cQ$ is a disjoint collection, and for any $Q \in \cQ$ there is $\Fth Q \in \Ch^*(Q_0)$ with $Q \in \Ch(\Fth Q)$ and 
\[
\dashint_{Q} |f| \,dm > \lambda, \qquad \dashint_{\Fth Q} |f| \,dm \leq \lambda.
\]
Moreover, for any $x \in Q_0 \setminus \bigcup \cQ$ and any $Q \in \Ch^*(Q_0)$ with $x \in Q$, we must have $\dashint_Q |f| \,dm \leq \lambda$ (otherwise, $Q$ or some ancestor thereof would be stopping, thus contradicting the fact that $x \notin \bigcup \cQ$). So, if we set $\Omega \defeq X \setminus \bigcup \cQ$, then
\[
(\cM_{Q_0} f)|_\Omega \leq \lambda, \qquad (\cM_{Q_0} f)|_{X \setminus \Omega} > \lambda
\]
pointwise. This concludes the proof of \eqref{eq:stopping_basic}.

As for \eqref{eq:stopping_adv}, the first statement is due to the disjointness of $\cQ$ and the weak-type $(1,1)$ bound for $\cM_{Q_0} f$ from Proposition \ref{prp:Mf}\ref{en:Mf_bdd}, while the second one follows from the pointwise bound from Proposition \ref{prp:Mf}\ref{en:Mf_ptw}.
\end{proof}

\begin{prp}\label{prp:CZdec}
Let $(X,\cE,m,\Adm)$ be a CZ space. Let $f \in L^1(X)$ and let $\lambda > \|f\|_1/m(X)$. Then, there exists a decomposition
\[
f = g + \sum_{j \in J} b_j,
\]
where $J$ is at most countable,
satisfying the following properties:
\begin{enumerate}[label=(\alph*)]
\item\label{en:CZdec_abQ} for any $j \in J$,
\[
b_j = a_j - \left(\int a_j \,dm\right) \frac{1}{m(\Fth Q_j)} \chr_{\Fth Q_j}
\]
for some $Q_j,\Fth Q_j \in \Adm$ with $Q_j \in \Ch(\Fth Q_j)$ and $a_j \in L^1(X)$ with $\{ a_j \neq 0 \} \subseteq Q_j$;
\item\label{en:CZdec_est} the estimates
\[
\sum_{j \in J} m(Q_j) \leq \|f\|_1/\lambda, \qquad \sum_{j \in J} \|a_j\|_1 \leq \|f\|_1
\]
hold true;
\item\label{en:CZdec_estp}
for any $p \in [1,\infty)$, 
\[
\|g\|_p^p \lesssim_p \lambda^{p-1} \|f\|_1
\]
with universal implicit constants;
\item\label{en:CZdec_est2} for all $j \in J$,
\[
\{ b_j \neq 0 \} \subseteq \Fth Q_j, \quad \int_{\Fth Q_j} b_j \,dm = 0, \quad \|b_j\|_1 \leq 2\|a_j\|_1, \quad
\sum_{j \in J} \|b_j\|_1 \leq 2 \|f\|_1.
\]
\end{enumerate}
\end{prp}
\begin{proof}
Notice that properties \ref{en:CZdec_est2} immediately follow from \ref{en:CZdec_abQ} and \ref{en:CZdec_est}. So we shall only discuss \ref{en:CZdec_abQ}, \ref{en:CZdec_est} and \ref{en:CZdec_estp} in the rest of the proof.

\medskip

We first show that it is enough to prove the result when $f$ is well-supported. Let $(Q_n)_{n \in \NN}$ be an exhausting sequence of admissible sets. Up to replacing it with a subsequence, we may assume that
\[
\|f \chr_{X \setminus Q_n} \|_1 \leq 2^{-n-1} \|f\|_1 \qquad \forall n \in \NN.
\]
In particular, if we set $f_0 \defeq f \chr_{Q_0}$ and $f_n \defeq f \chr_{Q_n \setminus Q_{n-1}}$ for $n > 0$, then
\[
f = \sum_{n \in \NN} f_n,
\]
where each $f_n$ is well-supported and
\begin{equation}\label{eq:exp_decay}
\|f_n\|_1 \leq 2^{-n} \|f\|_1 \qquad \forall n \in \NN.
\end{equation}
As $\lambda > \|f\|_1/m(X) \geq \|f_n\|_1/m(X)$, if the result is known for well-supported functions, for any $n \in \NN$ we obtain a decomposition
\[
f_n = g_n + \sum_{j \in J_n} b_{n,j}
\]
with the above properties. By summing in $n \in \NN$, we also have the decomposition
\[
f = g + \sum_{n \in \NN} \sum_{j \in J_n} b_{n,j}, \qquad g \defeq \sum_{n \in \NN} g_n,
\]
and since $\|f\|_1 = \sum_{n \in \NN} \|f_n\|_1$, we immediately obtain that such a decomposition of $f$ satisfies properties \ref{en:CZdec_abQ} and \ref{en:CZdec_est};
moreover, by \eqref{eq:exp_decay}, for all $p \in [1,\infty)$,
\[
\|g\|_p^p \leq \left(\sum_{n \in \NN} \|g_n\|_p\right)^p \lesssim_p \lambda^{p-1}\|f\|_1  \left( \sum_{n \in \NN} 2^{-n/p} \right)^p \lesssim_p \lambda^{p-1} \|f\|_1,
\]
thus property \ref{en:CZdec_estp} holds too.

\medskip

Thus, it remains to prove the result when $f \in L^1(X)$ is well-supported. As $m(X) > \|f\|_1/\lambda$ and $f$ is well-supported, we can find $Q_0 \in \Adm$ such that $\{ f \neq 0 \} \subseteq Q_0$ and $m(Q_0) > \|f\|_1/\lambda$. In particular,
\[
\dashint_{Q_0} |f| \,dm \leq \frac{1}{m(Q_0)} \|f\|_1 < \lambda.
\]
Thus, we can apply Lemma \ref{lem:stoppingsets} to obtain a disjoint collection $\cQ \subseteq \Ch^*(Q_0) \setminus \{Q_0\}$ and, for each $Q \in \cQ$, a $\Fth Q \in \Ch^*(Q_0)$ with $Q \in \Ch(\Fth Q)$, with the stated properties. Consequently, as $\{ f \neq 0 \} \subseteq Q_0$, if we set $\Omega \defeq X \setminus \bigcup \cQ$, from \eqref{eq:stopping_adv} we deduce that
\begin{equation}\label{eq:fOmega_estp}
\|f\chr_\Omega\|_p^p \leq \|f\|_1 \|f \chr_\Omega\|_\infty^{p-1} \leq \lambda^{p-1} \|f\|_1
\end{equation}
for all $p \in [1,\infty)$, and moreover
\[
\sum_{Q \in \cQ} m(Q) \leq \|f\|_1/\lambda.
\]
In addition, from \eqref{eq:stopping_basic} we have
\[
\dashint_{\Fth Q} |f| \,dm \leq \lambda \qquad \forall Q \in \cQ.
\]

Now, if we define
\[
a_Q \defeq f \chr_Q, \qquad b_Q \defeq a_Q - \left(\int a_Q \,dm \right) \frac{1}{m(\Fth Q)} \chr_{\Fth Q}
\]
for any $Q \in \cQ$, then clearly $\{a_Q \neq 0\} \subseteq Q$ and
\[
\sum_{Q \in \cQ} \|a_Q\|_1 \leq \|f\|_1.
\]
So, if we define
\begin{equation}\label{eq:gdef}
g \defeq f - \sum_{Q \in \cQ} b_Q = f \chr_\Omega + \sum_{Q \in \cQ} \left( \int_Q f \,dm \right) \frac{1}{m(\Fth Q)} \chr_{\Fth Q},
\end{equation}
then we obtain a decomposition
\[
f = g + \sum_{Q \in \cQ} b_Q
\]
satisfying \ref{en:CZdec_abQ} and \ref{en:CZdec_est}.

It only remains to verify \ref{en:CZdec_estp}. Set
\[
h \defeq \sum_{Q \in \cQ} \left( \int_Q |f| \,dm \right) \frac{1}{m(\Fth Q)} \chr_{\Fth Q},
\]
so $h$ is nonnegative and, by \eqref{eq:gdef},
\[
|g| \leq |f| \chr_\Omega + h.
\]
Thus, property \ref{en:CZdec_estp} follows from \eqref{eq:fOmega_estp} and the estimate
\[
\|h\|_p^p \lesssim_p \lambda^{p-1} \|f\|_1
\]
for all $p \in [1,\infty)$, which in turn is a consequence of Lemma \ref{lem:higherint}.
\end{proof}

Proposition \ref{prp:CZdec} is the analogue for CZ spaces of the decompositions used in \cite{CPW,CADMMRV2025}, where a similar statement is formulated in the special case of martingales. We can now prove Theorem \ref{thm:mainCZdec}, which upgrades the estimate for the good part to match exactly the classical case.

\begin{proof}[Proof of Theorem \ref{thm:mainCZdec}.] 
Set $f^{(0)} \defeq f$. We now construct inductively a sequence of functions $f^{(n)}$. Namely, for any $n \in \NN$, we apply Proposition \ref{prp:CZdec} to the function $f^{(n)}$ in place of $f$, and the level $2^{-n}\lambda$ in place of $\lambda$, in order to obtain a decomposition
\[
f^{(n)} = g^{(n)} + \sum_{j \in J_n} b^{(n)}_{j},
\]
where $b^{(n)}_j = a^{(n)}_j - (\int a^{(n)}_j \,dm) \chr_{\Fth Q^{(n)}_j}/m(\Fth Q^{(n)}_j)$ for some $Q^{(n)}_j,\Fth Q^{(n)}_j \in \Adm$ with $Q^{(n)}_j \in \Ch(\Fth Q^{(n)}_j)$ and $a^{(n)}_j \in L^1(X)$ with $\{ a^{(n)}_j \neq 0 \} \subseteq Q^{(n)}_j$, such that
\[
\sum_{j \in J_n} m(Q^{(n)}_j) \leq 2^n \|f^{(n)}\|_1/\lambda , \qquad \sum_{j \in J_n} \|a^{(n)}_j\|_1 \leq \|f^{(n)}\|_1,
\]
and
\[
\|g^{(n)}\|_1 \leq 3 \|f^{(n)}\|_1, \qquad  \|g^{(n)}\|_2^2 \leq C 2^{-n} \lambda \|f^{(n)}\|_1,
\]
where $C$ is the universal implicit constant in the estimate for $p=2$ of Proposition \ref{prp:CZdec}\ref{en:CZdec_estp}.
If we now set $h^{(n)} := g^{(n)} \chr_{\{ |g^{(n)}| \leq 4 C\lambda/2^n\}}$ and $f^{(n+1)} \defeq g^{(n)}-h^{(n)}$, then we obtain the decomposition
\begin{equation}\label{eq:2termdec}
f^{(n)} = f^{(n+1)} + h^{(n)} + \sum_{j \in J_n} b^{(n)}_{j},
\end{equation}
where
\[
\|h^{(n)}\|_1 \leq 3 \|f^{(n)}\|_1, \qquad \|h^{(n)}\|_\infty \leq 4 C \lambda/2^n, 
\]
and
\[
\|f^{(n+1)}\|_1 \leq m(\{|g^{(n)}| > 4 C\lambda/2^n\})^{1/2} \|g^{(n)}\|_2 \leq \|f^{(n)}\|_1 / 4.
\]
Notice that, inductively, one gets
\[
\|f^{(n)}\|_1 \leq 4^{-n} \|f\|_1
\]
for any $n \in \NN$; in particular, $\|f^{(n)}\|_1/m(X) \leq 4^{-n} \|f\|_1/m(X) < 2^{-n} \lambda$, so Proposition \ref{prp:CZdec} is indeed applicable at each step. Moreover, $f^{(n)} \to 0$ in $L^1(X)$ as $n \to \infty$; so, by combining the decompositions \eqref{eq:2termdec} of the $f^{(n)}$, we obtain the decomposition
\[
f = g + \sum_{n \in \NN} \sum_{j \in J_n} b^{(n)}_j
\]
of $f$, where $g \defeq \sum_{n \in \NN} h^{(n)}$ and
\[
\|g\|_1 \leq \sum_{n \in \NN} \|h^{(n)}\|_1 \leq 4 \|f\|_1, \qquad \|g\|_\infty \leq \sum_{n \in \NN} \|h^{(n)}\|_\infty \leq 8 C \lambda,
\]
while
\[
\sum_{n \in \NN} \sum_{j \in J_n} m(Q^{(n)}_j) \leq 2 \|f\|_1/\lambda, \qquad \sum_{n \in \NN} \sum_{j\in J_n} \|a^{(n)}_j\|_1 \leq \frac{4}{3} \|f\|_1.
\]
From the above estimates it is readily seen that this decomposition of $f$ satisfies the properties stated in Theorem \ref{thm:mainCZdec}.
\end{proof}

With the Calder\'on--Zygmund decomposition at hand, we are also ready to prove Theorem \ref{thm:mainCZwt}.

\begin{proof}[Proof of Theorem \ref{thm:mainCZwt}]
By interpolation, it is enough to prove the weak-type $(1,1)$ bound.

Let $f \in L^1(X) \cap L^q(X)$ and $\lambda > 0$. If $\lambda \leq \|f\|_1/m(X)$, then trivially
\[
m\{ |Tf| > \lambda \} \leq m(X) \leq \frac{\|f\|_1}{\lambda}.
\]

Assume instead $\lambda > \|f\|_1/m(X)$, and let
\[
f = g + \sum_{j \in J} b_j
\]
be the decomposition of $f$ at level $\lambda$ given by Proposition \ref{prp:CZdec}; recall that, for any $j \in J$,
\[
b_j = a_j - \left( \int a_j \,dm \right) \frac{1}{m(\Fth Q_j)} \chr_{\Fth Q_j}
\]
for some $Q_j,\Fth Q_j \in \Adm$ with $Q_j \in \Ch(\Fth Q_j)$, and $a_j \in L^1(X)$ with $\{ a_j \neq 0 \} \subseteq Q_j$.

Set $b \defeq \sum_{j \in J} b_j$ and $\Omega \defeq \bigcup_{j \in J} Q_j^*$. Then,
\[\begin{split}
m\{ |Tf| > \lambda \} &\leq m\{ |Tg| > \lambda/2 \} + m\{ |Tb| > \lambda/2 \} \\
&\leq m\{ |Tg| > \lambda/2 \} + m(\Omega) + m(\{|Tb| > \lambda/2\} \setminus \Omega).
\end{split}\]

Now, by Markov's inequality,
\[
m\{ |Tg| > \lambda/2 \} \leq (2/\lambda)^q \|Tg\|_q^q \leq (2/\lambda)^q \|T\|_{q \to q}^q \|g\|_q^q \lesssim_{q} \|T\|_{q \to q}^q \frac{\|f\|_1}{\lambda},
\]
where we used the $L^q$-bound for $g$ from Proposition \ref{prp:CZdec}\ref{en:CZdec_estp}.

Moreover,
\[
m(\Omega) \leq \sum_{j \in J} m(Q^*_j) \leq \kappa \sum_{j \in J} m(Q_j) \leq \kappa \frac{\|f\|_1}{\lambda},
\]
where we used \eqref{eq:extQ} and the first bound from Proposition \ref{prp:CZdec}\ref{en:CZdec_est}.

Finally, again by Markov's inequality,
\[
m(\{|Tb| > \lambda/2\} \setminus \Omega) \leq \frac{2}{\lambda} \int_{X \setminus \Omega} |Tb| \,dm \leq \frac{2}{\lambda} \sum_{j \in J} \int_{X \setminus Q_j^*} |Tb_j| \,dm.
\]

Now, for each $j \in J$,
\begin{multline*}
\int_{X \setminus Q_j^*} |Tb_j| \,dm \\
\leq \int_{X \setminus \Fth Q_j^*} |Tb_j| \,dm 
+ \int_{\Fth Q_j^* \setminus Q_j^*} |Ta_j| \,dm 
+ \frac{\|a_j\|_1}{m(\Fth Q_j)} \int_{\Fth Q_j^* \setminus Q_j^*} |T\chr_{\Fth Q_j}| \,dm
\end{multline*}

Since $\int b_j \,dm =0$ and $\{ b_j \neq 0 \} \subseteq \Fth Q_j$, if $z_j$ is any element of $\Fth Q_j$,
\begin{multline*}
\int_{X \setminus \Fth Q_j^*} |Tb_j| \,dm = \int_{X \setminus \Fth Q_j^*} \left|\int_{\Fth Q_j} (K(x,y)-K(x,z_j)) b_j(y) \,dm(y) \right| \,dm(x) \\
\leq C_{H}(T) \|b_j\|_1,
\end{multline*}
so, by Proposition \ref{prp:CZdec}\ref{en:CZdec_est2},
\[
\sum_{j \in J} \int_{X \setminus \Fth Q_j^*} |Tb_j| \,dm \leq C_H(T) \sum_{j \in J} \|b_j\|_1 \leq 2 C_H(T) \|f\|_1.
\]

Moreover, as $\{ a_j \neq 0\} \subseteq Q_j$,
\[
\int_{\Fth Q_j^* \setminus Q_j^*} |Ta_j| \,dm  = \int_{\Fth Q_j^* \setminus Q_j^*} \left|\int_{Q_j} K(x,y) a_j(y) \,dm(y) \right| \,dm(x) \leq C_S(T) \|a_j\|_1,
\]
so, by Proposition \ref{prp:CZdec}\ref{en:CZdec_est},
\[
\sum_{j \in J} \int_{\Fth Q_j^* \setminus Q_j^*} |Ta_j| \,dm \leq C_S(T) \sum_{j \in J} \|a_j\|_1 \leq C_S(T) \|f\|_1.
\]

Finally, by H\"older's inequality and \eqref{eq:extQ},
\begin{multline*}
\frac{\|a_j\|_1}{m(\Fth Q_j)} \int_{\Fth Q_j^* \setminus Q_j^*} |T\chr_{\Fth Q_j}| \,dm \leq \frac{\|a_j\|_1}{m(\Fth Q_j)} m(\Fth Q_j^*)^{1/q'} \|T \chr_{\Fth Q_j}\|_q \\
\leq \frac{\|a_j\|_1}{m(\Fth Q_j)} m(\Fth Q_j^*)^{1/q'} \|T\|_{q \to q} m(\Fth Q_j)^{1/q} \leq \kappa^{1/q'} \|T\|_{q \to q} \|a_j\|_1,
\end{multline*}
so, again by Proposition \ref{prp:CZdec}\ref{en:CZdec_est},
\[
\sum_{j \in J} \frac{\|a_j\|_1}{m(\Fth Q_j)} \int_{\Fth Q_j^* \setminus Q_j^*} |T\chr_{\Fth Q_j}| \,dm \leq
\kappa^{1/q'} \|T\|_{q \to q} \sum_{j \in J} \|a_j\|_1 \leq \kappa^{1/q'} \|T\|_{q \to q} \|f\|_1.
\]

By combining the above estimates, we obtain the required weak-type $(1,1)$ bound for $T$.
\end{proof}

\begin{rem}\label{rem:size_smoothness}
Let $(X,\cE,m,\Adm)$ be a CZ space with enlargements as in \eqref{eq:extQ}.
Let $T : L^q(X) \to L^q(X)$ be bounded for some $q \in (1,\infty)$ and let $K_T : X \times X \to \CC$ satisfy \eqref{eq:CZop_offdiag} and \eqref{eq:CZop_horm}. Notice that, for any $Q,\Fth Q \in \Adm$ with $Q \in \Ch(\Fth Q)$,
\[
\esssup_{y \in Q} \int_{\Fth Q^* \setminus Q^*} |K_T(x,y)| \,dm(x) 
= \sup_{\|f\|_1 \leq 1} \|\chr_{\Fth Q^* \setminus Q^*} T (\chr_Q f)\|_{1},
\]
and moreover
\[\begin{split}
&\|\chr_{\Fth Q^* \setminus Q^*} T (\chr_Q f)\|_{1} \\
&\leq \left\|\chr_{\Fth Q^* \setminus Q^*} T \left(\left(f-\dashint_Q f \,dm\right)\chr_Q\right)\right\|_{1}
+ \left|\dashint_Q f \,dm\right| \|\chr_{\Fth Q^* \setminus Q^*} T\chr_Q\|_{1}.
\end{split}\]
Now, as $(f-\dashint_Q f \,dm)\chr_Q$ has zero integral, for any $z \in Q$,
\[\begin{split}
&\left\|\chr_{\Fth Q^* \setminus Q^*} T \left(\left(f-\dashint_Q f \,dm\right)\chr_Q\right)\right\|_{1} \\
&= \int_{\Fth Q^* \setminus Q^*} \left| \int_Q (K_T(x,y)-K_T(x,z)) \left(f(y)-\dashint_Q f \,dm\right) \,dm(y) \right| \,dm(x) \\
&\leq 2 C_H(T) \, \|f\|_1 .
\end{split}\]
Therefore, if $C_S(T)$ is defined as in \eqref{eq:CZop_size}, then from the above estimates we see that
\[
C_S(T) \leq 2 C_H(T) + \sup_{\Fth Q \in \Adm} \sup_{Q \in \Ch(\Fth Q)} m(Q)^{-1} \|\chr_{\Fth Q^* \setminus Q^*} T\chr_Q\|_{1},
\]
whence it follows that, under \eqref{eq:CZop_horm}, the size condition \eqref{eq:CZop_size} is equivalent to
\[
\widetilde{C}_S(T) \defeq \sup_{\Fth Q \in \Adm} \sup_{Q \in \Ch(\Fth Q)} m(Q)^{-1} \|\chr_{\Fth Q^* \setminus Q^*} T\chr_Q\|_{1} < \infty.
\]

Notice now that, by H\"older's inequality and the $L^q$-boundedness of $T$,
\[
\|\chr_{\Fth Q^* \setminus Q^*} T\chr_Q\|_{1} \leq m(\Fth Q^* \setminus Q^*)^{1/q'} \|T \chr_Q\|_{q} \leq \kappa m(\Fth Q)^{1/q'} \|T\|_{q \to q} m(Q)^{1/q},
\]
and therefore
\[
\widetilde{C}_S(T) \leq \kappa \|T\|_{q \to q} \sup_{\Fth Q \in \Adm} \sup_{Q \in \Ch(\Fth Q)} (m(\Fth Q)/m(Q))^{1/q'}.
\]
Thus, when the CZ space $X$ is regular in the sense of Definition \ref{dfn:regularCZ},
\[
\widetilde{C}_S(T) \leq \kappa \rho(X)^{1/q'} \|T\|_{q \to q},
\]
so the size condition \eqref{eq:CZop_size} is automatically satisfied.

Therefore, for a regular CZ space, one can drop the assumption \eqref{eq:CZop_size} from Theorem \ref{thm:mainCZwt}, as well as in Theorem \ref{thm:mainHardybdd},
and still obtain the claimed bounds for the operator $T$, with implicit constants also depending on the constant $\rho(X)$ in Definition \ref{dfn:regularCZ}. In other words, in the case of regular CZ spaces, those results match other results in the literature for the boundedness of singular integral operators (see, e.g., \cite{HS} or \cite[Theorem 4.1]{V}), where only the H\"ormander-type smoothness condition \eqref{eq:CZop_horm} is required.
\end{rem}

\section{Hardy and BMO spaces}\label{sec:sec2}

In this section, we prove Theorems \ref{thm:mainHardyequiv}, \ref{thm:mainHardyinterp} and \ref{thm:mainHardybdd}. We start by developing the necessary $H^1$--$\BMO$ theory including John--Nirenberg, duality and interpolation, and we finish proving the endpoint estimates for integral operators. 

\subsection{BMO spaces}
We begin with the bounded mean oscillation spaces.

\begin{dfn}
Let $(X,\cE,m,\Adm)$ be a CZ space and $q \in [1,\infty)$. We set
\[
L^q_\loc(X) \defeq \{ f \in L^0(X) \tc f \chr_Q \in L^q(X) \ \forall Q \in \Adm \}.
\]
We define $\BMO_q(X)$ as the set of functions $f \in L^q_\loc(X)$ such that
\[
\|f\|_{\BMO_q} \defeq \sup_{\Fth Q \in \Adm} \sup_{Q \in \Ch(\Fth Q)} \left( \dashint_Q \left|f - \dashint_{\Fth Q} f \, dm\right|^q \, dm\right)^{1/q} < \infty.
\]
\end{dfn}

\begin{prp}\label{prp:BMO}
Let $(X,\cE,m,\Adm)$ be a CZ space and $q \in [1,\infty)$.
\begin{enumerate}[label=(\roman*)]
\item\label{en:BMO_Lqloclctvs} $L^q_\loc(X)$ is a locally convex topological vector space with the seminorms
\[
f \mapsto \left(\int_Q |f|^q \,dm \right)^{1/q}, \qquad Q \in \Adm.
\]
\item\label{en:BMO_preleq} For all $f \in L^q_\loc(X)$ and $Q \in \Adm$,
\begin{multline*}
\inf_{c \in \CC} \left(\dashint_Q |f-c|^q \,dm \right)^{1/q} \leq \left(\dashint_Q \left|f-\dashint_Q f \,dm\right|^q \,dm \right)^{1/q} \\
\leq 2 \inf_{c \in \CC} \left(\dashint_Q |f-c|^q \,dm \right)^{1/q} .
\end{multline*}
\item\label{en:BMO_preleqch} For all $f \in L^q_\loc(X)$, $\Fth Q \in \Adm$ and $Q \in \Ch(\Fth Q)$,
\begin{multline*}
\left(\dashint_Q \left|f-\dashint_{\Fth Q} f \,dm\right|^q \,dm \right)^{1/q} \\
\simeq \max \left\{ \left(\dashint_Q \left|f-\dashint_Q f \, dm\right|^q \,dm \right)^{1/q} , \left| \dashint_Q f \,dm -\dashint_{\Fth Q} f \,dm \right| \right\} ,
\end{multline*}
with universal implicit constants.
\item\label{en:BMO_normeq} For all $f \in L^q_\loc(X)$,
\[\begin{split}
&\|f\|_{\BMO_q} \\
&\simeq \max \left\{ \sup_{Q \in \Adm} \left( \dashint_Q \left| f- \dashint_Q f \,dm \right|^q \,dm\right)^{1/q} , \sup_{\Fth Q \in \Adm} \sup_{Q \in \Ch(\Fth Q)} \left| \dashint_Q f \,dm - \dashint_{\Fth Q} f \,dm \right| \right\} \\
&\simeq \max \left\{ \sup_{Q \in \Adm_\Ch} \left( \dashint_Q \left| f- \dashint_Q f \,dm \right|^q \,dm\right)^{1/q} , \sup_{\Fth Q \in \Adm} \sup_{Q \in \Ch(\Fth Q)} \left| \dashint_Q f \,dm - \dashint_{\Fth Q} f \,dm \right| \right\},
\end{split}\]
with universal implicit constants, where $\Adm_\Ch \defeq \bigcup_{Q \in \Adm} \Ch(Q)$.
\item\label{en:BMO_const} For all $f \in L^q_\loc(X)$, we have $\|f\|_{\BMO_q} = 0$ if and only if $f$ is constant.
\item\label{en:BMO_vss} $\BMO_q(X)$ is a vector subspace of $L^q_\loc(X)$, and $\|\cdot\|_{\BMO_q}$ is a seminorm on $\BMO_q(X)$.
\item\label{en:BMO_Linfty} $L^\infty(X) \subseteq \BMO_q(X)$ and $\|f\|_{\BMO_q} \leq 2 \|f\|_\infty$ for all $f \in L^\infty(X)$.
\item\label{en:BMO_incl} $\BMO_r(X) \subseteq \BMO_q(X)$ whenever $r \in [q,\infty)$, and $\|f\|_{\BMO_q} \leq \|f\|_{\BMO_r}$ for all $f \in \BMO_r(X)$.
\end{enumerate}
\end{prp}
\begin{proof}
\ref{en:BMO_Lqloclctvs}. This is clear from the definitions.

\ref{en:BMO_preleq}. The first inequality is obvious. As for the second one, by the triangle inequality and H\"older's inequality,
\begin{equation}\label{eq:c_ineq}
\left(\dashint_Q \left|f-\dashint_Q f \,dm\right|^q \,dm \right)^{1/q} \leq \left(\dashint_Q |f|^q \,dm \right)^{1/q} + \dashint_Q |f| \,dm \leq 2 \left(\dashint_Q |f|^q \,dm \right)^{1/q}.
\end{equation}
By applying the previous estimate with $f$ replaced by $f-c$ for $c \in \CC$, the left-hand side does not change, so we deduce that
\[
\left(\dashint_Q \left|f-\dashint_Q f \,dm\right|^q \,dm \right)^{1/q} \leq 2 \inf_{c \in \CC} \left(\dashint_Q |f-c|^q \,dm \right)^{1/q},
\]
as desired.

\ref{en:BMO_preleqch}. By part \ref{en:BMO_preleq},
\[
\left(\dashint_Q \left|f-\dashint_{\Fth Q} f \,dm\right|^q \,dm \right)^{1/q} \geq \frac{1}{2} \left(\dashint_Q \left|f-\dashint_Q f \,dm\right|^q \,dm \right)^{1/q},
\]
and moreover, by H\"older's inequality,
\[
\left(\dashint_Q \left|f-\dashint_{\Fth Q} f \,dm\right|^q \,dm \right)^{1/q} \geq \dashint_Q \left|f-\dashint_{\Fth Q} f\,dm\right| \,dm \geq \left|\dashint_{Q} f \,dm - \dashint_{\Fth Q} f \,dm\right|.
\]
This proves the estimate $\gtrsim$. As for the opposite estimate, by the triangle inequality,
\[
\left(\dashint_Q \left|f-\dashint_{\Fth Q} f \,dm\right|^q \,dm \right)^{1/q} \leq \left(\dashint_Q \left|f-\dashint_{Q} f \,dm\right|^q \,dm \right)^{1/q} + \left| \dashint_Q f \,dm - \dashint_{\Fth Q} f\,dm \right|.
\]

\ref{en:BMO_normeq}. By part \ref{en:BMO_preleqch},
\[\begin{split}
&\|f\|_{\BMO_q} \\
&\simeq \sup_{\Fth Q \in \Adm} \sup_{Q \in \Ch(\Fth Q)} \max \left\{ \left(\dashint_Q \left|f-\dashint_Q f \, dm\right|^q \,dm \right)^{1/q} , \left| \dashint_Q f \,dm -\dashint_{\Fth Q} f \,dm \right| \right\}\\
&= \max \left\{ \sup_{Q \in \Adm_\Ch} \left( \dashint_Q \left| f- \dashint_Q f \,dm \right|^q \,dm\right)^{1/q} , \sup_{\Fth Q \in \Adm} \sup_{Q \in \Ch(\Fth Q)} \left| \dashint_Q f \,dm - \dashint_{\Fth Q} f \,dm \right| \right\} \\
&\leq \max \left\{ \sup_{Q \in \Adm} \left( \dashint_Q \left| f- \dashint_Q f \,dm \right|^q \,dm\right)^{1/q} , \sup_{\Fth Q \in \Adm} \sup_{Q \in \Ch(\Fth Q)} \left| \dashint_Q f \,dm - \dashint_{\Fth Q} f \,dm \right| \right\}.
\end{split}\]
In order to conclude, it is therefore enough to prove that
\[
\sup_{Q \in \Adm} \left( \dashint_Q \left| f- \dashint_Q f \,dm \right|^q \,dm\right)^{1/q} \leq \|f\|_{\BMO_q}.
\]
On the other hand, for any $\Fth Q \in \Adm$,
\begin{multline*}
\int_{\Fth Q} \left|f- \dashint_{\Fth Q} f \,dm\right|^q \,dm = \sum_{Q \in \Ch(\Fth Q)} \int_Q \left| f - \dashint_{\Fth Q} f \,dm \right|^q \,dm \\
\leq \sum_{Q \in \Ch(\Fth Q)} m(Q) \|f\|_{\BMO_q}^q = m(\Fth Q) \|f\|_{\BMO_q}^q,
\end{multline*}
thus showing that
\[
\left( \dashint_{\Fth Q} \left|f- \dashint_{\Fth Q} f \,dm\right|^q \,dm \right)^{1/q} \leq \|f\|_{\BMO_q}
\]
for any $\Fth Q \in \Adm$, as required.

\ref{en:BMO_const}. If $\|f\|_{\BMO_q} = 0$, then $f$ is constant on any $Q \in \Adm$. By using an exhausting sequence, we see that $f$ is constant on the whole $X$.

\ref{en:BMO_vss}. This is easily checked.

\ref{en:BMO_Linfty}. This follows immediately from the definition.

\ref{en:BMO_incl}. This is clear by H\"older's inequality.
\end{proof}

\begin{rem}\label{rem:regularBMO}
Assume that $(X,\cE,m,\Adm)$ is regular, as in Definition \ref{dfn:regularCZ}. Then, for all $\Fth Q \in \Adm$ and $Q \in \Ch(\Fth Q)$,
\begin{multline*}
\left| \dashint_Q f \,dm - \dashint_{\Fth Q} f \,dm \right| \leq \dashint_Q \left| f  - \dashint_{\Fth Q} f \,dm \right| \,dm \leq \frac{m(\Fth Q)}{m(Q)} \dashint_{\Fth Q} \left| f  - \dashint_{\Fth Q} f \,dm \right| \,dm \\
\leq \rho(X) \, \dashint_{\Fth Q} \left| f  - \dashint_{\Fth Q} f \,dm \right| \,dm.
\end{multline*}
Thus, by Proposition \ref{prp:BMO}\ref{en:BMO_normeq}, for all $q \in [1,\infty)$ and $f \in L^q_\loc(X)$,
\[
\|f\|_{\BMO_q} \simeq_{\rho(X)} \sup_{Q \in \Adm} \left( \dashint_Q \left| f- \dashint_Q f \,dm \right|^q \,dm\right)^{1/q}.
\]
This shows that, when the CZ space is regular, the above-defined $\BMO$ spaces coincide with the ``more usual'' ones (cf., e.g., \cite[Section 3]{V}).
\end{rem}

\begin{lem}\label{lem:glueing}
Let $(X,\cE,m,\Adm)$ be a CZ space and $q \in [1,\infty]$. Assume that, for any $Q \in \Adm$, a function $g_Q \in L^q(Q)$ is given, so that $g_Q - g_R$ is constant on $Q \cap R$ for all $Q,R \in \Adm$. Then, there exists a function $g \in L^q_\loc(X)$ such that $g-g_Q$ is constant on $Q$ for all $Q \in \Adm$.
\end{lem}
\begin{proof}
Let $(Q_n)_{n \in \NN}$ be an exhausting sequence of admissible sets. For all $n \in \NN$, let $c_n \in \CC$ be the constant value of $g_{Q_0}-g_{Q_n}|_{Q_0}$, and define $g_n \defeq g_{Q_n}+c_n \in L^q(Q_n)$. Then $g_{n+1}|_{Q_n}-g_n = (g_{Q_{n+1}}|_{Q_n}-g_{Q_n}) + c_{n+1}-c_n$ is constant and vanishes on $Q_0$; hence $g_{n+1}|_{Q_n} = g_n$ for all $n \in \NN$, and therefore the $g_n$ can be ``glued together'' to define a single function $g \in L^0(X)$.

Similarly, we see that, for any $Q \in \Adm$, $(g-g_Q)|_{Q \cap Q_n} = (g_n-g_Q)|_{Q \cap Q_n}$ is constant for all $n \in \NN$; as $(Q_n)_n$ is an exhausting sequence, we conclude that $(g-g_Q)|_Q$ is constant. In particular, $g|_Q \in L^q(Q)$ for all $Q \in \Adm$, i.e., $g \in L^q_\loc(X)$. 
\end{proof}

\begin{prp}\label{prp:BMOcompl}
Let $(X,\cE,m,\Adm)$ be a CZ space and $q \in [1,\infty)$.
The quotient space $\BMO_q(X)/\CC$ with the norm induced by $\|\cdot\|_{\BMO_q}$ is a Banach space, continuously embedded in the locally convex topological vector space $L^q_\loc(X)/\CC$.
\end{prp}
\begin{proof}
By Proposition \ref{prp:BMO}\ref{en:BMO_const}, $\|f\|_{\BMO_q}$ vanishes if and only if $f$ is constant, so the seminorm $\|\cdot\|_{\BMO_q}$ descends to the quotient and induces a norm on $\BMO_q(X)/\CC$. In addition, from Proposition \ref{prp:BMO}\ref{en:BMO_Lqloclctvs}, we know that $L^q_\loc(X)$ is a locally convex topological vector space, and constant functions form a closed subspace of $L^q_\loc(X)$. Therefore, the quotient $L^q_\loc(X)/\CC$ is also a locally convex topological vector space with seminorms
\[
f \mapsto \inf_{c \in \CC} \left(\int_Q |f-c|^q \,dm\right)^{1/q}, \qquad Q \in \Adm.
\]
So, in light of Proposition \ref{prp:BMO}\ref{en:BMO_preleq} and the definition of $\|\cdot\|_{\BMO_q}$, it is clear that the inclusion $\BMO_q(X)/\CC \subseteq L^q_\loc(X)/\CC$ is continuous. It only remains to prove that $\BMO_q(X)/\CC$ is complete.

Let $(f_n+\CC)_{n \in \NN}$ be a Cauchy sequence in $\BMO_q(X)/\CC$. Then, by the characterisation of the $\BMO_q$-seminorm from Proposition \ref{prp:BMO}\ref{en:BMO_normeq}, for any $Q \in \Adm$, the sequence $(f_n|_Q-\dashint_Q f_n \,dm)_n$ is a Cauchy sequence in $L^q(Q)$, and therefore it converges in $L^q(Q)$ to a function $g_Q \in L^q(Q)$ satisfying $\dashint_Q g_Q \,dm= 0$.

Moreover, for any $Q,R \in \Adm$, the difference
\[
\left.\left(f_n|_Q-\dashint_Q f_n \,dm\right)\right|_{Q \cap R} - \left.\left(f_n|_R-\dashint_R f_n \,dm\right)\right|_{Q \cap R} = \dashint_R f_n \,dm - \dashint_Q f_n \,dm
\]
is constant; thus, by taking the limit as $n \to \infty$, we deduce that $g_Q|_{Q \cap R} - g_R|_{Q \cap R}$ is constant too.

We can then apply Lemma \ref{lem:glueing} to construct a function $g \in L^q_\loc(X)$ 
such that $g - g_Q$ is constant on $Q$ for all $Q \in \Adm$. As $\dashint_Q g_Q \,dm = 0$, we must have $(g-g_Q)|_Q = \dashint_Q g \,dm$ for all $Q \in \Adm$.

Now, for any $\Fth Q \in \Adm$, we have
\[
g|_{\Fth Q} - \dashint_{\Fth Q} g \,dm = g_{\Fth Q} = \lim_{n \to \infty} \left(f_n|_{\Fth Q}-\dashint_{\Fth Q} f_n \,dm\right)
\]
in $L^q(\Fth Q)$; if $Q \in \Ch(\Fth Q)$, by restricting to $Q$ we also obtain that
\begin{equation}\label{eq:lim_gQ}
g|_{Q} - \dashint_{\Fth Q} g \,dm = \lim_{n \to \infty} \left(f_n|_{Q}-\dashint_{\Fth Q} f_n \,dm\right)
\end{equation}
in $L^q(Q)$. In particular,
\begin{multline*}
\|g\|_{\BMO_q} = \sup_{\Fth Q \in \Adm} \sup_{Q \in \Ch(\Fth Q)} \left( \dashint_Q \left| g - \dashint_{\Fth Q} g \,dm \right|^{q} \,dm \right)^{1/q} \\
\leq \sup_{\Fth Q \in \Adm} \sup_{Q \in \Ch(\Fth Q)} \sup_n \left( \dashint_Q \left| f_n - \dashint_{\Fth Q} f_n \,dm \right|^{q} \,dm \right)^{1/q} = \sup_n \|f_n\|_{\BMO_q} < \infty,
\end{multline*}
because Cauchy sequences are bounded. This shows that $g \in \BMO_q(X)$.

Moreover, as $(f_n+\CC)_{n}$ is a Cauchy sequence in $\BMO_q(X)/\CC$, for any given $\epsilon>0$ we can find $N \in \NN$ such that
\[
\left( \dashint_Q \left| \left(f_j-\dashint_{\Fth Q} f_j \,dm \right) - \left(f_k-\dashint_{\Fth Q} f_k \,dm\right) \right|^q \,dm \right)^{1/q} \leq \epsilon
\]
for all $\Fth Q \in\Adm$ and $Q \in \Ch(\Fth Q)$, and all $j,k \geq N$; if we take the limit as $k \to \infty$ in the previous inequality and use \eqref{eq:lim_gQ}, we conclude that
\[
\left( \dashint_Q \left| \left(f_j-\dashint_{\Fth Q} f_j \,dm \right) - \left(g-\dashint_{\Fth Q} g \,dm\right) \right|^q \,dm \right)^{1/q} \leq \epsilon
\]
for all $\Fth Q \in\Adm$ and $Q \in \Ch(\Fth Q)$, and all $j \geq N$. As $\epsilon > 0$ was arbitrary, this proves that $\|f_j-g\|_{\BMO_q} \to 0$ as $j \to \infty$, i.e., the sequence $(f_n+\CC)_n$ converges to $g+\CC$ in $\BMO_q(X)/\CC$.
\end{proof}

We shall prove that the spaces $\BMO_q(X)$ for $q \in [1,\infty)$ are all the same, with equivalent seminorms. We first establish a density result for $L^\infty(X)$ in $\BMO_1(X)$.

\begin{prp}\label{prp:BMOlip}
Let $(X,\cE,m,\Adm)$ be a CZ space.
\begin{enumerate}[label=(\roman*)]
\item\label{en:BMOlip_lip} If $F : \CC \to \CC$ is an $L$-Lipschitz function and $f \in \BMO_1(X)$, then $F \circ f \in \BMO_1(X)$ and $\|F \circ f\|_{\BMO_1} \leq 2 L \|f\|_{\BMO_1}$.
\item\label{en:BMOlip_mod} If $f \in \BMO_1(X)$, then $|f| \in \BMO_1(X)$ and $\| |f| \|_{\BMO_1(X)} \leq 2 \| f \|_{\BMO_1(X)}$.
\item\label{en:BMOlip_approx}
If $f \in \BMO_1(X)$, then there exist functions $f_n \in L^\infty(X)$ such that
\[
|f_n| \leq |f|, \qquad \|f_n\|_\infty \leq n, \qquad \|f_n\|_{\BMO_1} \leq 2 \|f\|_{\BMO_1}
\]
for all $n \in \Npos$, and moreover $f_n \to f$ pointwise.
\end{enumerate}
\end{prp}
\begin{proof}
\ref{en:BMOlip_lip}. We first observe that, for all $\Fth Q \in \Adm$, $Q \in \Ch(\Fth Q)$ and $g \in L^1_\loc(X)$, by the triangle inequality
\[
\dashint_Q \left|g-\dashint_{\Fth Q} g \,dm\right| \,dm \leq \dashint_Q |g-c| \,dm + \dashint_{\Fth Q} |g-c| \,dm
\]
for all $c \in \CC$. If we apply this inequality with $g = F \circ f$ and $c = F\left(\dashint_{\Fth Q} f \,dm\right)$, then we obtain
\[\begin{split}
&\dashint_Q \left|F \circ f -\dashint_{\Fth Q} F \circ f \,dm\right| \,dm \\
&\leq \dashint_Q \left|F \circ f- F\left(\dashint_{\Fth Q} f \,dm\right)\right| \,dm 
+ \dashint_{\Fth Q} \left|F \circ f-F\left(\dashint_{\Fth Q} f \,dm\right)\right| \,dm\\
&\leq L \left[ \dashint_Q \left|f - \dashint_{\Fth Q} f \,dm\right| \,dm + \dashint_{\Fth Q} \left|f - \dashint_{\Fth Q} f \,dm\right| \,dm\right] \\
&\leq L \left[ \dashint_Q \left|f - \dashint_{\Fth Q} f \,dm\right| \,dm + \sum_{R \in \Ch(\Fth Q)} \frac{m(R)}{m(\Fth Q)} \dashint_{R} \left|f - \dashint_{\Fth Q} f \,dm\right| \,dm\right] \\
&\leq 2 L \|f\|_{\BMO_1},
\end{split}\]
as required.

\ref{en:BMOlip_mod}. This follows from part \ref{en:BMOlip_lip} and the fact that $z \mapsto |z|$ is $1$-Lipschitz.

\ref{en:BMOlip_approx}. For $n \in \Npos$, let $F_n : \CC \to \CC$ be defined by
\[
F_n(z) = \begin{cases}
z &\text{if } |z| \leq n,\\
n z/|z| &\text{if } |z| \geq n.
\end{cases}
\]
It is easily seen that $F_n$ is continuous and bounded, with $\|F_n\|_\infty = n$ and $|F_n(z)| \leq |z|$ for all $z \in \CC$; actually, as the (real) differential of $z \mapsto z/|z|$ is
\[
D (z/|z|) = |z|^{-1} (I - z z^T/|z|^2),
\]
one sees that $\|D F_n(z)\|_{\mathrm{op}} \leq 1$ whenever $|z| \neq n$, so $F_n$ is $1$-Lipschitz. As a consequence, from part \ref{en:BMOlip_lip} it follows that, if $f \in \BMO_1(X)$ and we set $f_n \defeq F_n \circ f$, then the $f_n$ have the required properties.
\end{proof}

We now prove a version of the John--Nirenberg inequality.

\begin{prp}\label{prp:BMOJN}
Let $(X,\cE,m,\Adm)$ be a CZ space.
\begin{enumerate}[label=(\roman*)]
\item\label{en:BMOJN_Linfty} If $f \in L^\infty(X)$ and $\|f\|_{\BMO_1} > 0$, then
\[
\sup_{\Fth Q \in \Adm}\sup_{Q \in \Ch(\Fth Q)} \dashint_Q \exp\left( \frac{t}{\|f\|_{\BMO_1}} \left| f - \dashint_{\Fth Q} f \,dm \right| \right) \,dm \lesssim_t 1
\]
for all $t \in (0,2^{-1} \log 2)$, with universal implicit constants.

\item\label{en:BMOJN_BMO} If $f \in \BMO_1(X)$ and $\|f\|_{\BMO_1} > 0$, then
\[
\sup_{\Fth Q \in \Adm}\sup_{Q \in \Ch(\Fth Q)} \dashint_Q \exp\left( \frac{t}{\|f\|_{\BMO_1}} \left| f - \dashint_{\Fth Q} f \,dm \right| \right) \,dm \lesssim_t 1
\]
for all $t \in (0,4^{-1} \log 2)$, with universal implicit constants.
\end{enumerate}
\end{prp}
\begin{proof}
\ref{en:BMOJN_Linfty}. For any $Q_0 \in \Adm$ and $R \in \Ch(Q_0)$, by definition we have
\[
\int_R \left| f - \dashint_{Q_0} f \,dm \right| \,dm \leq m(R) \|f\|_{\BMO_1}.
\]
We can now apply Lemma \ref{lem:stoppingsets} to the function $f_R = \left| f - \dashint_{Q_0} f \,dm \right| \chr_{R} \in L^1(X)$, the admissible set $R \in \Adm$, and the level $\lambda = 2\|f\|_{\BMO_1}$,
in order to obtain a disjoint collection $\cQ_R \subseteq \Ch^*(R) \setminus \{R\}$ such that:
\begin{itemize}
\item for any $Q \in \cQ_R$, there is $\Fth Q \in \Ch^*(R)$ with $Q \in \Ch(\Fth Q)$ and
\begin{equation}\label{eq:stopping_JN}
\left| \dashint_{\Fth Q} f \,dm - \dashint_{Q_0} f \,dm \right| \leq \dashint_{\Fth Q} \left| f - \dashint_{Q_0} f \,dm \right| \,dm \leq 2 \|f\|_{\BMO_1};
\end{equation}
\item if $Z_R \defeq \bigcup \cQ_R$, then
\begin{equation}\label{eq:halving_JN}
m(Z_R) \leq \frac{\left\| \left| f - \dashint_{Q_0} f \,dm \right| \chr_{R} \right\|_1}{2 \|f\|_{\BMO_1}} \leq \frac{m(R)}{2}
\end{equation}
and
\begin{equation}\label{eq:Linfty_JN}
\left\| \left| f - \dashint_{Q_0} f \,dm \right| \chr_{R\setminus Z_R} \right\|_\infty \leq 2 \|f\|_{\BMO_1}.
\end{equation}
\end{itemize}
Notice now that, for any $Q \in \cQ_R$,
\[
\left| \dashint_{\Fth Q} f \,dm - \dashint_{Q_0} f \,dm \right| \leq \dashint_{\Fth Q} \left| f - \dashint_{Q_0} f \,dm \right| \,dm \leq 2 \|f\|_{\BMO_1}.
\]
So, for all $t > 0$,
\[\begin{split}
&\dashint_{R} \exp\left( \frac{t}{\|f\|_{\BMO_1}} \left| f - \dashint_{Q_0} f \,dm \right| \right) \,dm 
= \frac{1}{m(R)} \left[\int_{R \setminus Z_R} + \sum_{Q \in \cQ_R} \int_Q \right] \\
&\leq \exp(2t) \left[ 1 + \sum_{Q \in \cQ_R} \frac{m(Q)}{m(R)} \dashint_Q \exp\left( \frac{t}{\|f\|_{\BMO_1}} \left| f - \dashint_{\Fth Q} f \,dm \right| \right) \,dm \right],
\end{split}\]
where we used \eqref{eq:stopping_JN} and \eqref{eq:Linfty_JN}; thus, by \eqref{eq:halving_JN},
\begin{equation}\label{eq:ind_ineq}
\begin{split}
&\dashint_{R} \exp\left( \frac{t}{\|f\|_{\BMO_1}} \left| f - \dashint_{Q_0} f \,dm \right| \right) \,dm \\
&\leq \exp(2t) \left[ 1 + \frac{m(Z_R)}{m(R)} \sup_{Q \in \cQ_R} \dashint_Q \exp\left( \frac{t}{\|f\|_{\BMO_1}} \left| f - \dashint_{\Fth Q} f \,dm \right| \right) \,dm \right] \\
&\leq \exp(2t) \left[ 1 + \frac{1}{2} \sup_{Q \in \cQ_R} \dashint_Q \exp\left( \frac{t}{\|f\|_{\BMO_1}} \left| f - \dashint_{\Fth Q} f \,dm \right| \right) \,dm \right].
\end{split}
\end{equation}

Thus, if we set
\[
F(t) \defeq \sup_{\Fth Q \in \Adm} \sup_{Q \in \Ch(\Fth Q)} \dashint_{Q} \exp\left( \frac{t}{\|f\|_{\BMO_1}} \left| f - \dashint_{\Fth Q} f \,dm \right| \right) \,dm,
\]
then, by taking the supremum over $Q_0 \in \Adm$ and $R \in \Ch(Q_0)$ in \eqref{eq:ind_ineq}, we obtain
\begin{equation}\label{eq:ind_ineq_sup}
F(t) \leq e^{2t} [ 1+ F(t)/2 ].
\end{equation}
Notice that, as $f \in L^\infty(X)$, we have
\[
F(t) \leq \exp\left( \frac{2t \|f\|_\infty}{\|f\|_{\BMO_1}} \right) < \infty
\]
for all $t > 0$. Thus, from \eqref{eq:ind_ineq_sup} we obtain that
\[
F(t) \leq e^{2t} (1-e^{2t}/2)^{-1}
\]
whenever the right-hand side is positive, i.e., whenever $t \in (0,2^{-1} \log 2)$.

\ref{en:BMOJN_BMO}. By Proposition \ref{prp:BMOlip}\ref{en:BMOlip_approx}, we can approximate $f \in \BMO_1$ by functions $f_n \in L^\infty(X)$ such that $\|f_n\|_{\BMO_1} \leq 2\|f\|_{\BMO_1}$ and $f_n \to f$ pointwise. By applying part \ref{en:BMOJN_Linfty} to $f_n$, we then obtain that, for all $\Fth Q \in \Adm$ and $Q \in \Ch(\Fth Q)$,
\[
\dashint_Q \exp\left( \frac{2t}{\|f\|_{\BMO_1}} \left| f_n - \dashint_{\Fth Q} f_n \,dm \right| \right) \,dm \lesssim_t 1
\]
for all $t \in (0,2^{-1} \log 2)$, with an implicit constant independent of $\Fth Q$ and $Q$. By Fatou's lemma, we thus conclude that
\[
\dashint_Q \exp\left( \frac{2t}{\|f\|_{\BMO_1}} \left| f - \dashint_{\Fth Q} f \,dm \right| \right) \,dm \lesssim_t 1
\]
too, and the desired estimate follows by relabelling $2t$ as $t$.
\end{proof}

We can now prove the equality between the spaces $\BMO_q(X)$.

\begin{cor}\label{cor:BMOequiv}
Let $(X,\cE,m,\Adm)$ be a CZ space. For all $q \in [1,\infty)$, we have $\BMO_q(X) = \BMO_1(X)$ and, for all $f \in \BMO_1(X)$,
\[
\|f\|_{\BMO_1} \leq \|f\|_{\BMO_q} \lesssim_q \|f\|_{\BMO_1}
\]
with universal implicit constants.
\end{cor}
\begin{proof}
From Proposition \ref{prp:BMO}\ref{en:BMO_incl}, we already know that $\BMO_q(X) \subseteq \BMO_1(X)$ and $\|f\|_{\BMO_1} \leq \|f\|_{\BMO_q}$.

Conversely, if $f \in \BMO_1(X)$, then, by Proposition \ref{prp:BMOJN}\ref{en:BMOJN_BMO} applied with $t = 8^{-1} \log 2$, we deduce that, for all $\Fth Q \in \Adm$ and $Q \in \Ch(\Fth Q)$,
\[
\frac{1}{\|f\|_{\BMO_1}^q} \dashint_Q \left| f - \dashint_{\Fth Q} f \,dm \right|^q \,dm \lesssim_q \dashint_Q \exp\left( \frac{t}{\|f\|_{\BMO_1}} \left| f - \dashint_{\Fth Q} f \,dm \right| \right) \,dm \lesssim 1,
\]
where we used the inequality $u^q \lesssim_q e^u$ for $u \in [0,\infty)$. Thus $f \in \BMO_q(X)$ and
\[
\|f\|_{\BMO_q} = \sup_{\Fth Q \in \Adm} \sup_{Q \in \Ch(\Fth Q)} \left(\dashint_Q \left| f - \dashint_{\Fth Q} f \,dm \right|^q \,dm \right)^{1/q} \lesssim_q \|f\|_{\BMO_1},
\]
as claimed.
\end{proof}

We finally discuss a property of $\BMO_q$-spaces which is valid when $m(X) < \infty$.

\begin{prp}\label{prp:BMOfin}
Let $(X,\cE,m,\Adm)$ be a CZ space with $m(X) < \infty$. Let $q \in [1,\infty)$.
\begin{enumerate}[label=(\roman*)]
\item\label{en:BMOfin_Lq} $\BMO_q(X) \subseteq L^q(X) \subseteq L^1(X)$.
\item\label{en:BMOfin_direct} The map $f \mapsto (f+\CC,\dashint_X f \,dm)$ is a linear isomorphism between $\BMO_q(X)$ and the direct product $(\BMO_q(X)/\CC) \times \CC$.
\item\label{en:BMOfin_norm} If we set
\[
\vvvert f \vvvert_{\BMO_q} \defeq \max\left\{ \|f\|_{\BMO_q}, \left| \dashint_X f \,dm \right| \right\}
\]
for all $f \in \BMO_q(X)$, then $\vvvert \cdot \vvvert_{\BMO_q}$ is a norm on $\BMO_q(X)$ which makes it into a Banach space.
\end{enumerate}
\end{prp}
\begin{proof}
\ref{en:BMOfin_Lq}. Assume first $q = 1$.
By Proposition \ref{prp:BMOlip}\ref{en:BMOlip_mod}, if $f \in\BMO_1(X)$, then $|f| \in \BMO_1(X)$ too. In particular, by Proposition \ref{prp:BMO},
\[
\sup_{Q \in \Adm} \dashint_Q \left| |f| - \dashint_Q |f| \,dm \right| \,dm \lesssim \| |f| \|_{\BMO_1}.
\]
Let $(Q_n)_n$ be an exhausting sequence of admissible sets. Then, for all $n \in \NN$,
\begin{multline*}
\frac{m(Q_0)}{m(X)}
\left| \dashint_{Q_0} |f| \,dm - \dashint_{Q_n} |f| \,dm \right|
 \leq \frac{m(Q_0)}{m(Q_n)} \dashint_{Q_0} \left| |f| - \dashint_{Q_n} |f| \,dm \right| \,dm \\
\leq \dashint_{Q_n} \left| |f| - \dashint_{Q_n} |f| \,dm \right| \,dm \lesssim \| |f| \|_{\BMO_1},
\end{multline*}
thus also
\[
\dashint_{Q_n} |f| \,dm \lesssim \frac{m(X)}{m(Q_0)} \| |f| \|_{\BMO_1} + \dashint_{Q_0} |f| \,dm.
\]
As the right-hand side is independent of $n$, by taking the limit as $n \to \infty$ one concludes that $f \in L^1(X)$.

Now, for an arbitrary $q \in [1,\infty)$, by Corollary \ref{cor:BMOequiv} we have 
\[
\BMO_q(X) = \BMO_1(X) \subseteq L^1(X).
\]
Moreover, again by Proposition \ref{prp:BMO}, for any $f \in \BMO_q(X)$,
\[
\sup_{Q \in \Adm} \left(  \dashint_Q \left| |f| - \dashint_Q |f| \,dm \right|^q \,dm\right)^{1/q} \lesssim \| |f| \|_{\BMO_q},
\]
and in particular
\[
\left(  \dashint_{Q_n} |f|^q \,dm \right)^{1/q} \lesssim_q \| |f| \|_{\BMO_q} + \dashint_{Q_n} |f| \,dm \leq \| |f| \|_{\BMO_q} + \frac{1}{m(Q_n)}\|f\|_1;
\]
again, by taking the limit as $n \to \infty$, we conclude that $f \in L^q(X)$ too.

\ref{en:BMOfin_direct}. By part \ref{en:BMOfin_Lq}, the map $f \mapsto (f + \CC,\dashint_X f \,dm)$ is well defined, and it is immediately checked that it is a linear isomorphism.

\ref{en:BMOfin_norm}. In light of part \ref{en:BMOfin_direct}, the quantity $\vvvert \cdot \vvvert_{\BMO_q}$ defined above corresponds, via the linear isomorphism, to the usual direct product norm on $(\BMO_q(X)/\CC) \times \CC$. As the direct product of Banach spaces is a Banach space, and $\BMO_q(X)/\CC$ is complete by Proposition \ref{prp:BMOcompl}, we conclude that $\BMO_q(X)$ with $\vvvert \cdot \vvvert_{\BMO_q}$ is Banach.
\end{proof}

\subsection{Hardy spaces}
We now introduce atomic Hardy spaces on a CZ space.

\begin{dfn}\label{dfn:H1}
Let $(X,\cE,m,\Adm)$ be a CZ space. Let $q \in (1,\infty]$.
\begin{enumerate}[label=(\alph*)]
\item A \emph{$q$-atom} is a function of the form
\[
a-\left(\int a \,dm \right) \frac{1}{m(\Fth Q)} \chr_{\Fth Q},
\]
for some $Q,\Fth Q \in \Adm$ with $Q \in \Ch(\Fth Q)$ and $a \in L^q(X)$ with $\{ a \neq 0 \} \subseteq Q$ and $\|a\|_q \leq m(Q)^{-1/q'}$.

\item\label{en:H1_H1} We define $H^{1,q}(X)$ as the set of functions $f \in L^1(X)$ that can be written as
\[
f = \sum_{n \in \NN} \lambda_n b_n,
\]
with convergence in $L^1(X)$, where the $b_n$ are $q$-atoms and $(\lambda_n)_n \in \ell^1$. For $f \in H^{1,q}(X)$, we define
\[
\|f\|_{H^{1,q}} = \inf \left\{ \sum_{n \in \NN} |\lambda_n| \tc f = \sum_{n \in \NN} \lambda_n b_n, \ b_n \text{ $q$-atoms}, \ (\lambda_n)_n \in \ell^1 \right\}.
\]
\end{enumerate}
\end{dfn}

\begin{prp}\label{prp:H1}
Let $(X,\cE,m,\Adm)$ be a CZ space. Let $q \in (1,\infty]$.
\begin{enumerate}[label=(\roman*)]
\item\label{en:H1_atom} If $b$ is a $q$-atom, then $\|b\|_1 \leq 2$ and $\int b \,dm = 0$; moreover, $b$ is also an $r$-atom for any $r \in (1,q]$.
\item\label{en:H1_L1} $H^{1,q}(X)$ is a vector subspace of $L^1(X)$ and $\|\cdot\|_{H^{1,q}}$ is a norm on $H^{1,q}(X)$. Moreover, $\|f\|_1 \leq 2 \|f\|_{H^{1,q}}$ for all $f \in H^{1,q}(X)$.
\item\label{en:H1_incl} For any $r \in (1,q]$, we have $H^{1,q}(X) \subseteq H^{1,r}(X)$ and $\|f\|_{H^{1,r}} \leq \|f\|_{H^{1,q}}$ for all $f \in H^{1,q}(X)$.
\item\label{en:H1_dense} For any $r \in [q,\infty]$, $H^{1,r}(X)$ is dense in $H^{1,q}(X)$.
\end{enumerate}
\end{prp}
\begin{proof}
Part \ref{en:H1_atom} is an immediate consequence of the definition of $q$-atom and H\"older's inequality. Parts \ref{en:H1_L1} and \ref{en:H1_incl} easily follow from part \ref{en:H1_atom} and the definitions.

As for part \ref{en:H1_dense}, due to the inclusions from part \ref{en:H1_incl}, it is enough to show that $H^{1,\infty}(X)$ is dense in $H^{1,q}(X)$. Indeed, let
\[
b = a - \left(\int a \,dm\right) \frac{1}{m(\Fth Q)}\chr_{\Fth Q}
\]
be a $q$-atom, where $Q,\Fth Q \in \Adm$, $Q \in \Ch(\Fth Q)$, $a \in L^q(X)$, $\|a\|_q \leq m(Q)^{-1/q'}$ and $\{a \neq 0\} \subseteq Q$. For all $n \in \Npos$, if $a_n \defeq a \chr_{\{ |a| \leq n\}}$, we have $a_n \in L^\infty(X)$ and $\{ a_n \neq 0\} \subseteq Q$, so
\[
b_n \defeq a_n - \left(\int a_n \,dm\right) \frac{1}{m(\Fth Q)}\chr_{\Fth Q}
\]
is a multiple of an $\infty$-atom, with $\|b_n\|_{H^{1,\infty}} \leq m(Q) n$. Moreover,
\[
b-b_n = (a-a_n) - \left(\int (a-a_n) \,dm\right) \frac{1}{m(\Fth Q)}\chr_{\Fth Q},
\]
so
\[
\|b-b_n\|_{H^{1,q}} \leq m(Q)^{1/q'} \|a-a_n\|_{q} \to 0 \qquad\text{as } n \to \infty.
\]
This shows that $b \in \overline{H^{1,\infty}(X)}^{H^{1,q}(X)}$. As $b$ is an arbitrary $q$-atom, we conclude that $H^{1,\infty}(X)$ is dense in $H^{1,q}(X)$.
\end{proof}

\begin{prp}\label{prp:H1complete}
Let $(X,\cE,m,\Adm)$ be a CZ space and $q \in (1,\infty]$. The space $H^{1,q}(X)$ with the norm $\|\cdot\|_{H^{1,q}}$ is a Banach space.
\end{prp}
\begin{proof}
Let $(f_n)_{n \in \NN}$ be a Cauchy sequence in $H^{1,q}(X)$. We can then find a subsequence $(f_{n_k})_{k \in \NN}$ with bounded variation; in other words, if $g_0 \defeq f_{n_0}$ and $g_{k+1} \defeq f_{n_{k+1}}-f_{n_{k}}$ for all $k \in \NN$, then $\sum_{k \in \NN} \|g_k\|_{H^{1,q}} < \infty$.

For each $k \in \NN$, since $g_k \in H^{1,q}(X)$, we can write $g_k = \sum_{j \in \NN} \lambda_{k,j} b_{k,j}$ for some $q$-atoms $b_{k,j}$ and scalars $\lambda_{k,j}$ with $\sum_{j \in \NN} |\lambda_{k,j}| \leq 2 \|g_k\|_{H^{1,q}}$. In particular, $\sum_{k,j \in \NN} |\lambda_{k,j}| \leq 2 \sum_{k \in \NN} \|g_k\|_{H^{1,q}} < \infty$. So, if we set 
\[
h \defeq \sum_{k \in \NN} g_k = \sum_{k \in \NN} \sum_{j \in \NN} \lambda_{k,j} b_{k,j},
\]
then the series converges absolutely in $L^1(X)$, and defines an element of $H^{1,q}(X)$. Moreover, for all $\ell \in \NN$,
\[
h-f_{n_\ell} = \sum_{k = \ell+1}^\infty g_k = \sum_{k=\ell+1}^\infty \sum_{j \in \NN} \lambda_{k,j} b_{k,j},
\]
thus
\[
\|h-f_{n_\ell}\|_{H^{1,q}} \leq \sum_{k=\ell+1}^\infty \sum_{j \in \NN} |\lambda_{k,j}| \leq 2 \sum_{k=\ell+1}^\infty \|g_k\|_{H^{1,q}},
\]
which tends to $0$ as $\ell \to \infty$.

This proves that the subsequence $(f_{n_k})_k$ tends to $h$ in $H^{1,q}(X)$. As $(f_n)_n$ is a Cauchy sequence, we conclude that the whole sequence $(f_n)_n$ tends to $h$.
\end{proof}

We now discuss an $L^q$-variant of the Calder\'on--Zygmund decomposition of Proposition \ref{prp:CZdec}, involving the Hardy space $H^{1,q}(X)$.

\begin{prp}\label{prp:H1CZ}
Let $(X,\cE,m,\Adm)$ be a CZ space. Let $q \in (1,\infty)$, $f \in L^q(X)$ and $\lambda > \|f\|_q/m(X)^{1/q}$. Then, there exists a decomposition
\[
f = g + b
\]
such that:
\begin{enumerate}[label=(\alph*)]
\item\label{en:H1CZ_good} $\|g\|_p \lesssim_{p,q} \lambda^{1-q/p} \|f\|_q^{q/p}$ for all $p \in [q,\infty)$, with universal implicit constants;
\item\label{en:H1CZ_bad} $\|b\|_{H^{1,q}} \leq \lambda^{1-q} \|f\|_q^q$.
\end{enumerate}
\end{prp}
\begin{proof}
Much as in the proof of Proposition \ref{prp:CZdec}, we may reduce to the case where $f$ is well-supported, which we shall assume in the rest of the proof.

As $\lambda > \|f\|_q/m(X)^{1/q}$ and $f$ is well-supported, we may find $Q_0 \in \Adm$ with $\{ f \neq 0 \} \subseteq Q_0$ and $\lambda > \|f\|_q/m(Q_0)^{1/q}$.

So we have $|f|^q \in L^1(X)$, $\{ |f|^q \neq 0 \} \subseteq Q_0$ and $\lambda^q > \| |f|^q \|_1/m(Q_0)$. 
By Lemma \ref{lem:stoppingsets}, applied to the function $|f|^q$ and the level $\lambda^q$ in place of $f$ and $\lambda$, we can find a disjoint family $\cQ \subseteq \Ch^*(Q_0) \setminus \{Q_0\}$ such that
\begin{equation}\label{eq:Qmeas_est}
m\left( \bigcup \cQ \right) \leq \| f \|_q^q /\lambda^q,
\end{equation}
and, if $\Omega \defeq X\setminus \bigcup \cQ$, then
\begin{equation}\label{eq:OmegaLinfty}
\| f \chr_\Omega \|_\infty = \| |f|^q \chr_{Q_0 \setminus \bigcup \cQ}\|_\infty^{1/q} \leq \lambda;
\end{equation}
moreover, for any $Q \in \cQ$, there is $\Fth Q \in \Ch^*(Q_0)$ such that $Q \in \Ch(\Fth Q)$ and
\begin{equation}\label{eq:stopping_Lq}
\dashint_{\Fth Q} |f|^q \,dm \leq \lambda^q.
\end{equation}

We now define, for any $Q \in \cQ$,
\[
b_Q \defeq f \chr_Q - \left(\int_Q f \,dm\right) \frac{1}{m(\Fth Q)} \chr_{\Fth Q}.
\]
As $f \in L^q(X)$, clearly $b_Q$ is a multiple of a $q$-atom and
\[
\|b_Q\|_{H^{1,q}} \leq m(Q)^{1/q'} \|f \chr_Q\|_q,
\]
so, by H\"older's inequality and the disjointness of $\cQ$,
\begin{multline}\label{eq:H1q_est}
\sum_{Q \in \cQ} \|b_Q\|_{H^{1,q}} \leq \sum_{Q \in \cQ} m(Q)^{1/q'} \|f \chr_Q\|_q \leq m\left(\bigcup \cQ\right)^{1/q'} \|f\|_q \\
\leq \lambda^{-q/q'} \| f \|_q^{q/q'+1} = \lambda^{1-q} \|f\|_q^q,
\end{multline}
where we used \eqref{eq:Qmeas_est} and the fact that $q/q' = q-1$. This shows that, if we set $b \defeq \sum_{Q \in \cQ} b_Q$, then $b \in H^{1,q}(X)$ and the estimate \ref{en:H1CZ_bad} holds true.

We thus have the decomposition $f = g+b$, provided we set
\[
g \defeq f-b = f \chr_\Omega + \sum_{Q \in \cQ} \left( \int_Q f \,dm \right) \frac{1}{m(\Fth Q)} \chr_{\Fth Q},
\]
where $\Omega \defeq X \setminus \bigcup \cQ = \{ \cM_{Q_0} |f|^q \leq \lambda^q \}$, and it remains to show that $g$ satisfies the estimates \ref{en:H1CZ_good}.

Actually, from \eqref{eq:OmegaLinfty} we deduce that,
for any $p \in [q,\infty)$, by H\"older's inequality,
\[
\|f \chr_\Omega\|_p \leq \|f\|_q^{q/p} \|f \chr_\Omega\|_\infty^{1-q/p} \leq \lambda^{1-q/p} \|f\|_q^{q/p}.
\]
Therefore, if we set
\[
h \defeq \sum_{Q \in \cQ} \left( \int_Q |f| \,dm \right) \frac{1}{m(\Fth Q)} \chr_{\Fth Q},
\]
then $|g| \leq |f| \chr_\Omega + h$, and the proof of \ref{en:H1CZ_good} is reduced to that of analogous estimates for $h$ in place of $g$.

Notice now that, for any $Q \in \cQ$, by \eqref{eq:stopping_Lq} and H\"older's inequality,
\[
\dashint_{\Fth Q} |f| \,dm \leq \left( \dashint_{\Fth Q} |f|^q \,dm \right)^{1/q} \leq \lambda.
\]
Therefore, by Lemma \ref{lem:higherint}, H\"older's inequality and \eqref{eq:Qmeas_est}, it follows that
\[
\|h\|_p^p \lesssim_p \lambda^{p-1} \|f \chr_{\bigcup \cQ} \|_1 \leq \lambda^{p-1} m\left(\bigcup \cQ\right)^{1/q'} \|f\|_q \leq \lambda^{p-q} \|f\|_q^q,
\]
as required.
\end{proof}

\begin{rem}
Much as in the proof of Theorem \ref{thm:mainCZdec}, an iterative argument would allow one to upgrade property \ref{en:H1CZ_good} in Proposition \ref{prp:H1CZ} to an $L^\infty$-bound $\|g\|_\infty \lesssim_q \lambda$. We shall recover this improved decomposition from the interpolation properties of the $H^{1,q}$-spaces, as shown in Corollary \ref{cor:H1Linfty} below.
\end{rem}

We briefly discuss a few auxiliary spaces that will be useful in the development of the theory.

\begin{dfn}
Let $(X,\cE,m,\Adm)$ be a CZ space.
\begin{enumerate}[label=(\alph*)]
\item For any $q \in [1,\infty]$ and $E \in \cE$ with $0 < m(E) < \infty$, we define
\[
L^q_0(E) \defeq \left\{ f \in L^q(X) \tc \{ f \neq 0 \} \subseteq E, \ \int f \,dm = 0 \right\}.
\]
\item For any $q \in [1,\infty]$, we denote by $L^q_\comp(X)$ the space of functions $f \in L^q(X)$ such that the set $\{ f \neq 0\}$ is contained in a finite union of admissible sets.
\item For any $q \in (1,\infty]$, we denote by $H^{1,q}_\fin(X)$ the subspace of $H^{1,q}(X)$ made of finite linear combinations of $q$-atoms.
\end{enumerate}
\end{dfn}

\begin{prp}\label{prp:Lq0H1}
Let $(X,\cE,m,\Adm)$ be a CZ space and $q \in [1,\infty]$.
\begin{enumerate}[label=(\roman*)]
\item\label{en:Lq0H1_clos} Let $E,F \in \cE$ with $0 < m(E),m(F) < \infty$. Then $L^q_0(E)$ is a closed subspace of $L^q(X)$, and $L^q_0(E) \cap L^q_0(F) = L^q_0(E \cap F)$.
\end{enumerate}
Assume now that $q \in (1,\infty]$.
\begin{enumerate}[label=(\roman*),resume]
\item\label{en:Lq0H1_Hfin} $H^{1,q}_\fin(X)$ is a dense subspace of $H^{1,q}(X)$, contained in $L^q_\comp(X)$.
\item\label{en:Lq0H1_incl} Let $\Fth Q \in \Adm$. Then $L^q_0(\Fth Q)$ is contained in $H^{1,q}(X)$, and
\[
\|f\|_{H^{1,q}} \leq m(\Fth Q)^{1/q'} \|f\|_q
\]
for all $f \in L^q_0(\Fth Q)$.
\end{enumerate}
\end{prp}
\begin{proof}
\ref{en:Lq0H1_clos}.
This is clear from the definitions.

\ref{en:Lq0H1_Hfin}. This is clear from the definitions.

\ref{en:Lq0H1_incl}.
If $f \in L^q_0(\Fth Q)$, then $\{ f \neq 0\} \subseteq \Fth Q$ and $\int f \,dm = 0$, so we can write
\[
f = \sum_{Q \in \Ch(\Fth Q)} f \chr_{Q} = \sum_{Q \in \Ch(\Fth Q)} b_Q, \qquad b_Q \defeq f \chr_{Q} - \left(\int_Q f \,dm\right) \frac{1}{m(\Fth Q)} \chr_{\Fth Q},
\]
where the sums over $Q \in \Ch(\Fth Q)$ converge in $L^1(X)$.
As $f \in L^q(X)$, each $b_Q$ is a multiple of a $q$-atom, thus $b_Q \in H^{1,q}(X)$. Moreover, by H\"older's inequality,
\[
\sum_{Q \in \Ch(\Fth Q)} \|b_Q\|_{H^{1,q}} \leq \sum_{Q \in \Ch(\Fth Q)} m(Q)^{1/q'} \|f\chr_{Q}\|_q \leq m(\Fth Q)^{1/q'} \|f\|_q,
\]
thus the sum $\sum_{Q \in \Ch(\Fth Q)} b_Q$ converges in $H^{1,q}(X)$ too; this proves that $f \in H^{1,q}(X)$ and that the desired norm bound holds.
\end{proof}

\begin{rem}\label{rem:regularHardy}
For $q \in (1,\infty]$, let us call ``ordinary $q$-atom'' any function $v \in L^q(X)$ such that $\{ v \neq 0 \} \subseteq Q$ for some $Q \in \Adm$, $\int_X v \,dm = 0$ and $\|v\|_q \leq m(Q)^{-1/q'}$. By replacing ``$q$-atom'' with ``ordinary $q$-atom'' in Definition \ref{dfn:H1}\ref{en:H1_H1}, one obtains another space $H^{1,q}_\ord(X)$ with norm $\|\cdot\|_{H^{1,q}_\ord}$, which corresponds to the definition of the atomic Hardy spaces elsewhere in the literature (cf., e.g., \cite[Section 2]{V}).

Notice that, if $v$ is an ordinary $q$-atom, then, by Proposition \ref{prp:Lq0H1}\ref{en:Lq0H1_incl}, $v \in H^{1,q}(X)$ and $\|v\|_{H^{1,q}} \leq 1$. From this it is easily seen that $H^{1,q}_\ord(X) \subseteq H^{1,q}(X)$ and $\|f\|_{H^{1,q}} \leq \|f\|_{H^{1,q}_\ord}$ for all $f \in H^{1,q}_\ord(X)$.

Conversely, let $b = a - (\int a \,dm) \chr_{\Fth Q}/m(\Fth Q)$ be a $q$-atom, where $\Fth Q \in \Adm$, $Q \in \Ch(\Fth Q)$, $a \in L^q(X)$, $\{ a \neq 0 \} \subseteq Q$ and $\|a\|_q \leq m(Q)^{-1/q'}$. Then $\|b\|_q \leq  2m(Q)^{-1/q'} < \infty$ and $\{ b \neq 0\} \subseteq \Fth Q$, so $b$ is a multiple of an ordinary $q$-atom, and therefore
\begin{equation}\label{eq:ordatom}
\|b\|_{H^{1,q}_\ord} \leq m(\Fth Q)^{1/q'} \|b\|_q \leq 2 (m(\Fth Q)/m(Q))^{1/q'}.
\end{equation}
This shows that $H^{1,q}_\fin(X) \subseteq H^{1,q}_\ord(X)$, but in general there may not be a uniform inequality of norms, as $m(\Fth Q)/m(Q)$ need not be uniformly bounded, and the inclusion $H^{1,q}_\ord(X) \subseteq H^{1,q}(X)$ need not be an equality.

However, if we assume that the CZ space $(X,\cE,m,\Adm)$ is regular, as in Definition \ref{dfn:regularCZ}, then the above inequality \eqref{eq:ordatom} shows that $\|b\|_{H^{1,q}_\ord} \leq 2 \rho(X)^{1/q'}$ for all $q$-atoms $b$. So in this case we may conclude that $H^{1,q}(X) = H^{1,q}_\ord(X)$ and
\[
\|f\|_{H^{1,q}} \leq \|f\|_{H^{1,q}_\ord} \leq 2 \rho(X)^{1/q'} \|f\|_{H^{1,q}}
\]
for all $f \in H^{1,q}(X)$.
\end{rem}

We now introduce a variant of the spaces $H^{1,q}(X)$, which is specific to the case $m(X)<\infty$.

\begin{dfn}
Let $(X,\cE,m,\Adm)$ be a CZ space with $m(X) < \infty$, and let $q \in (1,\infty]$.
\begin{enumerate}[label=(\alph*)]
\item We define 
\[
\lH^{1,q}(X) \defeq \{ c + f \tc c \in \CC, \ f \in H^{1,q}(X) \},
\]
i.e., $\lH^{1,q}(X)$ is the linear span of $H^{1,q}(X)$ and constant functions.
\item For any $f \in \lH^{1,q}(X)$, we set
\[
\|f\|_{\lH^{1,q}} \defeq \left| \int_X f \,dm \right| + \left\| f - \dashint_X f \,dm \right\|_{H^{1,q}}.
\]
\end{enumerate}
\end{dfn}

\begin{prp}\label{prp:lH1}
Let $(X,\cE,m,\Adm)$ be a CZ space with $m(X) < \infty$. Let $q \in (1,\infty]$.
\begin{enumerate}[label=(\roman*)]
\item\label{en:lH1_L10} $H^{1,q}(X)$ is contained in $L^1_0(X)$.
\item\label{en:lH1_direct} The map $f \mapsto (f - \dashint_X f \,dm, \int_X f \,dm)$ defines a linear isomorphism between $\lH^{1,q}(X)$ and the direct sum $H^{1,q}(X) \oplus \CC$.
\item\label{en:lH1_norm} $\|\cdot\|_{\lH^{1,q}}$ is a norm on $\lH^{1,q}(X)$, which restricts to $\|\cdot\|_{H^{1,q}}$ on $H^{1,q}(X)$, and $H^{1,q}(X)$ is a closed subspace of $\lH^{1,q}(X)$.
\item\label{en:lH1_banach} $\lH^{1,q}(X)$ with the norm $\|\cdot\|_{\lH^{1,q}}$ is a Banach space.
\item\label{en:lH1_L1} $\lH^{1,q}(X) \subseteq L^1(X)$ and $\|f\|_1 \leq 2 \|f\|_{\lH^{1,q}}$ for all $f \in \lH^{1,q}(X)$.
\item\label{en:lH1_Lqincl} If 
$X \in \Adm$ or 
$q \in (1,\infty)$, then $L^q(X) \subseteq \lH^{1,q}(X)$ and $\|f\|_{\lH^{1,q}} \leq 3 m(X)^{1/q'} \|f\|_q$ for all $f \in L^q(X)$.
\item\label{en:lH1_CZ} Assume that $q \in (1,\infty)$. Let $f \in L^q(X)$ and $\lambda > 0$. Then, there exists a decomposition
\[
f = g + b
\]
such that:
\begin{enumerate}[label=(\alph*)]
\item\label{en:lH1_CZ_good} $\|g\|_p \lesssim_{p,q} \lambda^{1-q/p} \|f\|_q^{q/p}$ for all $p \in [q,\infty)$, with universal implicit constants;
\item\label{en:lH1_CZ_bad} $\|b\|_{\lH^{1,q}} \leq 3 \lambda^{1-q} \|f\|_q^q$.
\end{enumerate}
\end{enumerate}
\end{prp}
\begin{proof}
\ref{en:lH1_L10}. Clearly $q$-atoms are in $L^1_0(X)$; as the latter is a closed subspace of $L^1(X)$, we conclude that $H^{1,q}(X)$ is contained in $L^1_0(X)$ too.

\ref{en:lH1_direct}. If $u \defeq \chr_X / m(X)$, then by definition $\lH^{1,q}(X) = H^{1,q}(X) + \CC u$, and from part \ref{en:lH1_L10} it follows that $H^{1,q}(X) \cap \CC u = \{0\}$, so the sum is direct. Moreover, if $f = f_0+c u \in H^{1,q}(X) \oplus \CC u$, then $\int_X f \,dm = c$ and $f_0 = f - cu = f - \dashint_X f \,dm$, which proves that indeed $f \mapsto (f-\dashint_X f \,dm,\int_X f \,dm)$ is a linear isomorphism.

\ref{en:lH1_norm}. In light of part \ref{en:lH1_direct}, the quantity $\|\cdot\|_{\lH^{1,q}}$ corresponds via the linear isomorphism to the usual norm for the ($\ell^1$-)direct sum of normed spaces.

\ref{en:lH1_banach}. The direct sum of Banach spaces is also a Banach space, and $H^{1,q}(X)$ is Banach by Proposition \ref{prp:H1complete}.

\ref{en:lH1_L1}. By Proposition \ref{prp:H1}\ref{en:H1_L1}, for all $f \in \lH^{1,q}(X)$,
\begin{multline*}
\|f\|_1 \leq \left\|f - \dashint_X f \,dm \right\|_1 + m(X) \left|\dashint_X f \,dm \right| \\
\leq 2 \left\|f - \dashint_X f \,dm \right\|_{H^{1,q}} +  \left|\int_X f \,dm \right| \leq 2 \|f\|_{\lH^{1,q}}.
\end{multline*}

\ref{en:lH1_Lqincl}. 
Let $(Q_n)_n$ be an exhausting sequence. Thus, $L^q_0(Q_n) \subseteq L^q_0(Q_{n+1})$ for all $n \in \NN$, and therefore $V \defeq \bigcup_{n \in \NN} L^q_0(Q_n)$ is a vector subspace of $L^q_0(X)$. Moreover, by Proposition \ref{prp:Lq0H1}\ref{en:Lq0H1_incl}, $V \subseteq H^{1,q}(X)$ and
\[
\|f\|_{H^{1,q}} \leq m(X)^{1/q'} \|f\|_{q} \qquad\text{for all } f \in V.
\]

In particular, $V \oplus \CC u \subseteq H^{1,q}(X) \oplus \CC u = \lH^{1,q}(X)$. Moreover, for any $f = f_0 + cu \in V \oplus \CC u$, we have
\[
|c| \leq \|f\|_1 \leq m(X)^{1/q'} \|f\|_q, \quad \|f_0\|_q \leq \|f\|_q + |c| \|u\|_q \leq 2 \|f\|_q,
\]
and therefore
\[
\|f\|_{\lH^{1,q}} = \|f_0\|_{H^{1,q}} + |c| \leq 3 m(X)^{1/q'} \|f\|_q.
\]
This shows that the desired norm inequality holds for the elements of the subspace $V \oplus \CC u$ of $L^q(X)$. In the case $X \in \Adm$, we may take the trivial exhausting sequence given by $Q_n = X$ for all $n$, thus $V \oplus \CC u = L^q_0(X) \oplus \CC u = L^q(X)$, and we are done.

If instead $X \notin \Adm$ and $q<\infty$, in order to conclude, it is enough to show that $V \oplus \CC u$ is dense in $L^q(X)$: indeed, by density, the inclusion map $V \oplus \CC u \to \lH^{1,q}$ would have a unique bounded extension $T : L^q(X) \to \lH^{1,q}(X)$ (since $\lH^{1,q}(X)$ is complete by part \ref{en:lH1_banach}), but such $T$ would a fortiori be bounded from $L^q(X)$ to $L^1(X)$ (by part \ref{en:lH1_L1}), and therefore $T$ would coincide with the inclusion map $L^q(X) \to L^1(X)$.

To prove that $V \oplus \CC u$ is dense in $L^q(X)$, take $f \in L^q(X)$ and observe that, for all $n \in \NN$, we can write
\begin{equation}\label{eq:f_fn_cnu_rn_dec}
f = f_n + c_n u + r_n,
\end{equation}
where
\[
c_n \defeq \int_{Q_n} f \,dm, \qquad f_n \defeq f \chr_{Q_n} - c_n u_n, \qquad u_n \defeq \chr_{Q_n}/m(Q_n).
\]
Clearly $f_n \in L^q_0(Q_n)$, so $f_n + c_n u \in V \oplus \CC u$. Moreover,
\begin{equation}\label{eq:rn_tends_to_zero}
r_n = f - f_n - c_n u = f \chr_{X \setminus Q_n} + c_n (u_n-u) \to 0
\end{equation}
in $L^q(X)$ as $n \to \infty$, because $(Q_n)_n$ is an exhausting sequence and $q<\infty$.

\ref{en:lH1_CZ}. In the case $\lambda > \|f\|_q/m(X)^{1/q}$, Proposition \ref{prp:H1CZ} provides the required decomposition $f = g+b$, with $b \in H^{1,q}(X)$. Assume instead that $\lambda \leq \|f\|_q/m(X)^{1/q}$. Then, 
$\lambda^{1-q} \geq  m(X)^{1/q'} \|f\|_q^{1-q}$, 
so, by part \ref{en:lH1_Lqincl},
\[
\|f\|_{\lH^{1,q}} \leq 3 m(X)^{1/q'} \|f\|_q \leq 3 \lambda^{1-q} \|f\|_q^q
\]
and therefore the required decomposition is obtained by taking $g = 0$ and $b = f$.
\end{proof}

\begin{rem}
We will see in Corollary \ref{cor:lH1Lqincl} that the inclusion $L^q(X) \subseteq \lH^{1,q}(X)$ and the norm inequality stated in Proposition \ref{prp:lH1}\ref{en:lH1_Lqincl} hold in greater generality.
\end{rem}

We shall later prove that the spaces $H^{1,q}(X)$ do not depend on $q$. For this purpose, we shall exploit the duality between Hardy and BMO spaces, together with the already proved independence of $\BMO_q(X)$ on $q$.

\subsection{\texorpdfstring{$H^1$}{H1}--\texorpdfstring{$\BMO$}{BMO} duality}

We now study the duality between the Hardy and BMO spaces that we introduced.

\begin{prp}\label{prp:H1BMOprelim}
Let $(X,\cE,m,\Adm)$ be a CZ space and $q \in (1,\infty]$.
For any $g \in \BMO_{q'}(X)$, the functional $\phi_g : H^{1,q}_\fin(X) \to \CC$, defined by
\begin{equation}\label{eq:H1finBMO_pairing}
\phi_g(f) \defeq \int_X f \, g \,dm \qquad \forall f \in H^{1,q}_\fin(X),
\end{equation}
is linear and satisfies
\[
|\phi_g(f)| \leq \|f\|_{H^{1,q}} \|g\|_{\BMO_{q'}}
\]
for all $f \in H^{1,q}_\fin(X)$. Thus, it uniquely extends to a bounded linear functional $\phi_g : H^{1,q}(X) \to \CC$ with
\[
\|\phi_g\|_{H^{1,q}(X)'} \leq \|g\|_{\BMO_{q'}} \ .
\]
\end{prp}
\begin{proof}
As $\BMO_{q'}(X) \subseteq L^{q'}_\loc(X)$ and $H^{1,q}_\fin(X) \subseteq L^q_\comp(X)$, the pairing \eqref{eq:H1finBMO_pairing} is well defined and determines a linear functional $\phi_g : H^{1,q}_\fin(X) \to \CC$.

In addition, for any $q$-atom
\[
b = a - \left(\int a \,dm \right) \frac{1}{m(\Fth Q)} \chr_{\Fth Q},
\]
where $Q,\Fth Q \in \Adm$ satisfy $Q \in \Ch(\Fth Q)$, $\{ a \neq 0 \} \subseteq Q$ and $\|a\|_q \leq m(Q)^{-1/q'}$, we have
\[\begin{split}
|\phi_g(b)| 
&= \left| \int a \, g \,dm - \left( \int a \,dm \right) \dashint_{\Fth Q} g \,dm \right| 
= \left| \int a \left( g - \dashint_{\Fth Q} g \,dm \right) \,dm \right| \\
&\leq \|a\|_q \left( \int_Q \left| g - \dashint_{\Fth Q} g \,dm \right|^{q'} dm \right)^{1/q'} 
\leq \left( \dashint_Q \left| g - \dashint_{\Fth Q} g \,dm \right|^{q'} dm \right)^{1/q'} \\
&\leq \|g\|_{\BMO_{q'}},
\end{split}\]
where we used H\"older's inequality. We have thus proved the bound
\begin{equation}\label{eq:unif_atom_bd}
|\phi_g(b)| \leq \|g\|_{\BMO_{q'}} \text{ for all $q$-atoms $b$}
\end{equation}
for any $g \in \BMO_{q'}$.

In the case $g \in L^\infty(X)$, the pairing \eqref{eq:H1finBMO_pairing} actually makes sense for all $f \in L^1(X)$, and defines a continuous linear functional $\phi_g : L^1(X) \to \CC$; hence a fortiori $\phi_g : H^{1,q}(X) \to \CC$ is continuous and from the uniform bound \eqref{eq:unif_atom_bd}
we easily deduce that $\|\phi_g\|_{H^{1,q}(X)'} \leq \|g\|_{\BMO_{q'}}$.

Now, an arbitrary $g \in \BMO_{q'}(X) = \BMO_1(X)$ can be approximated by functions $g_n \in L^\infty(X)$ as in Proposition \ref{prp:BMOlip}\ref{en:BMOlip_approx}, so that $|g_n| \leq |g|$ and $g_n \to g$ pointwise, while $\sup_n \|g_n\|_{\BMO_{q'}} \lesssim_q \|g\|_{\BMO_{q'}} < \infty$. In particular, the functionals $\phi_{g_n} \in H^{1,q}(X)'$ are uniformly bounded, and moreover by dominated convergence one sees that $\phi_{g_n}(f) \to \phi_g(f)$ for all $f \in H^{1,q}_\fin(X)$. Thus 
\[
|\phi_g(f)| \leq \sup_n |\phi_{g_n}(f)| \leq \|f\|_{H^{1,q}} \sup_n \|g_n\|_{\BMO_{q'}}
\]
for all $f \in H^{1,q}_\fin(X)$, which shows that $\phi_g : H^{1,q}_\fin(X) \to \CC$ is indeed continuous with respect to the $H^{1,q}$-norm, and therefore it extends uniquely to a bounded linear functional $\phi_g \in H^{1,q}(X)'$. We can then apply again the uniform bound \eqref{eq:unif_atom_bd} to conclude that $\|\phi_g\|_{H^{1,q}(X)'} \leq \|g\|_{\BMO_{q'}}$.
\end{proof}

We would like to show that any bounded linear functional on $H^{1,q}(X)$ is of the form $\phi_g$ for some $g \in \BMO_{q'}(X)$. For this purpose, we first prove a dual characterisation of the space $\BMO_{q'}(X)$, effectively reversing the norm inequality proved in Proposition \ref{prp:H1BMOprelim}.

\begin{lem}\label{lem:BMOnormdual}
Let $(X,\cE,m,\Adm)$ be a CZ space and $q \in (1,\infty]$. Let $g \in L^{q'}_\loc(X)$ and $C \in [0,\infty)$ be such that
\[
\left| \int_X g \, b \,dm \right| \leq C \quad\text{ for all $q$-atoms $b$.}
\]
Then, $g \in \BMO_{q'}(X)$ and $\|g\|_{\BMO_{q'}} \leq C$.
\end{lem}
\begin{proof}
For all $\Fth Q \in \Adm$ and $Q \in \Ch(\Fth Q)$, by the $L^q$--$L^{q'}$ duality,
\begin{equation}\label{eq:Lqdualnorm}
\left( \int_Q \left| g - \dashint_{\Fth Q} g \,dm \right|^{q'} \,dm \right)^{1/q'} = \sup_{f \in L^q(X) \tc \|f\|_q \leq 1}\left| \int_Q \left( g - \dashint_{\Fth Q} g \,dm \right) f \, dm \right|.
\end{equation}
Thus, by the $L^1$--$L^\infty$ duality applied to weighted spaces of sequences $(\alpha_Q)_{Q \in \Ch(\Fth Q)}$ indexed by $\Ch(\Fth Q)$, we obtain that
\[\begin{split}
&\sup_{Q \in \Ch(\Fth Q)} \left( \dashint_Q \left| g - \dashint_{\Fth Q} g \,dm \right|^{q'} \,dm \right)^{1/q'} \\
&= \sup_{\|f\|_q \leq 1} \sup_{Q \in \Ch(\Fth Q)} \frac{1}{m(Q)^{1/q'}} \left| \int_Q \left( g - \dashint_{\Fth Q} g \,dm \right) f \, dm \right| \\
&= \sup_{\|f\|_q \leq 1} \sup_{\sum\limits_{Q \in \Ch(\Fth Q)} m(Q)^{1/q'} |\alpha_Q| \leq 1} \left| \sum_{Q \in \Ch(\Fth Q)} \alpha_Q \int_Q \left( g - \dashint_{\Fth Q} g \,dm \right) f \, dm \right| \\
&= \sup_{\|f\|_q \leq 1} \sup_{\sum\limits_{Q \in \Ch(\Fth Q)} m(Q)^{1/q'} |\alpha_Q| \leq 1} \left| \int_{\Fth Q} \left( g - \dashint_{\Fth Q} g \,dm \right) \left( f \sum_{Q \in \Ch(\Fth Q)} \alpha_Q \chr_Q\right) \, dm \right| .
\end{split}\]
As $g - \dashint_{\Fth Q} g \,dm$ has zero integral on $\Fth Q$, we deduce that
\[\begin{split}
&\sup_{Q \in \Ch(\Fth Q)} \left( \dashint_Q \left| g - \dashint_{\Fth Q} g \,dm \right|^{q'} \,dm \right)^{1/q'} \\
&= \sup_{\|f\|_q \leq 1} \sup_{\sum\limits_{Q \in \Ch(\Fth Q)} m(Q)^{1/q'} |\alpha_Q| \leq 1} \Biggl| \int_{\Fth Q} \left( g - \dashint_{\Fth Q} g \,dm \right) \\
&\qquad \times \left( f \sum_{Q \in \Ch(\Fth Q)} \alpha_Q \chr_Q - \sum_{Q \in \Ch(\Fth Q)} \alpha_Q \left( \int_Q f \,dm \right) \frac{1}{m(\Fth Q)} \chr_{\Fth Q} \right) \, dm \Biggr|.
\end{split}\]
Now, notice that
\begin{multline*}
f \sum_{Q \in \Ch(\Fth Q)} \alpha_Q \chr_Q - \sum_{Q \in \Ch(\Fth Q)} \alpha_Q \left( \int_Q f \,dm \right) \frac{1}{m(\Fth Q)} \chr_{\Fth Q} \\
= \sum_{Q \in \Ch(\Fth Q)} \alpha_Q \left( f \chr_Q - \left( \int_Q f \,dm \right) \frac{1}{m(\Fth Q)} \chr_{\Fth Q} \right)
\end{multline*}
has zero integral and vanishes outside $\Fth Q$, so
\[\begin{split}
&\sup_{Q \in \Ch(\Fth Q)} \left( \dashint_Q \left| g - \dashint_{\Fth Q} g \,dm \right|^{q'} \,dm \right)^{1/q'} \\
&= \sup_{\|f\|_q \leq 1} \sup_{\sum\limits_{Q \in \Ch(\Fth Q)} m(Q)^{1/q'} |\alpha_Q| \leq 1} \left| \int g 
 \sum_{Q \in \Ch(\Fth Q)} \alpha_Q \left( f \chr_Q - \left( \int_Q f \,dm \right) \frac{1}{m(\Fth Q)} \chr_{\Fth Q} \right) \, dm \right| \\
&= \sup_{\|f\|_q \leq 1} \sup_{\sum\limits_{Q \in \Ch(\Fth Q)} m(Q)^{1/q'} |\alpha_Q| \leq 1} \left| \sum_{Q \in \Ch(\Fth Q)} \alpha_Q \int g 
 \left( f \chr_Q - \left( \int_Q f \,dm \right) \frac{1}{m(\Fth Q)} \chr_{\Fth Q} \right) \, dm \right| .
\end{split}\]
As $f \in L^q(X)$, the functions $f \chr_Q - \left( \int_Q f \,dm \right) \frac{1}{m(\Fth Q)} \chr_{\Fth Q}$ are multiples of $q$-atoms, so, by our assumption on $g$,
\[\begin{split}
&\sup_{Q \in \Ch(\Fth Q)} \left( \dashint_Q \left| g - \dashint_{\Fth Q} g \,dm \right|^{q'} \,dm \right)^{1/q'} \\
&\leq \sup_{\|f\|_q \leq 1} \sup_{\sum\limits_{Q \in \Ch(\Fth Q)} m(Q)^{1/q'} |\alpha_Q| \leq 1} \sum_{Q \in \Ch(\Fth Q)} |\alpha_Q| \, C \, m(Q)^{1/q'} \|f \chr_{Q}\|_q \\
&\leq C.
\end{split}\]
By taking the supremum over $\Fth Q \in \Adm$, we finally deduce that $g \in \BMO_{q'}(X)$ and 
$\|g\|_{\BMO_{q'}} \leq C$, as required.
\end{proof}

\begin{prp}\label{prp:H1BMOdual}
Let $(X,\cE,m,\Adm)$ be a CZ space and $q \in (1,\infty]$.
The correspondence $g \mapsto \phi_g$ from Proposition \ref{prp:H1BMOprelim} defines an isometric linear map from $\BMO_{q'}(X)/\CC$ to $H^{1,q}(X)'$, which is surjective whenever $q < \infty$.
\end{prp}
\begin{proof}
By Proposition \ref{prp:H1BMOprelim}, we deduce that, for all $q \in (1,\infty]$, the correspondence $g \mapsto \phi_g$ determines a contractive linear map from $\BMO_{q'}(X)/\CC$ to $H^{1,q}(X)'$, which is actually an isometry by Lemma \ref{lem:BMOnormdual}. It only remains to show that the map is surjective in the case $q < \infty$.

Let $\phi \in H^{1,q}(X)'$. For any $Q \in \Adm$, from the continuous inclusion $L^q_0(Q) \subseteq H^{1,q}(X)$ proved in Proposition \ref{prp:Lq0H1}, we see that $\phi|_{L^q_0(Q)} : L^q_0(Q) \to \CC$ is a bounded linear functional, with
\[
|\phi(f)| \leq \|\phi\|_{H^{1,q}(X)'} \, m(Q)^{1/q'} \|f\|_q \qquad \forall f \in L^q_0(Q).
\]
Thus, there exists $g_{Q} \in L^{q'}(Q)$ such that
\begin{equation}\label{eq:gQdef}
\phi(f) = \int f \, g_Q \,dm \qquad \forall f \in L^q_0(Q);
\end{equation}
notice that $g_Q$ is not uniquely determined by the latter property, since it can be replaced by $g_Q + c$ for any constant $c \in \CC$.

Now, if $Q,R \in \Adm$ have non-negligible intersection, then
\[
\int f \, (g_Q - g_R) \, dm = \phi(f)-\phi(f) = 0 \qquad\forall f \in L^q_0(Q \cap R) = L^q_0(Q) \cap L^q_0(R).
\]
This implies that $(g_Q - g_R)|_{Q \cap R}$ is constant.

We can now apply Lemma \ref{lem:glueing} to construct a function $g \in L^{q'}_\loc(X)$ such that $(g-g_Q)|_Q$ is constant for all $Q \in \Adm$. So, from \eqref{eq:gQdef} we deduce that
\begin{equation}\label{eq:funct_Lq0}
\phi(f) = \int f \, g \,dm \qquad\forall f \in L^q_0(Q) \ \forall Q \in \Adm.
\end{equation}
In particular, for any $q$-atom $b$,
\[
\left| \int g \, b \, dm \right| = |\phi(b)| \leq \|\phi\|_{H^{1,q}(X)'}.
\]
Thus, from Lemma \ref{lem:BMOnormdual}, we deduce that $g \in \BMO_{q'}(X)$ and
\begin{equation}\label{eq:gnorm_ineq}
\|g\|_{\BMO_{q'}} \leq \|\phi\|_{H^{1,q}(X)'}.
\end{equation}

In particular, by Proposition \ref{prp:H1BMOprelim} we can consider the functional $\phi_g \in H^{1,q}(X)'$ associated with $g$. Now, as each $q$-atom is in $L^q_0(Q)$ for some $Q \in \Adm$, from \eqref{eq:funct_Lq0} we see that $\phi$ coincides with $\phi_g$ on $q$-atoms, hence on the whole $H^{1,q}(X)$ by linearity and continuity, so $\phi = \phi_g$.
\end{proof}

The proof of the surjectivity in Proposition \ref{prp:H1BMOdual} breaks down in the case $q=\infty$: indeed, it is not clear whether we can construct functions $g_Q \in L^1(Q)$ satisfying \eqref{eq:gQdef}, as the dual of $L^\infty(Q)$ is in general larger than $L^1(Q)$. 

Nevertheless, as we shall prove, one has surjectivity in the case $q=\infty$ too. It is convenient first to consider a particular class of CZ spaces.

\begin{dfn}
We say that a CZ space $(X,\cE,m,\Adm)$ is \emph{simple} if $X \in \Adm$ and $\Adm = \Ch^*(X)$.
\end{dfn}

It is also convenient to work with an (a priori) smaller space than $H^{1,\infty}(X)$, defined by using a subclass of $\infty$-atoms.

\begin{dfn}
Let $(X,\cE,m,\Adm)$ be a simple CZ space.
\begin{enumerate}[label=(\alph*)]
\item\label{en:dspec_f} We say that a function $f \in L^0(X)$ is \emph{special} if there exists a finite subcollection $\cS \subseteq \Ch^*(X)$ such that $f$ is measurable with respect to the $\sigma$-algebra on $X$ generated by $\cS$.
\item A \emph{special $\infty$-atom} is an $\infty$-atom which is special in the sense of part \ref{en:dspec_f}.
\item The atomic space $\widetilde{H}^{1,\infty}(X)$, with norm $\|\cdot\|_{\widetilde{H}^{1,\infty}}$, is defined exactly as $H^{1,\infty}(X)$ in Definition \ref{dfn:H1}\ref{en:H1_H1}, but replacing $\infty$-atoms with special $\infty$-atoms.
\end{enumerate}
\end{dfn}

\begin{prp}\label{prp:sH1}
Let $(X,\cE,m,\Adm)$ be a simple CZ space.
\begin{enumerate}[label=(\roman*)]
\item\label{en:sH1_H1} $\widetilde{H}^{1,\infty}(X) \subseteq H^{1,\infty}(X)$ and $\|f\|_{H^{1,\infty}} \leq \|f\|_{\widetilde{H}^{1,\infty}}$ for all $f \in \widetilde{H}^{1,\infty}(X)$.
\item\label{en:sH1_banach} $\widetilde{H}^{1,\infty}(X)$ is a Banach space with the norm $\|\cdot\|_{\widetilde{H}^{1,\infty}}$.
\item\label{en:sH1_Linfty0} Let $Q \in \Adm$. If $f \in L^\infty_0(Q)$ is special, then $f \in \widetilde{H}^{1,\infty}(X)$ and $\|f\|_{\widetilde{H}^{1,\infty}} \leq m(Q) \|f\|_\infty$.
\item\label{en:sH1_dense} $\widetilde{H}^{1,\infty}(X)$ is dense in $H^{1,q}(X)$ for any $q \in (1,\infty)$.
\item\label{en:sH1_dual} Let $f \in L^1(X)$. Then
\[
\|f\|_{\BMO_1} \leq \sup \left\{ \left| \int f \, b \,dm \right| \tc b \text{ special $\infty$-atom} \right\}.
\]
\end{enumerate}
\end{prp}
\begin{proof}
\ref{en:sH1_H1}. This follows from the fact that any special $\infty$-atom is an $\infty$-atom.

\ref{en:sH1_banach}. This is proved in the same way as Proposition \ref{prp:H1complete}.

\ref{en:sH1_Linfty0}. The proof follows the lines of that of Proposition \ref{prp:Lq0H1}\ref{en:Lq0H1_incl}.

\ref{en:sH1_dense}. We already know that $\widetilde{H}^{1,\infty}(X) \subseteq H^{1,\infty}(X) \subseteq H^{1,q}(X)$. Now, if $f$ is a $q$-atom, and we define
\[
f_n \defeq \sum_{Q \in \Ch^n(X)} \left(\dashint_Q f \,dm\right) \chr_Q,
\]
then $f_n \in L^\infty_0(X)$ and $\|f_n\|_q \leq \|f\|_q$ for all $n \in \NN$, and moreover $f_n \to f$ in $L^q(X)$ \cite[Corollary 9.23]{K}. 
Further, for all $n \in \NN$, we can find a finite subfamily $\cS_n \subseteq \Ch^n(X)$ so that $\|f \chr_{X \setminus \bigcup \cS_n}\|_q < 2^{-n-1}$; thus, if we set
\begin{equation}\label{eq:martingale_approx_mod}
\tilde f_n \defeq \sum_{Q \in \cS_n} \left(\dashint_Q f \,dm\right) \chr_Q + \left(\dashint_{X \setminus \bigcup \cS_n} f \,dm \right) \chr_{X \setminus \bigcup \cS_n},
\end{equation}
then $\tilde f_n$ is special, and moreover $\|f_n - \tilde f_n\|_q \leq 2\|f \chr_{X \setminus \bigcup \cS_n}\|_q < 2^{-n}$, so $\tilde f_n \to f$ in $L^q(X)$ too. From part \ref{en:sH1_Linfty0} we then deduce that each $\tilde f_n$ is in $\widetilde{H}^{1,\infty}(X)$, while from Proposition \ref{prp:Lq0H1}\ref{en:Lq0H1_incl} we see that $\|f-\tilde f_n\|_{H^{1,q}} \leq m(X)^{1/q'} \|f-\tilde f_n\|_q \to 0$.

\ref{en:sH1_dual}. The proof is a minor modification of that of Lemma \ref{lem:BMOnormdual}. The key observation is the fact that, when applying the $L^\infty$--$L^1$ duality in \eqref{eq:Lqdualnorm}, instead of taking the supremum over all $f \in L^\infty(X)$ with $\|f\|_\infty \leq 1$, one can restrict to the functions $f$ which are special. Indeed, for an arbitrary $f \in L^\infty(X)$, the functions $\tilde f_n$ defined in \eqref{eq:martingale_approx_mod} converge to $f$ in $L^q(X)$ for some $q \in (1,\infty)$, hence (up to extracting a subsequence) also a.e., and moreover $\|\tilde f_n\|_\infty \leq \|f\|_\infty \leq 1$; so one can use dominated convergence to show that, for any $g \in L^1(X)$,
\[
\left| \int_X g f \,dm \right| = \lim_{n \to \infty} \left| \int_X g \tilde f_n \,dm \right| \leq \sup_n \left| \int_X g \tilde f_n \,dm \right|
\]
and therefore
\[
\|g\|_1 = \sup_{\|f\|_\infty \leq 1} \left| \int_X g f \,dm \right| = \sup_{\substack{\|f\|_\infty \leq 1 \\ f \text{ special}}} \left| \int_X g f \,dm \right|.
\]
Once $f$ is assumed to be special, the functions $f \chr_Q - (\int_Q f \,dm) \chr_{\Fth Q}/m(\Fth Q)$ constructed in the proof of Lemma \ref{lem:BMOnormdual} become multiples of special $\infty$-atoms, and the result follows.
\end{proof}

We can now adapt the proof of Proposition \ref{prp:H1BMOdual} in order to characterise the dual of $\widetilde{H}^{1,\infty}(X)$. The restriction to special $\infty$-atoms allows us to work around the fact that the dual of $L^\infty(X)$ is a larger space than $L^1(X)$ in general; in essence, we shall instead use the fact that the dual of the space $c_0$ of sequences vanishing at infinity is the space $\ell^1$ of summable sequences.

\begin{prp}\label{prp:sH1duality}
Let $(X,\cE,m,\Adm)$ be a simple CZ space.
The correspondence $g \mapsto \phi_g|_{\widetilde{H}^{1,\infty}(X)}$ defines an isometric isomorphism from $\BMO_1(X)/\CC$ to $\widetilde{H}^{1,\infty}(X)'$.
\end{prp}
\begin{proof}
From Propositions \ref{prp:H1BMOprelim} and \ref{prp:sH1}\ref{en:sH1_H1}, we deduce that $g \mapsto \phi_g|_{\widetilde{H}^{1,\infty}(X)}$ defines a contractive linear map from $\BMO_1(X)/\CC$ to $\widetilde{H}^{1,\infty}(X)'$, which is actually an isometry by Proposition \ref{prp:sH1}\ref{en:sH1_dual}. It remains to show that this isometry is surjective. 

Take $\phi \in \widetilde{H}^{1,\infty}(X)'$. We shall show that there exists $g \in \BMO_1(X)$ with $\phi = \phi_g|_{\widetilde{H}^{1,\infty}(X)}$.

Let $\Sigma_n$ be the $\sigma$-algebra generated by $\Ch^n(X)$. As $\Ch^n(X)$ is at most countable, 
if we denote by $\overline{W_n}$ the closure in $L^\infty(X,\Sigma_n,m)$ of $W_n \defeq \Span \{ \chr_Q \tc Q \in \Ch^n(X) \}$, then the dual of $\overline{W_n}$ is isometrically isomorphic to $L^1(X,\Sigma_n,m)$. 

Moreover, let $P : L^\infty(X) \to L^\infty_0(X)$ be the projection defined by $P f \defeq f - \dashint_X f \,dm$, and notice that $\|Pf\|_\infty \leq 2 \|f\|_\infty$. If $f \in W_n$, then $Pf \in L^\infty_0(X,\Sigma_n,m)$ is special, so by Proposition \ref{prp:sH1}\ref{en:sH1_Linfty0} we have $Pf \in \widetilde{H}^{1,\infty}(X)$ and
\[
|\phi(Pf)| \leq 2 m(X) \|\phi\|_{\widetilde{H}^{1,\infty}(X)'} \|f\|_\infty.
\]
This shows that the mapping $f \mapsto \phi(Pf)$ extends to a bounded linear functional on $\overline{W_n}$, and therefore there exists $g_n \in L^1(X,\Sigma_n,m)$ such that
\[
\|g_n\|_1 \leq 2 m(X) \|\phi\|_{\widetilde{H}^{1,\infty}(X)'}
\]
and
\begin{equation}\label{eq:phirepn}
\int_X f \, g_n \,dm = \phi(Pf)
\end{equation}
for all $f \in W_n$.

Notice now that, for all $\Fth Q \in \Ch^*(X)$ and $Q \in \Ch(\Fth Q)$, we have $b_Q \defeq \chr_Q - \frac{m(Q)}{m(\Fth Q)} \chr_{\Fth Q} \in \widetilde{H}^{1,\infty}(X)$ and $\|b_Q\|_{\widetilde{H}^{1,\infty}} \leq m(Q)$. In particular, the sum $\sum_{Q \in \Ch(\Fth Q)} b_Q$ converges absolutely in $\widetilde{H}^{1,\infty}(X)$, and therefore $\sum_{Q \in \Ch(\Fth Q)} b_Q = 0$ in $\widetilde{H}^{1,\infty}(X)$. As $\phi \in \widetilde{H}^{1,\infty}(X)'$, we deduce that $\sum_{Q \in \Ch(\Fth Q)} \phi(b_Q) = 0$, i.e.,
\[
\sum_{Q \in \Ch(\Fth Q)} \phi(P \chr_Q) = \phi(P \chr_{\Fth Q}).
\]
This information, together with the additivity of integrals, allows one to extend the validity of the identity \eqref{eq:phirepn} to all $f \in W_n^* \defeq \Span \{ \chr_Q \tc Q \in \bigcup_{k \leq n} \Ch^k(X) \}$. We remark that any special function $f \in L^\infty(X)$ is in $W_n^*$ for any sufficiently large $n$.

In particular, as $\chr_X \in W_n^*$ and $P\chr_X = 0$, we deduce that $\int_X g_n \,dm = 0$ for all $n \in \NN$.
Moreover, for any $Q \in \Ch^n(X)$, we have $\chr_Q \in W_n \subseteq W_N^*$ for all $N \geq n$, so
\[
\int_Q g_N \,dm = \phi(P \chr_Q) = \int_Q g_n \,dm.
\]
Thus, $(g_n)_n$ is a martingale with respect to the filtration $(\Sigma_n)_n$.

Observe further that, if we set
\[
\pi_n f \defeq \sum_{Q \in \Ch^n(X)} \left( \dashint_Q f \,dm \right) \chr_Q,
\]
then $\pi_n \chr_R \in W_n^*$ for all $R \in \Ch^*(X)$: indeed, if $R \in \Ch^k(X)$ for some $k \leq n$, then $\pi_n \chr_R = \chr_R$; if instead $R \in \Ch^k(X)$ for some $k > n$, then the above sum with $f = \chr_R$ has at most one nonzero summand. Therefore $\pi_n f \in W_n^*$ for all special functions $f \in L^\infty(X)$. This shows that
\[
\int_X f g_n \,dm = \int_X (\pi_n f) g_n \,dm = \phi(P \pi_n f)
\]
for all special functions $f \in L^\infty(X)$.

Now let $b$ be a special $\infty$-atom. Then $b = a - (\int a \,dm) \chr_{\Fth Q}/m(\Fth Q)$ for some $\Fth Q \in \Ch^N(X)$, $N \in \NN$, $Q \in \Ch(\Fth Q)$, $a \in L^\infty(X)$ special, $\{ a \neq 0 \} \subseteq Q$, $\|a\|_\infty \leq m(Q)^{-1}$.
Now, if $N \geq n$, then $\pi_n b = 0$, because $\{ b \neq 0 \} \subseteq \Fth Q$ and $\int b \,dm = 0$, and therefore
\[
\int_X b \, g_n \,dm = \phi(P \pi_n b) = 0.
\]
If instead $N < n$, then $\pi_n a$ is also a special function with $\{ \pi_n a \neq 0 \} \subseteq Q$ and $\int \pi_n a \,dm = \int a \, dm$, so $\pi_n b = \pi_n a - (\int a \,dm) \chr_{\Fth Q}/m(\Fth Q)$ is a special $\infty$-atom too, and
\[
\int_X b \, g_n \,dm = \phi(P \pi_n b) = \phi(\pi_n b).
\]
As $\phi$ is uniformly bounded on special $\infty$-atoms, from Proposition \ref{prp:sH1}\ref{en:sH1_dual} we deduce that $g_n \in \BMO_1(X)$ and
\[
\|g_n\|_{\BMO_1} \leq \|\phi\|_{\widetilde{H}^{1,\infty}(X)'}
\]
for all $n \in \NN$. From Corollary \ref{cor:BMOequiv} and Proposition \ref{prp:BMOfin}, we then also deduce that, for all $q \in (1,\infty)$, $g_n \in \BMO_{q'}(X) \cap L^1_0(X,\Sigma_n,m) \subseteq L^{q'}_0(X,\Sigma_n,m)$ and
\[
\|g_n\|_{L^{q'}} \lesssim_q \|g_n\|_{\BMO_{q'}} \lesssim_q \|\phi\|_{\widetilde{H}^{1,\infty}(X)'}.
\]

In conclusion, we have a martingale $(g_n)_n$ which is $L^{q'}$-bounded for any $q \in (1,\infty]$; from \cite[Corollary 9.23]{K} we conclude that $g_n \to g$ in $L^{q'}(X)$ for any $q \in (1,\infty]$, where $g \in \bigcap_{q \in (1,\infty]} L^{q'}(X)$, and moreover $g$ closes the martingale $(g_n)_n$.

Now, for any special $\infty$-atom $b$, we have that $b \in W^*_n$ for some $n \in \NN$, so in particular $b$ is $\Sigma_n$-measurable, and therefore
\[
\int_X b \, g \,dm = \int_X b \, g_n \,dm = \phi(P b) = \phi(b).
\]
Another application of Proposition \ref{prp:sH1}\ref{en:sH1_dual} then shows that $g \in \BMO_1(X)$, with
\[
\|g\|_{\BMO_1} \leq \|\phi\|_{\widetilde{H}^{1,\infty}(X)'}.
\]

In particular, we can consider the functional $\phi_g|_{\widetilde{H}^{1,\infty}} \in \widetilde{H}^{1,\infty}(X)'$. As clearly $\phi$ and $\phi_g$ coincide on special $\infty$-atoms, we conclude that $\phi = \phi_g|_{\widetilde{H}^{1,\infty}(X)}$.
\end{proof}

\begin{cor}\label{cor:sH1equality}
Let $X$ be a simple CZ space. Then $\widetilde{H}^{1,\infty}(X) = H^{1,\infty}(X)$, with equality of norms. Moreover, $H^{1,\infty}(X) = H^{1,q}(X)$ for all $q \in (1,\infty)$, with equivalence of norms (and universal implicit constants).
\end{cor}
\begin{proof}
From Propositions \ref{prp:H1}\ref{en:H1_incl} and \ref{prp:sH1}\ref{en:sH1_H1}, 
we already know that $\widetilde{H}^{1,\infty}(X) \subseteq H^{1,\infty}(X) \subseteq H^{1,q}(X)$ with continuous (contractive) inclusions, and moreover, by Proposition \ref{prp:sH1}\ref{en:sH1_dense}, $\widetilde{H}^{1,\infty}(X)$ is dense in $H^{1,q}(X)$. Now, for any $f \in \widetilde{H}^{1,\infty}(X)$, by Proposition \ref{prp:sH1duality} and the Hahn--Banach theorem,
\[
\|f\|_{\widetilde{H}^{1,\infty}} = \max \{ |\phi_g(f)| \tc g \in \BMO_{1}(X), \ \|g\|_{\BMO_{1}} \leq 1 \}.
\]
On the other hand, by Corollary \ref{cor:BMOequiv}, we know that $\BMO_{q'}(X) = \BMO_{1}(X)$, with equivalence of norms (and universal implicit constants), thus
\[
\|f\|_{\widetilde{H}^{1,\infty}} \simeq_q \max \{ |\phi_g(f)| \tc g \in \BMO_{q'}(X), \ \|g\|_{\BMO_{q'}} \leq 1 \} = \|f\|_{H^{1,q}},
\]
where the last equality follows from Proposition \ref{prp:H1BMOdual} and the Hahn--Banach theorem.

This shows that indeed the norms $\|\cdot\|_{\widetilde{H}^{1,\infty}}$ and $\|\cdot\|_{H^{1,q}}$ are equivalent (with universal implicit constants) on $\widetilde{H}^{1,\infty}(X)$. As $\widetilde{H}^{1,\infty}(X)$ is complete (see Proposition \ref{prp:sH1}\ref{en:sH1_banach}), this implies that $\widetilde{H}^{1,\infty}(X)$ is closed in $H^{1,q}(X)$, and therefore, by density, $\widetilde{H}^{1,\infty}(X) = H^{1,q}(X)$.

Thus, the inclusions $\widetilde{H}^{1,\infty}(X) \subseteq H^{1,\infty}(X) \subseteq H^{1,q}(X)$ must all be equalities, and the respective norms are equivalent to one another. As we already know from Proposition \ref{prp:H1BMOdual} that $g \mapsto \phi_g$ gives an isometric isomorphism $\BMO_{q'}(X)/\CC \to H^{1,q}(X)'$, and also an isometric embedding $\BMO_{1}(X)/\CC \to H^{1,\infty}(X)'$, we conclude that the latter is an isometric isomorphism too. From the dual characterisation of the norms on $H^{1,\infty}(X)$ and $\widetilde{H}^{1,\infty}(X)$, we finally conclude that $H^{1,\infty}(X) = \widetilde{H}^{1,\infty}(X)$ isometrically.
\end{proof}

We now extend the result to arbitrary CZ spaces (not necessarily simple).

\begin{prp}\label{prp:H1equiv}
Let $(X,\cE,m,\Adm)$ be a CZ space. Let $q \in (1,\infty]$.
\begin{enumerate}[label=(\roman*)]
\item\label{en:H1equiv_H1} We have $H^{1,\infty}(X) = H^{1,q}(X)$ with equivalence of norms (and universal implicit constants). 
\item\label{en:H1equiv_lH1} If $m(X) < \infty$, we also have $\lH^{1,\infty}(X) = \lH^{1,q}(X)$ with equivalence of norms (and universal implicit constants).
\item\label{en:H1equiv_H1dual} The map $g \mapsto \phi_g$ induces an isometric isomorphism from $\BMO_{q'}(X)/\CC$ to $H^{1,q}(X)'$.
\end{enumerate}
\end{prp}
\begin{proof}
\ref{en:H1equiv_H1}. From Proposition \ref{prp:H1}\ref{en:H1_incl}, we already know that $H^{1,\infty}(X) \subseteq H^{1,q}(X)$, with continuous (contractive) embedding. To prove the equality, it is enough to show that, for any $q$-atom $b$ on $X$, we have $b \in H^{1,\infty}(X)$ and
\[
\sup \{ \|b\|_{H^{1,\infty}(X)} \tc b \text{ $q$-atom}\} \lesssim_q 1.
\]

Let $b = a - (\int a \,dm) \chr_{\Fth Q}/m(\Fth Q)$ be a $q$-atom, where $\Fth Q \in \Adm$, $Q \in \Ch(\Fth Q)$, $a \in L^q(X)$, $\{ a \neq 0 \} \subseteq Q$, $\|a\|_q \leq m(Q)^{-1/q'}$. Notice now that $\Fth Q$ itself can be equipped with the structure of a simple CZ space, with admissible sets $\Ch^*(\Fth Q)$. If we identify functions on $\Fth Q$ with their extensions by zero to $X$, then $b$ is a $q$-atom on $\Fth Q$ too; in particular, $b \in H^{1,q}(\Fth Q) = H^{1,\infty}(\Fth Q)$ and
\[
\|b\|_{H^{1,\infty}(\Fth Q)} \lesssim_q \|b\|_{H^{1,q}(\Fth Q)} \leq 1,
\]
by Corollary \ref{cor:sH1equality}.
On the other hand, under the same identification, any $\infty$-atom on $\Fth Q$ is also an $\infty$-atom on $X$, so $H^{1,\infty}(\Fth Q) \subseteq H^{1,\infty}(X)$ with continuous (contractive) embedding, and in particular
\[
\|b\|_{H^{1,\infty}(X)} \leq \|b\|_{H^{1,\infty}(\Fth Q)} \lesssim_q 1,
\]
as desired.

\ref{en:H1equiv_lH1}. In light of Proposition \ref{prp:lH1}, this follows immediately from part \ref{en:H1equiv_H1}. 

\ref{en:H1equiv_H1dual}. By Proposition \ref{prp:H1BMOdual}, we know that the map $g \mapsto \phi_g$ gives an isometric embedding $\BMO_{q'}(X)/\CC \to H^{1,q}(X)'$, which is an isomorphism for $q < \infty$. As $\BMO_1(X) = \BMO_{q'}(X)$ and $H^{1,\infty}(X) = H^{1,q}(X)$, the isometric embedding is surjective for $q=\infty$ too.
\end{proof}

We can now establish the inclusion from Proposition \ref{prp:lH1}\ref{en:lH1_Lqincl} in full generality.

\begin{cor}\label{cor:lH1Lqincl}
Let $(X,\cE,m,\Adm)$ be CZ space with $m(X) < \infty$. Let $q \in (1,\infty]$. Then $L^q(X) \subseteq \lH^{1,q}(X)$ and 
\begin{equation}\label{eq:normineq_lH1q_Lq}
\|f\|_{\lH^{1,q}} \leq 3 m(X)^{1/q'}\|f\|_q
\end{equation}
for all $f \in L^q(X)$.
\end{cor}
\begin{proof}
From Proposition \ref{prp:lH1}\ref{en:lH1_Lqincl} we already know that $L^q(X) \subseteq \lH^{1,q}(X)$ for all $q<\infty$, with the norm inequality \eqref{eq:normineq_lH1q_Lq}; since $L^\infty(X) \subseteq L^q(X)$, by Proposition \ref{prp:H1equiv}\ref{en:H1equiv_lH1} we conclude that $L^\infty(X) \subseteq \lH^{1,q}(X) = \lH^{1,\infty}(X)$ too.

It only remains to prove the inequality \eqref{eq:normineq_lH1q_Lq} in the case $q=\infty$. Recall that, in the proof of Proposition \ref{prp:lH1}\ref{en:lH1_Lqincl}, this inequality was already established for all $f \in V \oplus \CC u$, where $u = \chr_X / m(X)$ and $V = \bigcup_{n \in \NN} L^\infty_0(Q_n)$ for a given exhausting sequence $(Q_n)_n$.

Moreover, much as in \eqref{eq:f_fn_cnu_rn_dec}, any $f \in L^\infty(X)$ can be decomposed as $f = g_n + r_n$, where
\[
g_n = (f \chr_{Q_n} - c_n u_n) + c_n u, \quad c_n \defeq \int_{Q_n} f \,dm, \quad u_n \defeq \chr_{Q_n}/m(Q_n).
\]
Thus $g_n \in V \oplus \CC u$, and moreover, arguing as in \eqref{eq:rn_tends_to_zero}, one sees that $r_n = f\chr_{X \setminus Q_n} + c_n (u_n-u) \to 0$ in $L^q(X)$ for all $q<\infty$. In particular, $g_n \to f$ in $L^q(X)$ for any $q<\infty$, hence also in $\lH^{1,q}(X) = \lH^{1,\infty}(X)$, and therefore
\[
\|f\|_{\lH^{1,\infty}} = \lim_{n \to \infty} \|g_n\|_{\lH^{1,\infty}} \leq 3 m(X) \liminf_{n \to \infty} \|g_n\|_\infty,
\]
where we applied to the $g_n$ the inequality \eqref{eq:normineq_lH1q_Lq} for $q=\infty$.
On the other hand, as $g_n = f \chr_{Q_n} + c_n(u-u_n)$, we see that
\[
|g_n(x)| \leq \begin{cases}
|f(x)| + \|f\|_1 (1/m(Q_n)-1/m(X)) &\text{if } x \in Q_n,\\
\|f\|_1/m(X) \leq \|f\|_\infty &\text{if } x \notin Q_n,
\end{cases}
\]
and therefore
\[
\liminf_{n \to \infty} \|g_n\|_\infty \leq \lim_{n \to \infty} [\|f\|_\infty + \|f\|_1 (1/m(Q_n)-1/m(X))] = \|f\|_\infty,
\]
as desired.
\end{proof}

Recall that, in the case $m(X) < \infty$, $\BMO_{q'}(X)$ is a Banach space with the norm $\vvvert \cdot \vvvert_{\BMO_{q'}}$ from Proposition \ref{prp:BMOfin}; in that case, we can also characterize $\BMO_{q'}(X)$ as a dual space.

\begin{prp}\label{prp:fH1BMOduality}
Let $(X,\cE,m,\Adm)$ be a CZ space with $m(X) < \infty$. Let $q \in (1,\infty]$.
For any $g \in \BMO_{q'}(X)$, the functional $\phi_g \in H^{1,q}(X)'$ extends uniquely to a functional $\phi_g \in \lH^{1,q}(X)'$ satisfying $\phi_g(\chr_X) = \int_X g \,dm$. Moreover, the correspondence $g \mapsto \phi_g$ defines an isometric isomorphism from $\BMO_{q'}(X)$ to $\lH^{1,q}(X)'$.
\end{prp}
\begin{proof}
Due to the isometric isomorphism $\lH^{1,q}(X) \cong H^{1,q}(X) \oplus \CC$ from Proposition \ref{prp:lH1}, we see that
\[
\lH^{1,q}(X)' \cong H^{1,q}(X)' \times \CC;
\]
so, by Propositions \ref{prp:H1equiv} and \ref{prp:BMOfin}, we also have
\[
\lH^{1,q}(X)' \cong (\BMO_{q'}(X)/\CC) \times \CC \cong \BMO_{q'}(X)
\]
isometrically.
From the explicit expressions for these isomorphisms, it is not difficult to check that the functional $\phi_g \in \lH^{1,q}(X)'$ corresponding to $g \in \BMO_{q'}(X)$ extends the functional $\phi_g \in H^{1,q}(X)'$ from Proposition \ref{prp:H1BMOprelim}, and moreover, if $u \defeq \chr_X/m(X)$, then $\phi_g(u) = \dashint_X g \,dm$, so indeed $\phi_g(\chr_X) = \int_X g\,dm$ as claimed.
\end{proof}

Thanks to the independence of the spaces $H^{1,q}(X)$ and $\BMO_{q'}(X)$ from the index $q \in (1,\infty]$ established in Corollary \ref{cor:BMOequiv} and Proposition \ref{prp:H1equiv}, it is convenient to drop the index from our notation.

\begin{dfn}
Let $(X,\cE,m,\Adm)$ be a CZ space.
We shall write $H^1(X)$ and $\BMO(X)$ in place of $H^{1,\infty}(X)$ and $\BMO_1(X)$. Moreover, when $m(X) < \infty$, we shall write $\lH^1(X)$ for $\lH^{1,\infty}(X)$.
\end{dfn}

\subsection{Real interpolation between \texorpdfstring{$H^1$}{H1}, \texorpdfstring{$\BMO$}{BMO}, and \texorpdfstring{$L^p$}{Lp}}

We refer to \cite[Chapter 3]{BL} for the theory of real interpolation, and particularly for the definition of the $K$-functional and the real interpolation functors $[\cdot,\cdot]_{\theta,p}$.

\begin{lem}\label{lem:interprelim}
Let $(X,\cE,m,\Adm)$ be a CZ space. Assume that $1 < q < p < \infty$, and let $\theta = p'/q'$, so that
\[
\frac{1-\theta}{1} + \frac{\theta}{p} = \frac{1}{q}.
\]
Then:
\begin{enumerate}[label=(\roman*)]
\item\label{en:interprelim_inf} if $m(X) = \infty$, we have the continuous embedding
\[
L^q(X) \subseteq [H^{1,q}(X),L^p(X)]_{\theta,\infty},
\]
with universal implicit constants;
\item\label{en:interprelim_fin} if $m(X) < \infty$, we have the continuous embedding
\[
L^q(X) \subseteq [\lH^{1,q}(X),L^p(X)]_{\theta,\infty},
\]
with universal implicit constants.
\end{enumerate}
\end{lem}
\begin{proof}
Notice that, if $\theta = p'/q'$, then indeed
\[
\frac{1-\theta}{1} + \frac{\theta}{p} = 1-\frac{p'}{q'}+\frac{p'}{p q'} = 1-\frac{1}{q'} = \frac{1}{q}.
\]

We first prove part \ref{en:interprelim_inf}.
Let $f \in L^q(X)$. By Proposition \ref{prp:H1CZ}, for any $\lambda > 0$ we can decompose $f = g^\lambda + b^\lambda$, where
\[
\|b^\lambda\|_{H^{1,q}} \leq \lambda^{1-q} \|f\|_q^q, \qquad \|g^\lambda\|_{p} \lesssim_{p,q} \lambda^{1-q/p} \|f\|_q^{q/p}.
\]
Therefore, for any $t \in (0,\infty)$,
\[\begin{split}
&K(t,f,H^{1,q}(X),L^p(X)) \\
&= \inf \{ \|b\|_{H^{1,q}} + t \|g\|_p \tc f = b+g, \ b \in H^{1,q}(X), \ g \in L^p(X)\} \\
&\lesssim_{p,q} \|f\|_q^q \inf_{\lambda > 0} \left( \lambda^{1-q} + t \|f\|_q^{-q/p'} \lambda^{1-q/p} \right) \\
&\simeq_{p,q} \|f\|_q \, t^{p'/q'};
\end{split}\]
in the last step, we used the fact that, by optimising in $\lambda>0$, one gets
\[
\lambda^{1-q} \simeq_{p,q} t \|f\|_q^{-q/p'} \lambda^{1-q/p},
\]
which readily reduces to $\lambda^{q/p'} \simeq_{p,q} \|f\|_q^{q/p'}/t$, i.e., $\lambda \simeq_{p,q} t^{-p'/q} \|f\|_q$.
Thus,
\[
\|f\|_{[H^{1,q}(X),L^p(X)]_{\theta,\infty}} = \sup_{t>0} t^{-p'/q'} K(t,f,H^{1,q}(X),L^p(X)) \lesssim_{p,q} \|f\|_q,
\]
as claimed.

The proof of part \ref{en:interprelim_fin} is entirely analogous: here we replace $H^{1,q}(X)$ by $\lH^{1,q}(X)$ and use the decomposition $f=g^\lambda+b^\lambda$ from Proposition \ref{prp:lH1}\ref{en:lH1_CZ}.
\end{proof}

\begin{prp}\label{prp:H1Lpinterp}
Let $(X,\cE,m,\Adm)$ be a CZ space. Assume that $1 < q < p \leq \infty$, and let $\theta = p'/q'$.
\begin{enumerate}[label=(\roman*)]
\item\label{en:H1Lpinterp_inf} If $m(X) = \infty$, then
\[
L^q(X) = [H^1(X),L^p(X)]_{\theta,q},
\]
with equivalence of norms (and universal implicit constants).
\item\label{en:H1Lpinterp_fin} If $m(X) < \infty$, then
\[
L^q(X) = [\lH^1(X),L^p(X)]_{\theta,q},
\]
with equivalence of norms (and universal implicit constants).
\end{enumerate}
\end{prp}
\begin{proof}
We only prove part \ref{en:H1Lpinterp_inf}, as the proof of part \ref{en:H1Lpinterp_fin} is entirely analogous.

We first consider the case $p<\infty$.
Take $r \in (1,q)$ and $s \in (q,p)$. Then, both $r,s \in (1,p)$, so, for appropriate $\theta_0,\theta_1 \in (0,1)$, we have
\[
\frac{1}{r} = \frac{1-\theta_0}{1} + \frac{\theta_0}{p}, \qquad \frac{1}{s} = \frac{1-\theta_1}{1} + \frac{\theta_1}{p},
\]
and therefore, by Lemma \ref{lem:interprelim},
\begin{align*}
L^r(X) &\subseteq [H^{1,r}(X),L^p(X)]_{\theta_0,\infty} = [H^1(X),L^p(X)]_{\theta_0,\infty}, \\
L^s(X) &\subseteq [H^{1,s}(X),L^p(X)]_{\theta_1,\infty} = [H^1(X),L^p(X)]_{\theta_1,\infty},
\end{align*}
where Proposition \ref{prp:H1equiv} was used in the last equalities.

On the other hand, as $q \in (r,s)$, for an appropriate $\eta \in (0,1)$ we have
\[
\frac{1}{q} = \frac{1-\eta}{r} + \frac{\eta}{s},
\]
and therefore
\[\begin{split}
L^q(X) 
&= [L^r(X),L^s(X)]_{\eta,q} \\
&\subseteq \left[[H^1(X),L^p(X)]_{\theta_0,\infty},[H^1(X),L^p(X)]_{\theta_1,\infty}\right]_{\eta,q} \\
&= [H^1(X),L^p(X)]_{\theta,q} \\
&\subseteq [L^1(X),L^p(X)]_{\theta,q} = L^q(X),
\end{split}\]
where we used the real interpolation properties of Lebesgue spaces, the reiteration theorem \cite[Theorem 3.5.3]{BL}, and the inclusion $H^1(X) \subseteq L^1(X)$.

We now consider the case $p=\infty$. Let us choose $s \in (q,\infty)$, and notice that $H^1(X) \cap L^\infty(X) \subseteq L^1(X) \cap L^\infty(X) \subseteq L^q(X) \cap L^s(X)$. Moreover,
\[
[H^1(X), L^s(X)]_{s'/q',q} = L^q(X), \qquad [L^q(X),L^\infty(X)]_{1-q/s,s} = L^s(X).
\]
So, by T.\ Wolff's reiteration theorem \cite[Theorem 1]{W},
\[
[H^1(X), L^\infty(X)]_{1/q',q} = L^q(X),
\]
as required.
\end{proof}

\begin{cor}\label{cor:H1Linfty}
Let $(X,\cE,m,\Adm)$ be a CZ space and $q \in (1,\infty)$. Let $f \in L^q(X)$.
\begin{enumerate}[label=(\roman*)]
\item\label{en:H1Linfty_inf} For all $\lambda > \|f\|_q/m(X)^{1/q}$, there exist $b \in H^1(X)$ and $g \in L^\infty(X)$ such that $f = b+g$ and
\[
\|b\|_{H^{1}} \lesssim_q \lambda^{1-q} \|f\|_q^q, \qquad \|g\|_\infty \lesssim_q \lambda,
\]
with universal implicit constants.
\item\label{en:H1Linfty_fin} If $m(X) < \infty$, then, for all $\lambda > 0$, there exist $b \in \lH^{1}(X)$ and $g \in L^\infty(X)$ such that $f = b+g$ and
\[
\|b\|_{\lH^{1}} \lesssim_q \lambda^{1-q} \|f\|_q^q, \qquad \|g\|_\infty \lesssim_q \lambda,
\]
with universal implicit constants.
\end{enumerate}
\end{cor}
\begin{proof}
If $m(X) = \infty$, from the interpolation result $L^q(X) = [H^1(X),L^\infty(X)]_{1/q',q}$ contained in Proposition \ref{prp:H1Lpinterp}, one deduces that, for all $f \in L^q(X)$ and $\lambda > 0$, there exists a decomposition $f = b+g$ as in part \ref{en:H1Linfty_inf}
(see \cite[Theorem 3.1.2, eq.\ (6)]{BL} applied with $s=\lambda^{-q} \|f\|_q^q$).

Similarly, in the case $m(X) < \infty$, by replacing $H^{1}(X)$ with $\lH^{1}(X)$ one proves the decomposition $f=b+g$ as in part \ref{en:H1Linfty_fin} for all $\lambda > 0$. In particular, if we write $b = c u + b_0$, where $b_0 \in H^{1}(X)$ and $u = \chr_X/m(X)$, then
\[
\|g\|_\infty \lesssim_q \lambda, \qquad |c| + \|b_0\|_{H^{1}} = \|b\|_{\lH^{1}} \lesssim_q \lambda^{1-q} \|f\|_q^q,
\]
and therefore, if $\lambda > \|f\|_q/m(X)^{1/q}$, then
\[
\|g+cu\|_\infty \lesssim_q \lambda (1+\lambda^{-q} \|f\|_q^q / m(X)) \lesssim \lambda;
\]
thus, the decomposition $f = b_0 + (g+cu)$ has the properties stated in part \ref{en:H1Linfty_inf}.
\end{proof}

We now establish some containment and density results for intersections.

\begin{lem}\label{lem:denseinters}
Let $(X,\cE,m,\Adm)$ be a CZ space and $q \in (1,\infty)$.
\begin{enumerate}[label=(\roman*)]
\item\label{en:denseinters_H1} If $p \in [q,\infty]$, then $H^{1}(X) \cap L^p(X)$ is dense in $H^{1}(X)$ and, if $m(X) = \infty$, it is also densely contained in $L^q(X)$.
\item\label{en:denseinters_lH1} If $m(X) < \infty$ and $p \in [q,\infty]$, then $\lH^{1}(X) \cap L^p(X)$ coincides with $L^p(X)$ and is densely contained both in $\lH^{1}(X)$ and in $L^q(X)$.
\item\label{en:denseinters_L1BMO} If $p \in [1,q)$, then $L^p(X) \cap \BMO(X) \subseteq L^q(X)$.
\end{enumerate}
\end{lem}
\begin{proof}
\ref{en:denseinters_H1}. From Proposition \ref{prp:Lq0H1}\ref{en:Lq0H1_Hfin}, we know that $H^{1,p}_\fin(X) \subseteq H^1(X) \cap L^p(X)$ and $H^{1,p}_\fin(X)$ is dense in $H^1(X)$, so a fortiori $H^{1}(X) \cap L^p(X)$ is dense in $H^{1}(X)$.

Assume that $m(X) = \infty$. Clearly $H^{1}(X) \cap L^p(X) \subseteq L^1(X) \cap L^p(X) \subseteq L^q(X)$. In order to prove the density of $H^{1}(X) \cap L^p(X)$ in $L^q(X)$, in light of Proposition \ref{prp:CZbasic}\ref{en:CZbasic_approx}, it is enough to show that the closure $\overline{H^{1}(X) \cap L^p(X)}^{L^q(X)}$ contains all well-supported $f \in L^q(X)$. On the other hand, for any such function $f$, we have that $f \chr_{\{|f| \leq n\}} \in L^\infty_\comp(X) \subseteq L^p(X)$ is also well-supported for any $n \in \NN$, and moreover $f \chr_{\{|f| \leq n\}} \to f$ in $L^q(X)$ as $n \to \infty$, by dominated convergence; so it is actually enough to show that $\overline{H^{1}(X) \cap L^p(X)}^{L^q(X)}$ contains all well-supported $f \in L^p(X)$.

Let $f \in L^p(X)$ be well-supported. Then, for any $n \in \Npos$ we can find $Q_n \in \Adm$ with $\{ f \neq 0 \} \subseteq Q_n$ and $m(Q_n) \geq n$, and then decompose
\[
f = f_n + c u_n,
\]
where
\[
c \defeq \int_X f \,dm, \qquad u_n \defeq \chr_{Q_n}/m(Q_n).
\]
Then, clearly $f_n \in L^p_0(Q_n) \subseteq L^p(X) \cap H^{1}(X)$ by Proposition \ref{prp:Lq0H1}\ref{en:Lq0H1_incl}. On the other hand, $\|u_n\|_q = m(Q_n)^{-1/q'} \to 0$ as $n \to \infty$, so indeed $f_n \to f$ in $L^q(X)$.

\ref{en:denseinters_lH1}. Recall from Proposition \ref{prp:lH1} that $\lH^{1}(X) = H^{1}(X) \oplus \CC u$, where $u \defeq \chr_X/m(X)$. 
Since $m(X) < \infty$, constants belong to $L^p(X)$, so $\lH^{1}(X) \cap L^p(X) = (H^{1}(X) \cap L^p(X)) \oplus \CC u$. As $H^{1}(X) \cap L^p(X)$ is dense in $H^{1}(X)$ by part \ref{en:denseinters_H1}, we conclude that $\lH^{1}(X) \cap L^p(X)$ is dense in $\lH^{1}(X)$.

On the other hand, from Proposition \ref{prp:lH1} we also know that $L^p(X) \subseteq \lH^{1}(X)$, so here $\lH^{1}(X) \cap L^p(X) = L^p(X)$, which is densely contained in $L^q(X)$.

\ref{en:denseinters_L1BMO}. 
In the case $m(X) < \infty$, we already know that $L^p(X) \cap \BMO(X) = \BMO(X) \subseteq L^q(X)$ by Proposition \ref{prp:BMOfin}\ref{en:BMOfin_Lq}.

Suppose instead that $m(X) = \infty$, and let $g \in L^p(X) \cap \BMO(X)$. Then, by Proposition \ref{prp:H1BMOprelim}, the duality pairing
\[
\phi_g(f) = \int_X f \, g \,dm
\]
makes sense for any $f \in H^{1,\infty}_\fin(X) + L^{p'}(X)$ and extends to a bounded linear functional on the Banach space sum $H^{1}(X) + L^{p'}(X)$. On the other hand, from Proposition \ref{prp:H1Lpinterp} we know that
\[
L^{q'}(X) = [H^{1}(X),L^{p'}(X)]_{p/q,q'} \subseteq H^{1}(X) + L^{p'}(X)
\]
with continuous embeddings. Therefore, the functional $\phi_g$ is also bounded on $L^{q'}(X)$, thus proving that $g \in L^q(X)$.
\end{proof}

By duality, we can deduce some interpolation results involving $\BMO$ spaces.

\begin{rem}
Let $(X,\cE,m,\Adm)$ be a CZ space. Recall from Proposition \ref{prp:BMOcompl} that the quotient $\BMO(X)/\CC$ is continuously embedded in $L^1_\loc(X)/\CC$. In addition, for any $q \in [1,\infty)$, if $m(X) = \infty$, then $L^q(X)$ does not contain constant functions other than zero, so the restriction to $L^q(X)$ of the canonical projection $L^1_\loc(X) \to L^1_\loc(X)/\CC$ gives a continuous embedding of $L^q(X)$ into $L^1_\loc(X)/\CC$. Thus, when $m(X) = \infty$, the spaces $L^q(X)$ and $\BMO(X)/\CC$ are compatible for any $q \in [1,\infty)$ and can be interpolated.
\end{rem}

\begin{prp}\label{prp:BMOLpinterp}
Let $(X,\cE,m,\Adm)$ be a CZ space. Assume that $1 \leq p < q < \infty$, and let $\theta = p/q$.
\begin{enumerate}[label=(\roman*)]
\item\label{en:BMOLpinterp_inf} If $m(X) = \infty$, then
\[
L^q(X) = [\BMO(X)/\CC,L^p(X)]_{\theta,q},
\]
with equivalence of norms (and universal implicit constants).
\item\label{en:BMOLpinterp_fin} If $m(X) < \infty$, then
\[
L^q(X) = [\BMO(X),L^p(X)]_{\theta,q},
\]
with equivalence of norms (and universal implicit constants).
\end{enumerate}
\end{prp}
\begin{proof}
In the case $p > 1$, we can identify $L^p(X)$ with the dual of $L^{p'}(X)$; so, due to the density results in Lemma \ref{lem:denseinters} and the duality results in Propositions \ref{prp:H1equiv}\ref{en:H1equiv_H1dual} and \ref{prp:fH1BMOduality}, the claimed interpolation relations follow from Proposition \ref{prp:H1Lpinterp} and the duality theorem for real interpolation \cite[Theorem 3.7.1]{BL}.

In the case $p=1$, from Lemma \ref{lem:denseinters}\ref{en:denseinters_L1BMO} we know that $L^1(X) \cap \BMO(X) \subseteq L^s(X)$ for all $s \in (1,\infty)$. So, the interpolation result in this case can be recovered, much as in the proof of Proposition \ref{prp:H1Lpinterp}, via an application of T. Wolff's reiteration theorem \cite[Theorem 1]{W}.
\end{proof}

We finally deduce an interpolation result between $H^1$ and $\BMO$.

\begin{cor}\label{cor:H1BMOinterp}
Let $(X,\cE,m,\Adm)$ be a CZ space and $q \in (1,\infty)$.
\begin{enumerate}[label=(\roman*)]
\item If $m(X) = \infty$, then
\[
L^q(X) = [H^1(X),\BMO(X)/\CC]_{1/q',q},
\]
with equivalence of norms (and universal implicit constants).
\item If $m(X) < \infty$, then
\[
L^q(X) = [\lH^1(X),\BMO(X)]_{1/q',q},
\]
with equivalence of norms (and universal implicit constants).
\end{enumerate}
\end{cor}
\begin{proof}
As $H^1(X) \cap \BMO(X) \subseteq L^1(X) \cap \BMO(X) \subseteq L^p(X)$ for all $p \in (1,\infty)$, and analogous inclusions hold for $m(X)<\infty$ with $\lH^1(X)$ in place of $H^1(X)$, this is again a consequence of Propositions \ref{prp:H1Lpinterp} and \ref{prp:BMOLpinterp}, and T. Wolff's reiteration theorem \cite[Theorem 1]{W}.
\end{proof}

Proposition \ref{prp:H1Lpinterp} and Corollary \ref{cor:H1BMOinterp} together complete the proof of Theorem \ref{thm:mainHardyinterp}. We end this section with our $H^1$--$L^1$ and $L^\infty$--$\BMO$ boundedness results.

\begin{proof}[Proof of Theorem \ref{thm:mainHardybdd}]
The bounds for $T^*$ follow by duality from those for $T$, so we shall only discuss the latter.

We shall prove that
\begin{equation}\label{eq:T_atom_bound}
C(T) \defeq \sup \{ \|Tb\|_1 \tc b \text{ $q$-atom} \} < \infty .
\end{equation}
This is enough to conclude that $T$ is bounded from $H^{1,q}(X)$ to $L^1(X)$ (hence also from $H^{1}(X)$ to $L^1(X)$ by Proposition \ref{prp:H1equiv}). Indeed, we already know that $T$ is of weak type $(1,1)$, so a fortiori $T$ is continuous from $H^{1,q}(X)$ to $L^{1,\infty}(X)$. Now, for any $f \in H^{1,q}(X)$, if we decompose $f = \sum_{n \in \NN} \lambda_n b_n$ for some $q$-atoms $b_n$ and $(\lambda_n)_n \in \ell^1$, then
\[
Tf = \sum_{n \in \NN} \lambda_n T b_n,
\]
with convergence in $L^{1,\infty}(X)$; on the other hand, under the assumption \eqref{eq:T_atom_bound}, the series on the right-hand side converges absolutely in $L^1(X)$, hence $Tf \in L^1(X)$ and $\|Tf\|_1 \leq C(T) \sum_{n \in \NN} |\lambda_n|$. Taking the infimum over all such decompositions of $f$, we finally obtain that $\|Tf\|_1 \leq C(T) \|f\|_{H^{1,q}}$, as required. In addition, if $m(X) < \infty$ and $u \defeq \chr_X / m(X)$, then
\[
\| T u \|_{1} \leq m(X)^{1/q'} \|Tu\|_q \leq m(X)^{1/q'} \|T\|_{q \to q} \|u\|_q  \leq \|T\|_{q \to q},
\]
thus showing that in this case $T$ is also bounded from $\lH^{1,q}(X)$ to $L^1(X)$.

To prove \eqref{eq:T_atom_bound}, let
\[
b = a - \left( \int a \,dm \right) \frac{1}{m(\Fth Q)} \chr_{\Fth Q}
\]
be a $q$-atom, where $Q,\Fth Q \in \Adm$ and $a \in L^q(X)$ satisfy $Q \in \Ch(\Fth Q)$, $\{ a \neq 0 \} \subseteq Q$, $\|a\|_q \leq m(Q)^{-1/q'}$. Notice that, in particular, $\|a\|_1 \leq 1$ and $\|b\|_1 \leq 2$. Then,
\[
\|T b\|_1 \leq \int_{X \setminus \Fth Q^*} |Tb| \,dm + \int_{\Fth Q^* \setminus Q^*} |Ta| \,dm + \int_{Q^*} |Ta| \,dm + \frac{1}{m(\Fth Q)} \int_{\Fth Q^*} |T\chr_{\Fth Q}| \,dm.
\]

Now, as $\{ b \neq 0\} \subseteq \Fth Q$ and $\int b \,dm = 0$, if $z$ is any point in $\Fth Q$,
\begin{multline*}
\int_{X \setminus \Fth Q^*} |Tb| \,dm = \int_{X \setminus \Fth Q^*} \left|\int_{\Fth Q} (K(x,y) - K(x,z)) b(y) \,dm(y) \right| \,dm(x) \\
\leq C_H(T) \|b\|_1 \leq 2 C_H(T).
\end{multline*}

Moreover, as $\{ a \neq 0 \} \subseteq Q$,
\[
\int_{\Fth Q^* \setminus Q^*} |Ta| \,dm = \int_{\Fth Q^* \setminus Q^*} \left|\int_Q K(x,y) a(y) \,dm(y) \right| \,dm(x) \leq C_S(T) \|a\|_1 \leq C_S(T).
\]

Further, by H\"older's inequality and \eqref{eq:extQ},
\[
\int_{Q^*} |Ta| \,dm \leq m(Q^*)^{1/q'} \|T a\|_{q} \leq m(Q^*)^{1/q'} \|T\|_{q \to q} m(Q)^{-1/q'} \leq \kappa^{1/q'} \|T\|_{q \to q}.
\]

Finally, again by H\"older's inequality and \eqref{eq:extQ},
\begin{multline*}
\frac{1}{m(\Fth Q)} \int_{\Fth Q^*} |T\chr_{\Fth Q}| \,dm \leq \frac{1}{m(\Fth Q)} m(\Fth Q^*)^{1/q'} \|T \chr_{\Fth Q}\|_{q} \\
\leq \frac{1}{m(\Fth Q)} m(\Fth Q^*)^{1/q'} \|T\|_{q \to q} m(\Fth Q)^{1/q} \leq \kappa^{1/q'} \|T\|_{q \to q}.
\end{multline*}

The claimed bound thus follows by combining the above estimates.
\end{proof}

\section{Scenarios of application. Flow trees}\label{sec:sec3}

\subsection{Metric CZ spaces}
We now introduce a variant of the definition of CZ spaces, which is applicable when a metric structure is also available.

\begin{dfn}\label{dfn:mCZdef}
A \emph{metric CZ space} is a metric-measure space $(X,d,\cE,m)$, where $\cE$ is the Borel $\sigma$-algebra on $(X,d)$, equipped with a collection $\Adm \subseteq \cE$ of \emph{admissible sets} satisfying conditions \ref{en:CZdef_meas}, \ref{en:CZdef_exh} and \ref{en:CZdef_chil} of Definition \ref{dfn:CZdef}, together with:
\begin{enumerate}[label=(\alph*'),start=4]
\item\label{en:mCZdef_diam} For all $Q \in \Adm$ and $x \in Q$, we have $\inf \{ \diam(R) \tc R \in \Ch^*(Q), \ x \in R \} = 0$.
\item For all $Q \in \Adm$, a doubling enlargement $Q^* \in \cE$ is assigned so that
\begin{equation}\label{eq:metricQstar}
m(Q^*) \leq \kappa m(Q), \qquad \{ x \in X \tc d(x,Q) \leq \delta \diam(Q) \} \subseteq Q^*
\end{equation}
for some constants $\delta \in (0,1]$ and $\kappa \in [1,\infty)$.
\end{enumerate}
\end{dfn}

\begin{prp}\label{prp:metricCZisCZ}
Let $(X,d,\cE,m,\Adm)$ be a metric CZ space. Then, $(X,\cE,m,\Adm)$ is a CZ space with enlargements. In addition, 
for the quantity $C_H(T)$ appearing in the smoothness condition \eqref{eq:CZop_horm} on the singular integral operator $T$ in Theorems \ref{thm:mainCZwt} and \ref{thm:mainHardybdd}, we have
\begin{equation}\label{eq:metric_CHT}
\begin{split}
C_H(T) &\leq \sup_{Q \in \Adm} \esssup_{y,z \in Q} \int_{\{x \tc d(x,Q) > \delta \diam(Q) \}} |K_T(x,y)-K_T(x,z)| \,dm(x) \\
&\leq \esssup_{y,z \in X} \int_{\{x \tc d(x,y),d(x,z) > \delta d(y,z) \}} |K_T(x,y)-K_T(x,z)| \,dm(x).
\end{split}
\end{equation}
\end{prp}
\begin{proof}
To show that $(X,\cE,m,\Adm)$ is a CZ space, we only need to check condition \ref{en:CZdef_sigma} of Definition \ref{dfn:CZdef}.

Let $Q \in \Adm$ and let $A \subseteq X$ be open; we claim that
\begin{equation}\label{eq:AcapQ_dyadic}
A \cap Q = \bigcup \{ R \in \Ch^*(Q) \tc R \subseteq A \cap Q \}.
\end{equation}
Indeed, for any $x \in A \cap Q$, we can find $r>0$ such that $B(x,r) \subseteq A$; thus, by property \ref{en:mCZdef_diam} of Definition \ref{dfn:mCZdef}, there exists $R \in \Ch^*(Q)$ with $x \in R$ and $\diam(R) < r$, and in particular $R \subseteq A \cap Q$.

As $\Ch^*(Q)$ is at most countable, \eqref{eq:AcapQ_dyadic} shows that $A \cap Q$ is in the $\sigma$-algebra $\Sigma_Q$ on $Q$ generated by $\Ch^*(Q)$. As $\cE$ is the Borel $\sigma$-algebra on $X$, $\cE|_Q$ is generated by sets of the form $A \cap Q$ with $A$ open in $X$, and therefore $\cE|_Q$ is contained in $\Sigma_Q$ too. Thus, $(X,\cE,m,\Adm)$ is a CZ space.

Finally, the inequalities in \eqref{eq:metric_CHT} are immediate consequences of the definition of $C_H(T)$ in \eqref{eq:CZop_horm} and the condition \eqref{eq:metricQstar}.
\end{proof}

From \eqref{eq:metric_CHT} we see that, in the case of metric CZ spaces, the smoothness condition in Theorems \ref{thm:mainCZwt} and \ref{thm:mainHardybdd} admits a stronger version that is purely expressed in terms of the metric and independent of the family $\Adm$ of admissible sets.
The following lemma gives a convenient way in applications to check the metric version of the smoothness assumption.

\begin{lem}\label{lem:dec_horm_metric}
Let $(X,d,\cE,m)$ be a metric-measure space. Let $K : X \times X \to \CC$ be measurable. Assume that we may decompose
\begin{equation}\label{eq:Kn_dec_metric}
K = \sum_{n \in \ZZ} K_n,
\end{equation}
where the kernels $K_n$ satisfy, for some constants $\varkappa,\varepsilon,\alpha > 0$ and $b > 1$, the estimates
\begin{align}
\label{eq:size_dec_metric}
\int_X |K_n(x,y)| \, (\dist(x,y)/b^n)^{\varepsilon} \,dm(x) &\leq \varkappa \\
\label{eq:grad_dec_metric}
\int_X |K_n(x,y)-K_n(x,z)| \,dm(x) &\leq \varkappa (d(y,z)/b^n)^{\alpha}
\end{align}
for all $n \in \ZZ$ and almost all $y,z \in X$. Then, for all $A \in \cE$ and $\gamma > 0$,
\[
\esssup_{y,z \in A} \int_{\{x \tc \dist(x,A) \geq \gamma\diam(A)\}} |K(x,y)-K(x,z)| \,dm(x) \lesssim_{\varkappa,b,\varepsilon,\alpha,\gamma} 1.
\]
\end{lem}
\begin{proof}
Let $A \in \cE$ and $y,z \in A$. From the decomposition \eqref{eq:Kn_dec_metric} we deduce that
\[\begin{split}
&\int_{\{x \tc \dist(x,A) \geq \gamma\diam(A)\}} |K(x,y)-K(x,z)| \,dm(x) \\
&\leq \sum_{n \in \ZZ} \int_{\{x \tc \dist(x,A) \geq \gamma\diam(A)\}} |K_n(x,y)-K_n(x,z)| \,dm(x) \\
&= \sum_{\substack{n \in \ZZ \\ b^n < \diam(A)}} + \sum_{\substack{n \in \ZZ \\ b^n \geq \diam(A)}} \eqdef \rmI + \rmII.
\end{split}\]

Now, by the triangle inequality,
\[\begin{split}
\rmI
&\leq \sum_{\substack{n \in \ZZ \\ b^n < \diam(A)}} \int_{\{x \tc \dist(x,A) \geq \gamma\diam(A)\}} |K_n(x,y)| \,dm(x) \\
&\qquad+ \sum_{\substack{n \in \ZZ \\ b^n < \diam(A)}} \int_{\{x \tc \dist(x,A) \geq \gamma\diam(A)\}} |K_n(x,z)| \,dm(x).
\end{split}\]
Since $y \in A$, in the first sum we have $\dist(x,y) \geq \dist(x,A) \geq \gamma\diam(A)$; thus, by \eqref{eq:size_dec_metric}, for a.a.\ $y\in A$,
\[\begin{split}
&\sum_{\substack{n \in \ZZ \\ b^n < \diam(A)}} \int_{\{x \tc \dist(x,A) \geq \gamma\diam(A)\}} |K_n(x,y)| \,dm(x) \\
&\lesssim_\gamma
\sum_{\substack{n \in \ZZ \\ b^n < \diam(A)}} (\diam(A)/b^n)^{-\varepsilon} \int_{\{x \tc \dist(x,A) \geq \gamma\diam(A)\}} |K_n(x,y)| \,(\dist(x,y)/b^n)^{\varepsilon} \,dm(x) \\
&\lesssim_{\varkappa} \sum_{\substack{n \in \ZZ \\ b^n < \diam(A)}} (b^n/\diam(A))^{\varepsilon} \lesssim_{b,\varepsilon} 1,
\end{split}\]
and an analogous estimate holds with $z$ in place of $y$.

On the other hand, by \eqref{eq:grad_dec_metric},
\[
\rmII
\lesssim_\varkappa \sum_{\substack{n \in \ZZ \\ b^n \geq \diam(A)}} (\dist(y,z)/b^n)^{\alpha}
\leq \sum_{\substack{n \in \ZZ \\ b^n \geq \diam(A)}} (\diam(A)/b^n)^{\alpha} \lesssim_{b,\alpha} 1,
\]
for a.a.\ $y,z \in A$, and we are done.
\end{proof}

With the help of Definition \ref{dfn:mCZdef} and Proposition \ref{prp:metricCZisCZ}, we can now list several examples of CZ spaces, showing that many singular integral theories developed in a variety of geometric settings can be subsumed under our paradigm.
\begin{itemize}
\item If $(X,\cE,d,m)$ is a doubling metric-measure space, we can make it into a regular metric CZ space by taking as admissible sets the dyadic cubes belonging to the adjacent dyadic systems constructed in \cite[Theorem 4.1]{HK12}. Up to equivalence of norms, our $H^1$--BMO and singular integral theory reduces in this context to the classical ones on spaces of homogeneous type \cite{CW71,CW}.
\item Beyond the doubling realm, the families of admissible sets constructed in \cite{HS,V} on $ax+b$ groups, in \cite{MOV,DMLMV25} on rank-one extensions of Carnot groups, in \cite{LSTV} on locally doubling flow trees, and in \cite{MPl} on certain hypergroups are easily seen to satisfy the conditions in Definition \ref{dfn:mCZdef}, as well as the regularity condition in Definition \ref{dfn:regularCZ}. The Hardy--BMO and singular integral theory from the previous sections thus recovers those established in the cited works.
\item The upcoming work \cite{BeLi} shows that Damek--Ricci spaces can also be equipped with a regular metric CZ space structure, thus providing a framework for Hardy--BMO theory and endpoint bounds for singular integral operators related to a distinguished Laplacian in that setting.
\item As we mentioned in the Introduction, any atomic filtration of $\sigma$-algebras gives rise to a (possibly non-regular) CZ space structure (with the $\sigma$-algebra atoms as admissible sets). This allows one to immediately recover the endpoint theory from \cite{LSMP} and the recent $H^1$-bounds from \cite{CAW25}.
\item Elaborating on the item above, \cite{CP19} builds a filtration on $\mathbb{R}^n$ --- equipped with a measure of polynomial growth --- that turns it into a CZ space. Moreover, Calderón--Zygmund kernels in that setting are shown to satisfy the size and Hörmander conditions from Theorem \ref{thm:mainCZwt}, so our theory implies endpoint bounds for them.
\end{itemize}
The above examples allow us to recover several known results from the literature. In the next sections we shall discuss in detail a previously unknown example. There, the regularity condition may fail, and therefore we shall need the full force of the theory developed in the previous sections to prove endpoint bounds for singular integrals. 

\subsection{Flow trees}

Here we show how the theory in the previous sections can be applied to the case of flow trees.
We recall briefly the main definitions and notation that we will use, and refer to \cite{LSTV,MSTV} for further details.

We shall work with a (locally finite) flow tree $(\Tree,\mu)$. Every vertex $x \in \Tree$ has a predecessor $\prd(x)$, which is a neighbour of $x$; the remaining neighbours are the successors of $x$, i.e., the elements of $\scc(x)$, and $\prd(y) = x$ for all $y \in \scc(x)$. Moreover, $\mu$ is a flow measure on $\Tree$, i.e., $\mu(x) = \sum_{y \in \scc(x)} \mu(y) \in (0,\infty)$ for all $x \in \Tree$; the measure of a subset $A \subseteq \Tree$ is given by $\mu(A) = \sum_{x \in A} \mu(x)$. A level function $\lev : \Tree \to \ZZ$ is assigned, so that $\lev(\prd(x)) = \lev(x)+1$ for all $x \in \Tree$; the level sets of $\lev$ are called horocycles. The geodesic distance function on $\Tree$ is denoted by $\dist$.

We write $\CC^\Tree$ for the space of all complex-valued functions on $\Tree$ and $\czz(\Tree)$ for the subspace of finitely supported functions. Lebesgue spaces $L^p(\Tree)$ for $p \in [1,\infty]$ are defined with respect to the flow measure $\mu$.

The shift operator $\opSh$ and its formal adjoint $\opSh^*$ on $(\Tree,\mu)$ are given by
\[
\opSh f(x) = f(\prd(x)), \qquad \opSh^* f(x) = \frac{1}{\mu(x)} \sum_{y \in \scc(x)} f(y)\mu(y),
\]
for all $f \in \CC^\Tree$ and $x \in \Tree$. The flow gradient $\nabla$ and the flow Laplacian $\opL$ on $(\Tree,\mu)$ are given by
\[
\nabla  = I - \opSh, \qquad \opL = \frac{1}{2} \nabla^* \nabla = I - \frac{1}{2} (\opSh + \opSh^*).
\]
The flow Laplacian is a positive self-adjoint operator on $L^2(\Tree)$ with spectrum $[0,2]$.
The Riesz transform $\Rz$ on $(\Tree,\mu)$ is given by
\[
\Rz = \nabla \opL^{-1/2}.
\]

For any linear operator $T : \czz(\Tree) \to \CC^\Tree$, its integral kernel is denoted by $K_T$; thus
\[
T f(x) = \sum_{y \in \Tree} K_T(x,y) \, f(y) \,\mu(y)
\]
for all $f \in \czz(\Tree)$ and $x \in \Tree$.

\subsection{CZ theory on flow trees}

We shall now equip an arbitrary flow tree with a metric CZ space structure.
The definition of ``admissible set'' on a flow tree presented below is a minor modification of that of ``admissible trapezoid'' in \cite[Section 3]{LSTV} for locally doubling flow trees.

\begin{dfn}\label{dfn:CZtreedef}
Let $(\Tree,\mu)$ be a flow tree.
\begin{enumerate}[label=(\alph*)]
\item For any $x_0 \in \Tree$ and $h',h'' \in \NN$ with $h' < h''$, we define
\[
R_{h',h''}(x_0) \defeq \{ x \in \Tree \tc x \leq x_0, \ h' \leq \dist(x,x_0) < h'' \}.
\]
\item We say that a subset $Q$ of $\Tree$ is \emph{admissible} if
\begin{itemize}
\item either $Q = \{x_0\}$ for some $x_0 \in \Tree$,
\item or $Q = R_{h',h''}(x_0)$ for some $x_0 \in \Tree$ and $h',h'' \in \Npos$ and $2 \leq h''/h' < 8$.
\end{itemize}
We denote by $\Adm$ the collection of admissible subsets of $\Tree$.
\item\label{en:CZtreedef_ch} If $Q \in \Adm$, we declare the following sets to be the \emph{children} of $Q$:
\begin{itemize}
\item if $Q = \{x_0\}$, then the only child of $Q$ is $Q$ itself;
\item if $Q = R_{h',h''}(x_0)$ with $4 \leq h''/h' < 8$, then $R_{h',2h'}(x_0)$ and $R_{2h',h''}(x_0)$ are the children of $Q$;
\item if $Q = R_{h',h''}(x_0)$ with $2 \leq h''/h' < 4$ and $h' \geq 2$, then the children of $Q$ are the sets $R_{h'-1,h''-1}(y)$ with $y \in \scc(x_0)$;
\item if $Q = R_{1,2}(x_0)$, then the children of $Q$ are the sets $\{y\}$ with $y \in \scc(x_0)$;
\item if $Q = R_{1,3}(x_0)$, then the children of $Q$ are the sets $\{y\}$ and $R_{1,2}(y)$ with $y \in \scc(x_0)$.
\end{itemize}
We denote by $\Ch(Q)$ the collection of the children of $Q$.
\item For any $Q \in \Adm$, we define its enlargement $Q^* \subseteq \Tree$ as follows:
\begin{itemize}
\item if $Q = \{x_0\}$, then $Q^* = \{x_0\}$;
\item if $Q = R_{h',h''}(x_0)$ with $2 \leq h''/h' < 8$, then $Q^* = R_{0,h''+h'}(x_0)$.
\end{itemize}
\end{enumerate}
\end{dfn}

\begin{prp}\label{prp:treeCZ}
Let $(\Tree,\mu)$ be a flow tree with level function $\lev$ and geodesic distance $\dist$.
\begin{enumerate}[label=(\roman*)]
\item\label{en:treeCZ_ch} For any $Q \in \Adm$, the children of $Q$ are admissible and form a partition of $Q$.
\item\label{en:treeCZ_dil} For any $Q \in \Adm$, we have $Q \subseteq Q^*$ and $\mu(Q^*) \leq 3 \mu(Q)$.
\item\label{en:treeCZ_dildist} For any $Q \in \Adm$, we have 
\[
\{ x \in \Tree \tc \dist(x,Q) \leq \diam(Q)/16 \} \subseteq Q^*
\]
and
\[
\inf_{y \in Q} \dist(y,\Tree \setminus Q^*) \geq \# \lev(Q) / 7;
\]
moreover, for all $R \in \Ch(Q)$,
\[
\# \lev(Q^*) \leq 21 \# \lev(R).
\]
\item\label{en:treeCZ_CZ} The metric-measure space $(\Tree,\dist,\mu)$ with the family $\Adm$ of admissible sets and the associated enlargements is a metric CZ space.
\end{enumerate}
\end{prp}
\begin{proof}
\ref{en:treeCZ_ch}. It is easily checked that the children of $Q$ form a partition of $Q$, in each of the cases of Definition \ref{dfn:CZtreedef}\ref{en:CZtreedef_ch} above.

As for admissibility, we simply observe that, if $h',h'' \in \Npos$ and $4 \leq h''/h' < 8$, then clearly $2h'/h' = 2$, while $2 \leq h''/2h' < 4$, so both $2h'/h', h''/2h' \in [2,8)$; if instead $2 \leq h''/h' < 4$ and $h' \geq 2$, then $h'-1 \geq 1$ and $h''-1 \leq h'' < 4h' = 4(h'-1+1) \leq 8(h'-1)$, while $2(h'-1) \leq h''-2 \leq h''-1$, so again $(h''-1)/(h'-1) \in [2,8)$.

\ref{en:treeCZ_dil}. The inequality is trivially true when $Q = \{x_0\} = Q^*$. Moreover, since $\mu$ is a flow measure, $\mu(R_{h',h''}(x_0)) = (h''-h') \mu(x_0)$. So, if $Q = R_{h',h''}(x_0)$ with $2 \leq h''/h' < 8$, then $h''-h' \geq h'$ and $h''+h' = h''-h'+2h' \leq 3(h''-h')$, whence $\mu(Q^*) = (h'+h'')\mu(x_0) \leq 3\mu(Q)$.

\ref{en:treeCZ_dildist}. If $Q = \{x_0\}$, then $\diam(Q) = 0$, $\# \lev(Q) = 1$ and $Q^*= Q=R$ for any $R \in \Ch(Q)$, so the estimates are trivially true.

Assume instead that $Q = R_{h',h''}(x_0)$ with $h',h'' \in \Npos$ and $2 \leq h''/h' < 8$. Then $\diam(Q) \leq 2h''< 16 h'$, so 
\[
Q^* = R_{0,h'+h''}(x_0) = \{ x \in \Tree \tc \dist(x,Q) \leq h'\} \supseteq \{ x \in \Tree \tc \dist(x,Q) \leq \diam(Q)/16 \}.
\]
Further, $\# \lev(Q) = h''-h' < 7h'$; so, for all $y \in Q$, 
\[
\dist(y,\Tree \setminus Q^*) \geq h' \geq \# \lev(Q) /7,
\]
and moreover
\[
\# \lev(Q^*) = h'+h'' \leq 3(h''-h') \leq 3 \# \lev(Q).
\]
Finally, if $R \in \Ch(Q)$, then, by inspection of the cases in \ref{dfn:CZtreedef}\ref{en:CZtreedef_ch}, one sees that either $\#\lev(R)=\#\lev(Q)$, or $h'' \geq 3h'$ and $\# \lev(R) \geq h''-2h' \geq h' \geq \# \lev(Q)/7$. So in any case $\# \lev(Q^*) \leq 3 \# \lev(Q) \leq 21 \# \lev(R)$, as required.

\ref{en:treeCZ_CZ}. 
Let us check the conditions in Definition \ref{dfn:mCZdef}. Notice that, as $\dist$ induces the discrete topology on $\Tree$ and the latter is countable, any subset of $\Tree$ is Borel.

By definition, each $Q \in \Adm$ is finite and nonempty, so $0 < \mu(Q) < \infty$.

Take any $o \in \Tree$ with $\lev(o)=0$, and consider, for any $n \in \Npos$, the set $Q_n \defeq R_{n,3n}(\prd^{2n}(o))$. Then clearly $Q_n \in \Adm$ and $Q_n \subseteq Q_{n+1}$ for all $n \in \Npos$, and moreover $\bigcup_{n \in \Npos} Q_n = \Tree$. Indeed, for any $x \in \Tree$, there exists $k \in \NN$ such that $x \leq \prd^k(o)$, and then $x \in Q_n$ for any $n > \max\{ k, -\lev(x) \}$.

The fact that the children of any $Q \in \Adm$ are a (finite) partition of $Q$ made of admissible sets is proved in part \ref{en:treeCZ_ch}.

Furthermore, by inspection of the various cases in Definition \ref{dfn:CZtreedef}\ref{en:CZtreedef_ch} above, we see that the children of the children of an admissible set of the form $R_{h',h''}(x_0)$ are either singletons or sets $R_{k',k''}(y)$ with $k'' < h''$, and therefore $\Ch^{2h''}(R_{h',h''}(x_0))$ is made of singletons, which have diameter zero.

Finally, the required properties for the enlargement $Q^*$ of any $Q \in \Adm$ are proved in parts \ref{en:treeCZ_dil} and \ref{en:treeCZ_dildist}.
\end{proof}

In light of Proposition \ref{prp:treeCZ}, the theory developed in the previous section is applicable to flow trees. 
In particular, we can use Theorems \ref{thm:mainCZwt} and \ref{thm:mainHardybdd} to recover boundedness properties of singular integral operators in this setting. We state here a consequence of those theorems in the setting of flow trees,
where the assumptions on the operator do not directly refer to admissible sets. 

\begin{prp}\label{prp:treeCZop}
Let $(\Tree,\mu)$ be a flow tree, with level function $\lev$ and geodesic distance function $\dist$. Let
$q \in (1,\infty)$ and $T : L^q(\Tree) \to L^q(\Tree)$ be a bounded operator. Assume that, for some $\varkappa \in (0,\infty)$, the integral kernel $K_T : \Tree \times \Tree \to \CC$ of the operator $T$ satisfies the following conditions:
\begin{enumerate}[label=(\alph*)]
\item\label{en:treeCZop_size} for all $A \subseteq \Tree$ and $y \in \Tree$,
\[
\sum_{x \in A} |K_T(x,y)| \,\mu(x) \leq \varkappa \left(1 + \frac{\# \lev(A)}{1+\dist(y,A)}\right);
\]
\item\label{en:treeCZop_horm} for all $A \subseteq \Tree$,
\[
\sup_{y,z \in A} \sum_{\substack{x \in \Tree \\ \dist(x,A) \geq \diam(A)/16}} |K_T(x,y)-K_T(x,z)| \,\mu(x) \leq \varkappa.
\]
\end{enumerate}
Then, the operator $T$ is of weak type $(1,1)$, bounded on $L^p(\Tree)$ for any $p \in (1,q]$, and bounded from $H^1(\Tree)$ to $L^1(\Tree)$, with bounds only depending on $\varkappa$, $\|T\|_{q \to q}$, $q$, and $p$ (but not on the flow tree $\Tree$).
\end{prp}
\begin{proof}
In light of Proposition \ref{prp:treeCZ}, it is sufficient to check that the integral kernel $K_T$ satisfies the conditions in Theorem \ref{thm:mainCZwt}.

Now, if $Q \in \Adm$, then, by Proposition \ref{prp:treeCZ}\ref{en:treeCZ_dildist}, for all $y,z \in Q$,
\begin{multline*}
\sum_{x \in \Tree \setminus Q^*} |K_T(x,y) - K_T(x,z)| \,\mu(x) \\
\leq \sum_{\substack{ x \in \Tree \\ \dist(x,Q) \geq \diam(Q)/16}} |K_T(x,y) - K_T(x,z)| \,\mu(x) \leq \varkappa.
\end{multline*}
Moreover, if $\Fth Q \in \Adm$ and $Q \in \Ch(\Fth Q)$, then, again by Proposition \ref{prp:treeCZ}\ref{en:treeCZ_dildist}, for all $y \in Q$,
\begin{multline*}
\sum_{x \in \Fth Q^* \setminus Q^*} |K_T(x,y)| \,\mu(x) 
\leq \varkappa\left(1+\frac{\# \lev(\Fth Q^* \setminus Q^*)}{1+\dist(y,\Fth Q^* \setminus Q^*)} \right)\\
\leq \varkappa\left(1+\frac{\# \lev(\Fth Q^*)}{\dist(y,\Tree \setminus Q^*)} \right) \leq 148 \varkappa,
\end{multline*}
as required.
\end{proof}

\begin{rem}\label{rem:dec_horm_tree}
In order to verify assumption \ref{en:treeCZop_horm} in Proposition \ref{prp:treeCZop}, we may apply Lemma \ref{lem:dec_horm_metric} with $K = K_T$. In the case of a flow tree $\Tree$, the conditions \eqref{eq:size_dec_metric}-\eqref{eq:grad_dec_metric} with $\eta=1$ in Lemma \ref{lem:dec_horm_metric} can be equivalently rewritten as
\begin{align}
\label{eq:size_dec}
\sup_{y \in \Tree} \sum_{x \in \Tree} |K_n(x,y)| \, (\dist(x,y)/b^n)^\varepsilon \,\mu(x) &\leq \varkappa \\
\label{eq:grad_dec}
\sup_{y \in \Tree} \sum_{x \in \Tree} |K_n(x,y)-K_n(x,\prd(y))| \,\mu(x) &\leq \varkappa b^{-n}
\end{align}
for all $n \in \ZZ$.
\end{rem}

\subsection{Riesz transforms}
We now prove the boundedness results for Riesz transforms on flow trees stated in Theorem \ref{thm:maintrees}.
To that end, we shall need the following heat kernel estimates.

\begin{lem}
For all $t \geq 1$ and $\varepsilon > 0$,
\begin{equation}\label{prp:level_heatkernelest}
\sup_{y \in \Tree} \sup_{n \in \ZZ} \sum_{\substack{x \in \Tree \\ \lev(x) = n}} |K_{\nabla \exp(-t\opL)}(x,y)| \,\exp(\varepsilon \dist(x,y)/\sqrt{t}) \,\mu(x) \lesssim_\varepsilon t^{-1}.    
\end{equation}
\end{lem}
\begin{proof}
Estimate \eqref{prp:level_heatkernelest} is contained in \cite[Proposition 2.8]{MSV} in the case where $(\Tree,\mu)$ is a homogeneous tree with canonical flow, with an implicit constant that is independent of the degree of the tree. By transference \cite[Theorem 4.8]{MSTV}, the same estimate holds true on any flow tree.
\end{proof}

\begin{proof}[Proof of Theorem \ref{thm:maintrees} for Riesz transforms] The result is already known when $\Tree$ is locally doubling, as shown in \cite[Theorem 1.2]{MSTV}; indeed, in \cite[Theorem 1.2]{MSTV} the $L^p$-bounds for $p \in (1,\infty)$ were proved for arbitrary flow trees $\Tree$. What remains to prove are the endpoint bounds (weak-type $(1,1)$ and $H^1(\Tree) \to L^1(\Tree)$ boundedness) for an arbitrary flow tree $\Tree$. 

We shall thus prove size and H\"ormander estimates for the Riesz transform kernel, showing that $\Rz$ satisfies the assumptions of Proposition \ref{prp:treeCZop}. We first show that for any $y \in \Tree$ and $A \subseteq \Tree$,
\[
\sum_{x \in A} |K_{\Rz}(x,y)| \,\mu(x) \lesssim 1+\frac{\# \lev(A)}{1+\dist(y,A)}.
\]
From \eqref{prp:level_heatkernelest} it readily follows that, for all $t \geq 1$,
\begin{equation}\label{eq:level_heatkernelest_A}
\sum_{x \in A} |K_{\nabla\exp(-t\opL)}(x,y)| \,\mu(x) \lesssim \exp(-\dist(y,A)/\sqrt{t}) \frac{\# \lev(A)}{t}.
\end{equation}

From the subordination formula
\[
\Rz = \nabla \opL^{-1/2} = \frac{1}{\sqrt{\pi}} \int_0^\infty \nabla \exp(-t\opL) \,\frac{dt}{\sqrt{t}}
\]
we deduce that
\[
\sum_{x \in A} |K_{\Rz}(x,y)| \,\mu(x)
\lesssim \int_0^\infty \sum_{x \in A} |K_{\nabla\exp(-t\opL)}(x,y)| \,\mu(x) \,\frac{dt}{\sqrt{t}} = \int_0^1 + \int_1^\infty \eqdef \mathrm{I} + \mathrm{II}.
\]
Now,
\[
\mathrm{I} \leq \int_0^1 \| \nabla \exp(-t\opL) \|_{1 \to 1} \,\frac{dt}{\sqrt{t}} \lesssim 1,
\]
as $\sup_{t > 0} \|\nabla \exp(-t\opL)\|_{1 \to 1} \lesssim 1$. On the other hand, by \eqref{eq:level_heatkernelest_A} we obtain that
\[
\mathrm{II} \leq \int_1^\infty \exp(-\dist(y,A)/\sqrt{t}) \frac{\# \lev(A)}{t} \, \frac{dt}{\sqrt{t}} \simeq \frac{\# \lev(A)}{1+\dist(y,A)},
\]
and we are done. 

Finally, we show that for all $\gamma>0$,
\begin{equation}\label{eq:Rz_horm}
\sup_{y,z \in A} \sum_{\substack{x \in \Tree \\ \dist(x,A) \geq \gamma \diam(A)}} |K_{\Rz}(x,y)-K_{\Rz}(x,z)| \,\mu(x) \lesssim_\gamma 1.
\end{equation}
Much as in \cite[proof of Theorem 1.2, case $p \leq 2$]{MSTV}, we may decompose
\[
\Rz = \Rz^{(0)} + \Rz^{(\infty)} = \Rz^{(0)} + \sum_{n \in \NN} \Rz^{(\infty)}_n.
\]
As $\Rz^{(0)}$ is $L^1$-bounded,
\[
\sup_{y,z \in \Tree} \sum_{x \in \Tree} |K_{\Rz^{(0)}}(x,y)-K_{\Rz^{(0)}}(x,z)| \,\mu(x) \leq 2 \|\Rz^{(0)}\|_{1 \to 1} < \infty,
\]
so the estimate \eqref{eq:Rz_horm} with $K_{\Rz^{(0)}}$ in place of $K_{\Rz}$ holds trivially. The analogous estimate for $K_{\Rz^{(\infty)}}$ follows from Remark \ref{rem:dec_horm_tree}, thanks to the bounds in \cite[eq.\ (6.5)]{MSTV} for the kernels $K_{\Rz^{(\infty)}_n}$. 
\end{proof}

\subsection{Spectral multipliers} We now prove Theorem \ref{thm:maintrees} for spectral multipliers. Again, the locally doubling case was previously known in full \cite{MSTV}, while the $L^p$-bounds for $p\in (1,\infty)$ were also known for arbitrary flow trees. What remains to be proved are the endpoint bounds for arbitrary flow trees.

We shall need the following estimates, relating the functional calculus of the flow Laplacian $\opL$ on an arbitrary flow tree $\Tree$ and that of the discrete Laplacian $\Delta_\ZZ$ on $\ZZ$. As in \cite{MSTV}, we define
\[
\Delta_\ZZ f(n) = f(n) - \frac{1}{2}(f(n+1)+f(n-1)), \qquad \tilde\nabla_\ZZ f(n) = f(n-1)-f(n+1)
\]
for all $f \in \CC^\ZZ$ and $n \in \ZZ$; we also write $k_T$ for the convolution kernel of a translation-invariant operator $T : \czz(\ZZ) \to \CC^\ZZ$.

\begin{lem}
Let $s>3/2$. For all $F \in \sobolev{2}{s}(\RR)$ and all increasing weights $w : \NN \to [0,\infty)$,
\begin{align}
\sup_{y \in \Tree} \sum_{x \in \Tree} |K_{F(\opL)}(x,y)| \,w(\dist(x,y)) \,\mu(x) &\lesssim \sum_{n \in \NN} n w(n) \, |\tilde\nabla_\ZZ k_{F(\Delta_\ZZ)}(n)|, \\
\label{eq:redZ_w_lev}
\sup_{l \in \ZZ} \sup_{y \in \Tree} \sum_{\substack{x \in \Tree \\ \lev(x) = l}} |K_{F(\opL)}(x,y)| \,w(\dist(x,y)) \,\mu(x) &\lesssim \sum_{n \in \NN} w(n) \, |\tilde\nabla_\ZZ k_{F(\Delta_\ZZ)}(n)|.
\end{align}
\end{lem}
\begin{proof}
By transference \cite[Theorem 5.11]{MSTV}, it is enough to prove the estimates in the case where $(\Tree,\mu)$ is a homogeneous tree with canonical flow, with implicit constants independent of the degree of the tree.

In the homogeneous tree case, the first estimate is proved in \cite[Proposition 5.6]{MSTV} for all $F \in C[0,2]$; a straightforward modification of that proof, using \cite[eq.\ (2.17)]{MSV} in place of \cite[eq.\ (2.15)]{MSV}, gives the second estimate too.
\end{proof}

\begin{prp}\label{prp:L1est}
The following estimates hold.
\begin{enumerate}[label=(\roman*)]
\item\label{en:L1est_nosupp} For all $s > 3/2$ and $F \in \sobolev{2}{s}(\RR)$,
\[
\sup_{y \in \Tree} \sum_{x \in \Tree} |K_{F(\opL)}(x,y)| \,\mu(x) \lesssim_s \|F\|_{\sobolev{2}{s}}.
\]
\item\label{en:L1est_supp} For all $\alpha \in [0,\infty)$ and $\beta > \alpha+3/2$, for all $F \in \sobolev{2}{\beta}(\RR)$ supported in $[1/4,7/4]$, and for all $t \geq 1$,
\begin{align}
\sup_{y \in \Tree} \sum_{x \in \Tree} |K_{F(t\opL)}(x,y)| \left(1+\frac{\dist(x,y)}{\sqrt{t}}\right)^\alpha \,\mu(x) &\lesssim_{\alpha,\beta} \|F\|_{\sobolev{2}{\beta}},\\
\sup_{l \in \ZZ} \sup_{y \in \Tree} \sum_{\substack{x \in \Tree \\ \lev(x)=l}} |K_{F(t\opL)}(x,y)| \left(1+\frac{\dist(x,y)}{\sqrt{t}}\right)^{\alpha+1} \,\mu(x) &\lesssim_{\alpha,\beta} t^{-1/2} \|F\|_{\sobolev{2}{\beta}}.
\end{align}
\end{enumerate}
\end{prp}
\begin{proof}
\ref{en:L1est_nosupp}. This is proved in \cite[Proposition 5.8]{MSTV}.

\ref{en:L1est_supp}. The first estimate is proved in \cite[Proposition 6.9]{MSTV}. The proof of the second estimate is similar. Indeed, let $\varepsilon > 0$ be such that $\beta = \alpha+3/2+\varepsilon$. By applying \eqref{eq:redZ_w_lev} to the function $F(t\cdot) \in \sobolev{2}{\beta}$ and the weight $w(n) = (1+n/\sqrt{t})^{\alpha+1}$, we obtain that
\begin{multline*}
\sup_{l \in \ZZ} \sup_{y \in \Tree} \sum_{\substack{x \in \Tree \\ \lev(x) = l}} |K_{F(t\opL)}(x,y)| \left(1+\frac{\dist(x,y)}{\sqrt{t}}\right)^{\alpha+1} \,\mu(x) \\
\lesssim \sum_{n \in \NN} \left(1+\frac{n}{\sqrt{t}}\right)^{\alpha+1} |\tilde\nabla_\ZZ k_{F(t\Delta_\ZZ)}(n)|,
\end{multline*}
and the latter quantity is bounded by $t^{-1/2} \|F\|_{\sobolev{2}{\beta}}$ exactly as in the proof of \cite[Proposition 6.9]{MSTV}.
\end{proof}

We can now prove the desired kernel estimates for spectral multipliers of $\opL$, showing that they satisfy the assumptions of Proposition \ref{prp:treeCZop} under suitable smoothness assumptions.

\begin{prp}\label{prp:MHmult_kernelests}
Let $\chi \in C^\infty_c(\RR)$ be supported in $(3/4,5/4)$ and such that $0 \leq \chi \leq 1$ and $\chi(1) = 1$.
Let $s>3/2$. Let $F : \RR \to \CC$ be a Borel measurable function satisfying
\[
\|F\|_{L^2_{s,\sloc}} \defeq \sup_{0 < t \leq 2} \|F(t\cdot) \chi\|_{\sobolev{2}{s}} < \infty.
\]
Let $A \subseteq \Tree$.
\begin{enumerate}[label=(\roman*)]
\item\label{en:mult_size} For any $y \in \Tree$,
\begin{equation}\label{eq:mult_size}
\sum_{x \in A} |K_{F(\opL)}(x,y)| \,\mu(x) \lesssim_s \left( 1+\frac{\# \lev(A)}{1+\dist(y,A)}\right) \|F\|_{L^2_{s,\sloc}}.
\end{equation}
\item\label{en:mult_horm} For all $\gamma>0$,
\begin{equation}\label{eq:mult_horm}
\sup_{y,z \in A} \sum_{\substack{x \in \Tree \\ \dist(x,A) \geq \gamma \diam(A)}} |K_{F(\opL)}(x,y)-K_{F(\opL)}(x,z)| \,\mu(x) \lesssim_{s,\gamma} \|F\|_{L^2_{s,\sloc}}.
\end{equation}
\end{enumerate}
\end{prp}
\begin{proof}
As in the proof of \cite[Theorem 1.4]{MSTV}, we may decompose
\[
F(\opL) = F^{(0)}(\opL) + F^{(1)}(\opL) = \sum_{k \in \NN} F^{(0)}_k(2^k \opL) + F^{(1)}(\opL),
\]
where $\supp F^{(0)}_k \subseteq (1/4,1)$ for all $k \in \NN$, $\supp F^{(1)} \subseteq (1/4,\infty)$, and
\[
\|F^{(1)}\|_{\sobolev{2}{s}} + \sup_{k \in \NN} \|F^{(0)}_k\|_{\sobolev{2}{s}} \lesssim_s \|F\|_{L^2_{s,\sloc}}.
\]

In particular, Proposition \ref{prp:L1est}\ref{en:L1est_nosupp} implies that
\[
\sup_{y \in \Tree} \sum_{x \in \Tree} |K_{F^{(1)}(\opL)}(x,y)| \,\mu(x) \lesssim_s \|F\|_{L^2_{s,\sloc}},
\]
whence we deduce the analogues of both \eqref{eq:mult_size} and \eqref{eq:mult_horm} where $K_{F(\opL)}$ is replaced with $K_{F^{(1)}(\opL)}$. It only remains to prove the corresponding estimates for $K_{F^{(0)}(\opL)}$.

Of these, the analogue of \eqref{eq:mult_horm} for $K_{F^{(0)}(\opL)}$ follows from Lemma \ref{lem:dec_horm_metric} and Remark \ref{rem:dec_horm_tree}, together with the bounds for the kernels $K_{F^{(0)}_k(2^k\opL)}$ in \cite[eq. (6.17)]{MSTV}.

As for the remaining estimate, if $\varepsilon>0$ is such that $s > 3/2+\varepsilon$, then
by applying Proposition \ref{prp:L1est}\ref{en:L1est_supp} to the function $F^{(0)}_k$ and $t=2^k$, we deduce that
\[
\sup_{l \in \ZZ} \sup_{y \in \Tree} \sum_{\substack{x \in \Tree \\ \lev(x)=l}} |K_{F^{(0)}_k(2^k\opL)}(x,y)| \left(1+\frac{\dist(x,y)}{2^{k/2}}\right)^{1+\varepsilon} \,\mu(x) \lesssim_{s} 2^{-k/2} \|F\|_{L^2_{s,\sloc}},
\]
and therefore, if $A \subseteq \Tree$ and $y \in \Tree$, then
\[
\sum_{x \in A} |K_{F^{(0)}_k(2^k\opL)}(x,y)| \,\mu(x) \lesssim_s \frac{\# \lev(A)}{2^{k/2}} \left(1+\frac{\dist(y,A)}{2^{k/2}}\right)^{-1-\varepsilon} \|F\|_{L^2_{s,\sloc}}.
\]
Summing in $k \in \NN$ thus yields
\[
\sum_{x \in A} |K_{F^{(0)}(\opL)}(x,y)| \,\mu(x) \lesssim_s \frac{\# \lev(A)}{1+\dist(y,A)} \|F\|_{L^2_{s,\sloc}},
\]
where we used the fact that
\[
\sum_{k \in \NN} 2^{-k/2} ( 1+ 2^{-k/2} D)^{-1-\varepsilon} \simeq_\varepsilon \frac{1}{1+D}
\]
for any $D \in [0,\infty)$.
\end{proof}

\begin{proof}[Proof of Theorem \ref{thm:maintrees} for spectral multipliers]
In the case $F$ satisfies the more restrictive assumption 
\begin{equation*}
\sup_{0 < t \leq 2}\|F(t\cdot)\chi\|_{L^2_s}<\infty,
\end{equation*}
Proposition \ref{prp:MHmult_kernelests} shows that the integral kernel of $F(\opL)$ satisfies the assumptions of Proposition \ref{prp:treeCZop}, whence we deduce the weak-type $(1,1)$ and $H^1(\Tree) \to L^1(\Tree)$ bounds; moreover, as $F(\opL)^* = \overline{F}(\opL)$ and $\overline{F}$ satisfies the same smoothness assumption as $F$, by duality we also obtain the $L^\infty(\Tree) \to \BMO(\Tree)/\CC$ bound.

If instead $F$ satisfies the assumption 
\begin{equation*}
\sup_{0<t\leq 1}\|F(t\cdot)\chi\|_{L^2_s}+\sup_{0<t\leq 1}\|F(2-t\cdot)\chi\|_{L^2_s}<\infty,
\end{equation*}
then, much as in the proof of \cite[Theorem 1.5]{MSTV},  by using cutoff functions, we can decompose $F = F_0 + F_2$, where $\supp F_0 \subseteq (-\infty,5/4)$ and $\supp F_2 \subseteq (3/4,\infty)$, so that
\[
\sup_{t>0} \|F_0(t\cdot) \chi\|_{L^2_s} < \infty, \qquad \sup_{t>0} \|F_2(2-t\cdot) \chi\|_{L^2_s} < \infty,
\]
and
\[
F(\opL) = F_0(\opL) + F_2(\opL) = F_0(\opL) + \modE F_2(2-\opL) \modE,
\]
where $\modE$ is the modulation operator given by 
\[
\modE f(x) \defeq (-1)^{\lev(x)} f(x), \qquad \forall x \in \Tree.
\]
From the estimate above we then deduce that $F_0(\opL)$ and $F_2(2-\opL)$ are of weak type $(1,1)$, so the same is true for $\modE F_2(2-\opL) \modE$, and therefore for $F(\opL)$ too.
\end{proof}

\bibliographystyle{abbrvalpha}
\bibliography{bibliography}

\end{document}